\documentclass[times, doublespace]{nmeauth}

\usepackage{graphicx}
\usepackage{nicefrac}

\usepackage{amssymb}

\usepackage{amsmath}
\usepackage[hang,raggedright]{subfigure}

\usepackage{lineno}

\usepackage{xspace}

\usepackage[square, numbers, comma, sort&compress]{natbib}

\usepackage[squaren]{SIunits}

\usepackage{binary}
\usepackage{cite}
\usepackage{verbatim}  
\usepackage{color}

\usepackage{listings}
\usepackage{url}

\lstdefinelanguage{C++enhanced}[]{C++}{%
  morekeywords={volDiagTensorField, volSymmTensor4thOrderField, fvVectorMatrix,dimensionedScalar},
}
\definecolor{gray}{rgb}{0.5,0.5,0.5}
\usepackage{courier}
\usepackage{algorithm}
\usepackage[noend]{algpseudocode}

\newcommand{\etc} {\textit{etc.}\@\xspace}

\newcommand{\ie} {\textit{i.e.}\@\xspace}
\newcommand{\foam} {{OpenFOAM}\xspace}
\newcommand{\etal} {et~al.\xspace}
\renewcommand*{\refname}{}

\newcommand{\ra}[1]{\renewcommand{\arraystretch}{#1}}

\begin{document}


\runningheads{P. Cardiff \etal}{A Finite Volume Solid Mechanics Toolbox for OpenFOAM}

\title{An open-source finite volume toolbox for solid mechanics and fluid-solid interaction simulations}




\author
{
P. Cardiff\affil{1}\corrauth,
A. Kara\v{c}\affil{2},
P. De Jaeger\affil{1,3},
H. Jasak\affil{4},
J. Nagy\affil{5},
A. Ivankovi\'{c}\affil{1}
and
\v{Z}. Tukovi\'{c}\affil{4}
}

\address{
\affilnum{1}University College Dublin, Bekaert University Technology Centre, School of Mechanical and Materials Engineering, Belfield, Dublin, Ireland \break
\affilnum{2}University of Zenica, Polytechnic Faculty, Fakultetska 1, Zenica, Bosnia and Herzegovina \break
\affilnum{3}NV Bekaert SA, Belgium \break
\affilnum{4}University of Zagreb, Faculty of Mechanical Engineering and Naval Architecture, Croatia \break
\affilnum{5}Johannes Kepler University, Linz, Austria
}

\corraddr{Philip Cardiff, University College Dublin, Bekaert University Technology Centre, School of Mechanical and Materials Engineering, Belfield, Ireland. E-mail: philip.cardiff@ucd.ie}

\begin{abstract}
\noindent
Over the past 30 years, the cell-centred finite volume method has developed to become a viable alternative to the finite element method in the field of computational solid mechanics.
The current article presents an open-source toolbox for solid mechanics and fluid solid interaction simulations based on the finite volume library OpenFOAM$^{\tiny{\textregistered}}$.
The object oriented toolbox design is outlined, where emphasis has been given to code use, comprehension, maintenance and extension.
The toolbox capabilities are demonstrated on a number of representative test problems, where comparisons are given with finite element solutions.

\end{abstract}



\maketitle




\section{Introduction}


Numerical analysis techniques, such as the finite element (FE) and finite volume (FV) methods, have found widespread use in modern engineering industry and academia, ranging from component design to the analysis of physical mechanisms.
For multi-physics problems, there are a number of viable commercial software options; however, there has been relatively little development of open-source software for this purpose.
In particular, there are few open-source packages that offer significant fluid and solid analysis capabilities within the same framework.

The FE method dominates the field of computational solid mechanics (CSM), whereas the FV method is the most popular technique for computational fluid dynamics (CFD); this poses a significant challenge for problems which involve both fluid-like and solid-like materials. Coupling of separate analysis packages and techniques is an option; however, performing such multi-physics analyses within a consistent software framework offers a number of advantages in terms of code development and solver efficiency.

In the past decade, open-source software OpenFOAM has become one of the most popular packages in the realm of CFD and multi-physics simulations.
The OpenFOAM project emerged from the Field Operations And Manipulations (FOAM) concept at Imperial College London in the late-1980s.
As stated by \citet{Weller1998}, an aim of the project was ``\emph{to make it as easy as possible to develop reliable and efficient computational continuum-mechanics codes}"; this was achieved by exploiting newly developed object-oriented design paradigms, and making the top-level syntax closely resemble conventional mathematical notation for partial differential equations;
combined with native parallelisation support and the open-source GPL license \citep{GNUGPL}, this has lead to the widespread adoption of OpenFOAM in academia and industry, for example, \citep{Jakirlic2017, Islam2009}.

Although OpenFOAM has seen significant developments across the entire domain of CFD applications, the establishment of procedures for CSM and fluid-solid interactions (FSI) is still nascent, contributed to by the general lack of personnel with combined CFD and CSM expertise. In the formative publications on OpenFOAM, Weller \etal \citep{Weller1998, Jasak2000:linearElasticity} presented analysis of the classical \emph{hole in an elastic plate} problem, where the cell-centred FV discretisation built on the pioneering work of Demird\v{z}i\'{c} and co-workers \citep{Demirdzic1988, Demirdzic1995, Demirdzic1997}.
Subsequently, within the OpenFOAM framework, there have been a number of extensions to the cell-centred FV methods for CSM, as regards finite strains \citep{Tukovic2007:updatedLagrangainFAMENA, Cardiff2014:orthotropicPaper, Cardiff2016:metalForming, Leonard2012}, constitutive relations \citep{Tang2015, Cardiff2014:orthotropicPaper, Cardiff2016:metalForming, Safari2016}, boundary conditions \citep{Cardiff2012:contactPaper, Tukovic2010acp, Carolan2013, Cardiff2014:mappedMusclesPaper, Cardiff2014:hipPaper}, discretisation and solution methodologies \citep{Tukovic2012, Cardiff2016:blockCoupled, Elsafti2016, Haider2017}, as well as FSI approaches, notably \citep{Greenshields1999, Ivankovic2002, Karac2002, Tukovic2007:fsi, Karac2009, Degroote2009, Tukovic2014, Cardiff2015:hydraulicFractures, Gillebaart2016, Sekutkovski2016, Tukovic2017:fsi}.

Apart from \citet{Tukovic2014, Tukovic2017:fsi}, the previous works all focused on the development of bespoke numerical procedures, where no emphasis was placed on the creation of generalised code structures that may be easily combined, adapted and extended to related CSM and FSI problems.
In particular, little attention has been given to the development of modular designs, allowing straight-forward combination of differing solid and fluid procedures for FSI problems, or allowing for the addition of new solid constitutive relations, for example.
The official versions of OpenFOAM are released with only the most basic solid mechanics tools, specifically the \texttt{solidDisplacementFoam} solver, which is restricted to Hookean solids undergoing small strains and rotations.
The current article aims to overcome these shortcomings by presenting and sharing a solid mechanics and FSI toolbox built on OpenFOAM software framework.
The adopted approach aims to follow object-oriented design principles that allow for straight-forward use, comprehension, extension, and maintenance of the code.
Furthermore, the toolbox will be made freely available to the community via the FOAM Extend community fork \citep{OpenFOAMExtend} and the OpenFOAM Community Repository \citep{OpenFOAMCommunityRepo}, such that the implementations are fully open to academic scrutiny and allow future collaborative improvements.

The article is structured as follows: Section 2 gives an overview of the toolbox design; Section 3 summarises the mathematical models that are solved; Section 4 outlines the cell-centred FV discretisation and solution methodologies; in Section 5, a number of test cases are presented, highlighting the capabilities and applicability of the toolbox; finally, Section 6 briefly discusses the presented methods in light of conventional FE methods and indicates future directions.

\section{Toolbox Structure \& Design}

Over the past half century, there has been extensive literature published on the development of discretisations, solution methodologies and algorithms within the field of computational mechanics; however, relatively little importance has been given to code design: this is somewhat understandable given the procedures are typically developed by engineers and applied mathematicians as opposed to computer programmers.
As the complexity of modern parallelised multi-physics procedures continues to increase, good code design becomes ever more crucial to the solution of complex problems.
Historically, procedural programming implementations were favoured, for example, in FOTRAN, but due to the proliferation of object-oriented programming (OOP) languages, such as \texttt{C++} in the 1990s, there has been a gradual turn to the benefits of such abstract programming paradigms.
One of the key characteristics of procedural programming is that it relies on \emph{procedures} that operate on \emph{data}; in contrast, these two separate concepts are bundled into objects in OOP approaches: this makes it possible to create complex behaviour with less code and less repetition of code.

With the development of OpenFOAM, \citet{Weller1998} was one of the first to place specific emphasis on the design philosophy of an engineering analysis tool for CFD;
the philosophy of OOP was embraced, where conceptual constructs (a \emph{class}) are represented in the program, with specific details hidden behind an interface with limited, well-defined access;
in this way, data types are created that represent, for example, tensor fields allowing the governing partial differential equations to be constructed in a manner similar to their mathematical counterparts.
Numerical details - irrelevant at this level - are hidden by \emph{encapsulation} (contain and protect the data that make up the class).
There are a number of other notable works that have emphasised the software design element of an engineering tool, for example, \citet{Archer1996}, \citet{iMOOSE2004}, \citet{Bangerth2007}, and \citet{SU2:2013}. 


As OpenFOAM is released under the open-source GPL license \citep{GNUGPL}, it lends itself to the development of add-on toolboxes, such as a wave generation toolbox by \citet{Jacobsen2012}, a toolbox for modelling soil-structure interaction around marine structures by \citet{Elsafti2016}, and a porous media flow toolbox by \citet{Horgue2015}; in a similar fashion, the current article presents here a toolbox, titled \texttt{solids4foam}, which aims to generalise the OpenFOAM design further to allow straightforward implementation of advanced solid mechanics and fluid-solid interaction procedures. Throughout the design the \texttt{solids4foam}, significant emphasis has been placed on the following complementary factors:
\begin{itemize}
	\item \emph{usable}: the solvers should be intuitive for a new user to use;
	\item \emph{understandable}: the code structure should be easy to follow for a developer who is not familiar with the project, where esoteric structures are avoided; for example, class, data and function names should be descriptive, and a consistent style convention should be adopted;
	\item \emph{maintainable}: the design should enable straightforward maintenance, in particular as the code grows and features are added; modularisation should be embraced and code repetition should be avoided;
	\item \emph{extendable}: the design should allow key features to be extended in a straight-forward manner, particularly for unfamiliar developers, for example, adding a new constitutive relation, solution methodology or boundary condition.
\end{itemize}

In the subsequent sections, these four points are addressed through the design and structure of the presented \texttt{solids4foam} toolbox.

\subsection{Design of the class structure}


In the standard OpenFOAM public release, there are a large number of \emph{solver} executables, each designed for a specific group of problems, for example, \texttt{laplacianFoam} solves the transient Laplace equation (heat conduction), \texttt{icoFoam} solves the isothermal laminar Newtonian incompressible Navier-Stokes equations and \texttt{solidDisplacementFoam} solves the small strains/rotations Hookean elastic equations; this approach has a number of merits relative to the single executable design approach; however, it has its limitations when considering coupled physics problems.
Consider a fluid-solid interaction problem where we wish to combine any one of the fluid models with any one of the solid models, for example, a compressible turbulent multiphase fluid with a visco-elastic neo-Hookean solid.
To overcome such a limitation, the current article builds on the framework of Tukovi\'{c} \etal \citep{Tukovic2014, Tukovic2017:fsi} to propose a modular design structure where solution methodologies are encapsulated within classes; for example, the Pressure-Implicit with Splitting of Operators (PISO) solution algorithms for fluids, or the total Lagrangian nonlinear geometry solution algorithm for solids.

\begin{figure}[htb]
	\centering
	\includegraphics[width=\textwidth]{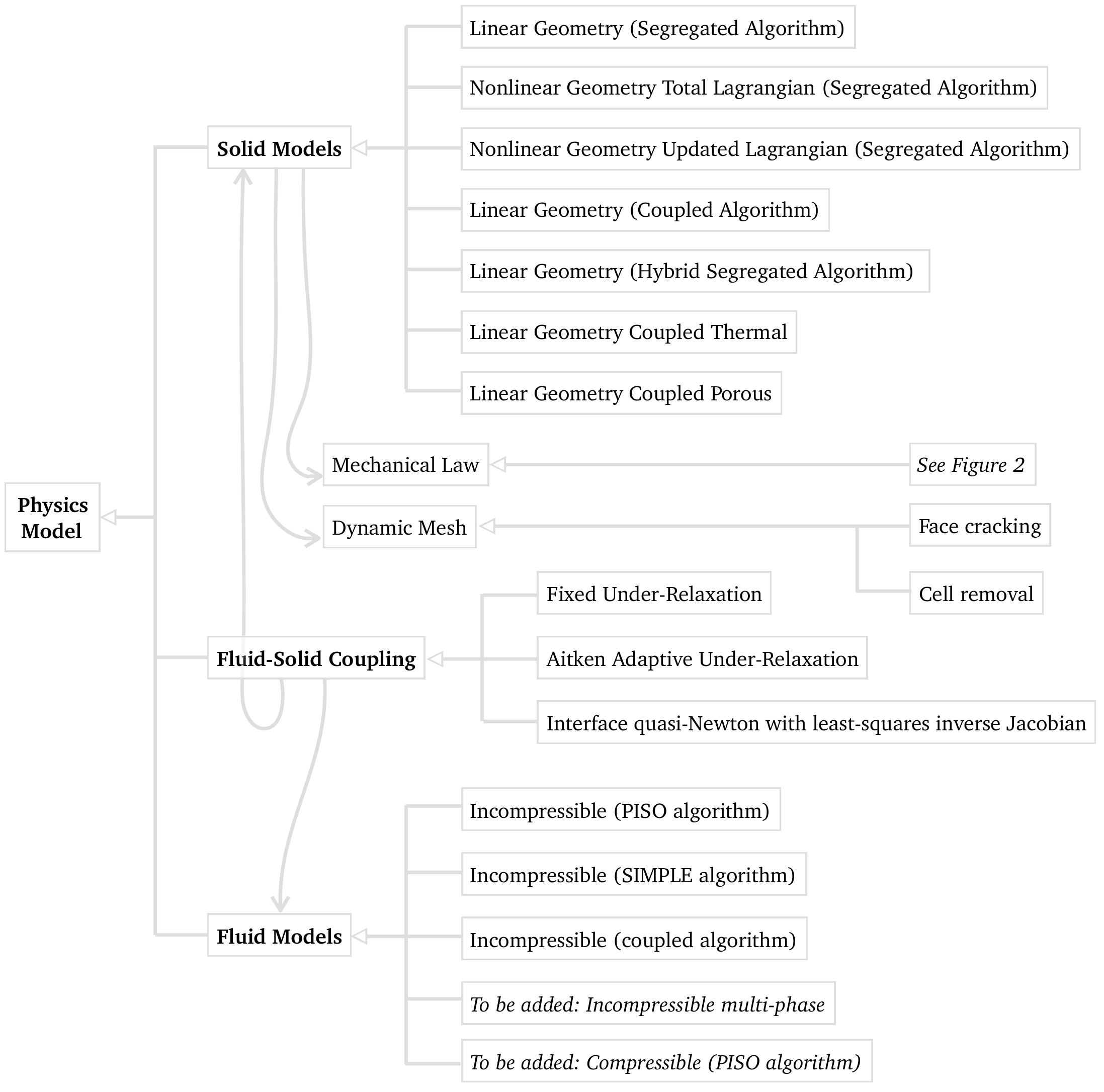}
	\caption{Class structure: fluid, solid, and fluid-solid interaction mathematical models}
	\label{fig:classStructure1}
\end{figure}

The top level class structure of \texttt{solids4foam} is shown in Figure \ref{fig:classStructure1}, where the UML notation \citep{Larman1998} is broadly followed;
as indicated, the \texttt{Physics Model} class represents a mathematical model (boundary value problem), for either a solid or fluid or fluid-solid problem.
At run-time the user specifies which specific implementation will be invoked, via user input files; in this way, the \texttt{Physics Model} is a so-called abstract base class, from which \texttt{Solid Models}, \texttt{Fluid Models} and \texttt{Fluid-Solid Interaction Models} are derived.
This demonstrates the OOP concept of \emph{polymorphism}, where classes can be considered conceptually equivalent and share the same interface but have different implementations;
for example, at run-time the user can switch between a mathematical model for a solid that assumes small strains (linear geometry) and a mathematical model valid for finite strains (nonlinear geometry), or switch between a segregated and block-coupled solution methodology.

\begin{figure}[htb]
	\centering
	\includegraphics[width=0.9\textwidth]{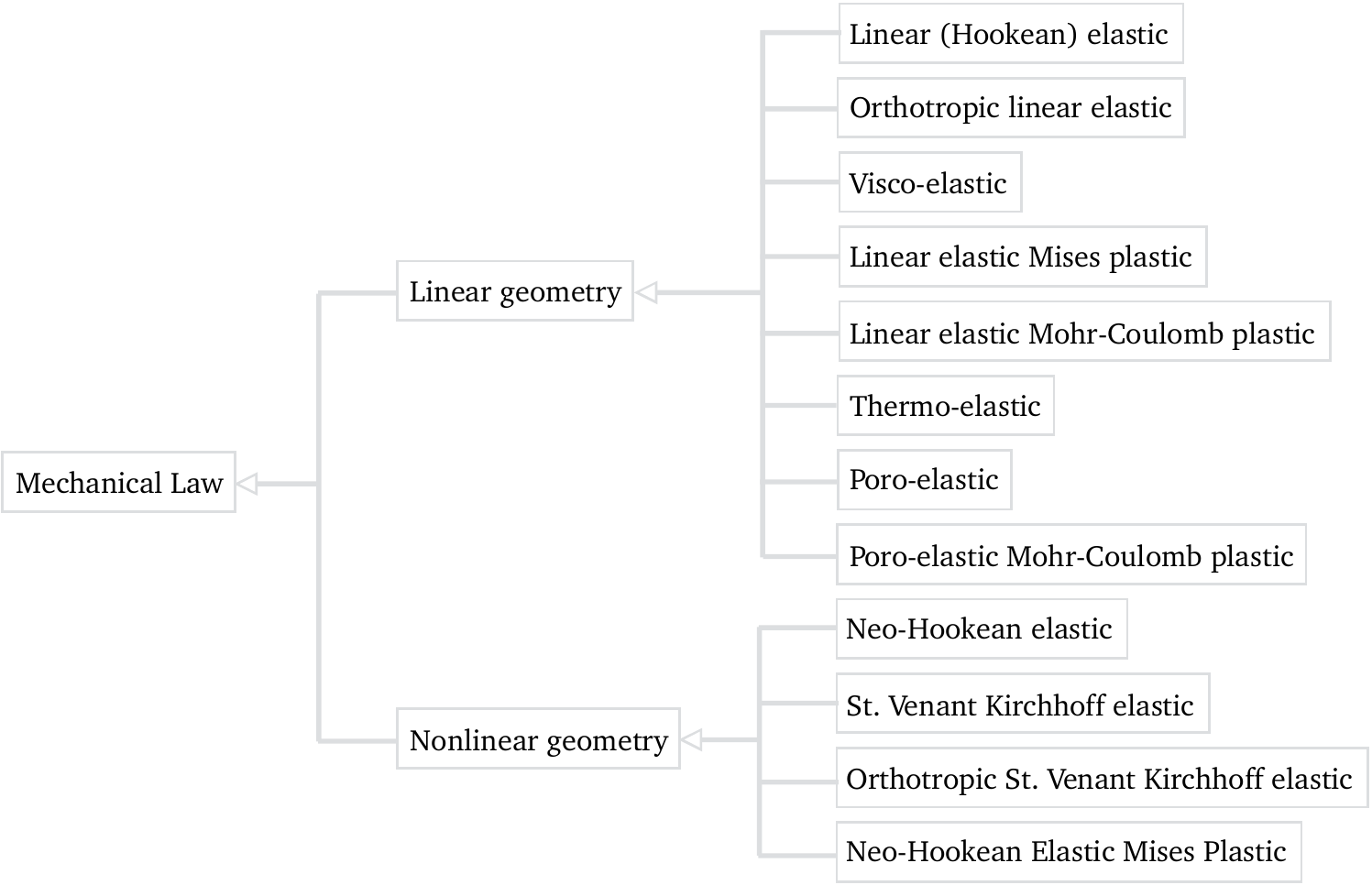}
	\caption{Class structure: mechanical laws and dynamic meshes }
	\label{fig:classStructure2}
\end{figure}

In Figure \ref{fig:classStructure1}, the lines with a closed, unfilled arrowhead indicate an \emph{inheritance} relationship between classes: for example, the \texttt{Linear Geometry (Segregated Algorithm)} class inherits the attributes of the \texttt{Solid Model} class, and similarly the \texttt{Solid Model} class inherits the attributes of the \texttt{Physics Model} class; in this way, inheritance enables commonality relationships to be expressed between classes.
It should be noted that the polymorphism and inheritance mechanisms are effective methods to minimise code duplication (\ie repeating segments of code in multiple places), which is a critical point to ensure the code maintainability. In addition, such a system allows straightforward addition of new derived classes, for example, consider the mechanical constitutive law classes in Figure \ref{fig:classStructure2}.
Lines ending in open arrowheads (for example, from the \texttt{Mechanical Law} to the \texttt{Solid Model}) indicates a uni-directional association; in such an association, the two classes are related, but only one class knows that the relationship exists; for example, the \texttt{Solid Model} class contains a \texttt{Mechanical Law}, which is required to calculate the stress tensor field; however, the \texttt{Mechanical Law} class need not know about how the specific solid model solution algorithm is implemented and should only be concerned with the calculation of the stress field when asked.
This demonstrates the concept of \emph{encapsulation} where a \texttt{Mechanical Law} derived class, such as \texttt{Linear (Hookean) elastic}, encapsulates and protects data specific to the law (for example, Young's modulus and Poisson's ratio) and only provides access through a well defined interface (for example, calculate the stress given the strain or deformation gradient field).

The class structure and design presented in Figures \ref{fig:classStructure1} and \ref{fig:classStructure2} takes advantage of design patterns commonly used in OOP; further discussion of design patterns can be found in \citep{Gamma1995}. The use of such design patterns can speed up the development process by providing tested, proven development paradigms; in addition, the reuse of design patterns helps improve code readability for those familiar with the patterns. In this case the so-called \emph{factory method} \citep{Gamma1995} is exploited multiple times, where a parent class defines an interface for creating objects, but lets subclasses decide which classes to instantiate; for example, the \texttt{Solid Model} parent class defines the interfaces but derived classes such as \texttt{Linear Geometry (Segregated Algorithm)} define the actual implementation.
Considerable time has been given to consider future extensions that may occur and ensuring the code design is sufficiently flexible to allow for this; however, the code base is not static and continues to evolve; as such, the habit of refactoring continuously will allow easier extension and maintenance; in this sense, refactoring refers to the process of restructuring existing computer code without changing its external behaviour, in an attempt to improve code readability and reduce complexity.


\subsection{Note on Style, Revision Control \& Documentation}
As the \texttt{solids4foam} toolbox is presented to the community as an open-source repository, users and developers from other institutions are invited to contribute; consequently, to ensure that the toolbox does not become a \emph{patchwork quilt} of contrasting styles making the code difficult to read, it is critical to enforce stringent coding standards;
strict enforcement of such coding standards ensures the software is easy to read, maintain and extend.
Following such standards immediately ensures code is more readable and understandable, allowing faster and more error-free developments; in addition, when coding standards are not followed, reading code generated by others becomes a tedious, painstaking and even an impossible task.
Fortunately, the developers of the public release of OpenFOAM have set out a comprehensive coding style \citep{OpenFOAM:codingStyle}; this style is strictly followed in the \texttt{solids4foam} toolbox.

A Git revision control system is used to track changes to the source code, both on the local repository as well as within the open-source published source code. 
Doxygen code documentation provides developers with the code structure and interfaces.

\subsection{Note on the open-source paradigm}
Commercial engineering software offers notable capabilities as regards the implementation of custom user routines; however, the development and implementation of novel solution procedures is not possible without access to the code base.
For those interested in such developments, whether for academic or industrial reasons, one option is to develop code \emph{in-house}.
This overcomes the limitations of commercial software; however, it has the drawback that the documentation, maintenance and extension of the code can be difficult due to the small \emph{in-house} user base.
A potentially more appealing option is starting from a suitable open-source library with a significant user base;
this potentially increases the speed of developing new procedures as well as allowing the reuse of code provided by the community.
An additional benefit of the open-source paradigm is that the code base is fully open to academic scrutiny, where normal users are granted a complete view into the source code; this provides access to which algorithms are implemented, and more importantly \emph{how} they are implemented.
Within such an open-source framework, large scale community collaborations become possible, and development of complex multi-physics procedures can proceed at an accelerated rate.

Open-source software is, of course, not without its weaknesses: it requires a certain level of system administration experience and technical expertise in order to manage and develop content; in addition, there are often no official code reviews or quality assurance processes in place, instead relying on users and developers to perform these verifications, checks and subsequent fixes.
Nevertheless, open-source software projects certainly possess a number of advantageous features suitable for the development of state-of-the-art numerical procedures.

\section{Mathematical Models}
In this section, an overview is given of the mathematical models that are solved in the \texttt{Solid Model}, \texttt{Fluid Model} and \texttt{Fluid-Solid Interaction} classes.
As the fluid models and FSI interaction procedures have been extensively described previously, for example, \citep{Weller1998, Tukovic2017:fsi}, emphasis is given here to the description of the solid model methods.

\subsection{Governing Equations}
Considering an arbitrary body of volume $\Omega$ bounded by surface $\Gamma$ with outward facing normal $\boldsymbol{n}$, the governing equations \-- conservation of mass, conservation of energy (heat equation form) and conservation of linear momentum \-- are given in strong integral form as:
\begin{eqnarray}
	\frac{\text{D}}{\text{D} t}  \int_\Omega \rho  \; \text{d}\Omega = 0
\end{eqnarray}
\begin{eqnarray}
	\frac{\text{D}}{\text{D} t}  \int_\Omega \rho C_p T  \; \text{d}\Omega
	\;=\; - \oint_\Gamma \boldsymbol{n} \cdot  \boldsymbol{q} \; \text{d}\Gamma
	\;+\; \int_\Omega \boldsymbol{\sigma} : \boldsymbol{\nabla} \boldsymbol{v} \; \text{d}\Omega
\end{eqnarray}
\begin{eqnarray}
	\frac{\text{D}}{\text{D} t}  \int_\Omega \rho  \boldsymbol{v} \; \text{d}\Omega
	\;=\; \oint_\Gamma \boldsymbol{n} \cdot  \boldsymbol{\sigma} \; \text{d}\Gamma
	\;+\; \int_\Omega \rho \,\boldsymbol{b} \; \text{d}\Omega
\end{eqnarray}
where $\rho > 0$ is the density, $C_p > 0$ is the specific heat capacity at constant pressure, $T$ is the temperature, $\boldsymbol{q}$ is the heat flux vector, $\boldsymbol{\sigma}$ is the Cauchy stress tensor, $\boldsymbol{v}$ is the velocity vector, and $\boldsymbol{b}$ is the body force per unit mass; symbol $\boldsymbol{\nabla}$ signifies the so called Hamilton operator, synonymous with the del or nabla operator.
In addition, the necessary and sufficient condition for the conservation of angular momentum is the symmetry of the Cauchy stress tensor.

Employing Reynold's transport theorem, the total derivatives may be replaced by partial derivatives as:
\begin{eqnarray} \label{eq:massEqn}
	\frac{\partial}{\partial t} \int_\Omega \rho  \; \text{d}\Omega
	\;+\; \oint_\Gamma \rho \left[ \boldsymbol{n} \cdot (\boldsymbol{v} - \boldsymbol{v}_\Gamma) \right]  \; \text{d}\Gamma
	= 0
\end{eqnarray}
\begin{eqnarray} \label{eq:energyEqn}
	\frac{\partial}{\partial t} \int_\Omega \rho  C T  \; \text{d}\Omega
	\;+\; \oint_\Gamma \rho C T \left[\boldsymbol{n} \cdot (\boldsymbol{v} - \boldsymbol{v}_\Gamma) \right]   \; \text{d}\Gamma
	\;=\; - \oint_\Gamma \boldsymbol{n} \cdot  \boldsymbol{q} \; \text{d}\Gamma
	\;+\; \int_\Omega \boldsymbol{\sigma} : \boldsymbol{\nabla} \boldsymbol{v} \; \text{d}\Omega
\end{eqnarray}
\begin{eqnarray} \label{eq:momentumEqn}
	\frac{\partial}{\partial t} \int_\Omega \rho  \boldsymbol{v}  \; \text{d}\Omega
	\;+\; \oint_\Gamma \rho \boldsymbol{v} \left[\boldsymbol{n} \cdot (\boldsymbol{v} - \boldsymbol{v}_\Gamma) \right]   \; \text{d}\Gamma
	\;=\; \oint_\Gamma \boldsymbol{n} \cdot  \boldsymbol{\sigma} \; \text{d}\Gamma
	\;+\; \int_\Omega \rho \,\boldsymbol{b} \; \text{d}\Omega
\end{eqnarray}
where $\boldsymbol{v}_\Gamma$ is the velocity of the control volume surface $\Gamma$.
As yet, no distinction has been made between solid and fluid materials; the difference comes in the assumed constitutive relations, for example, in the definition of the Cauchy stress $\boldsymbol{\sigma}$. In addition, for convenience a Lagrangian approach is typically adopted for the analysis of solids, with an Eulerian approach adopted for fluids; in the Lagrangian approach the velocity of the material is assumed equal to the velocity of the domain \ie $\boldsymbol{v} - \boldsymbol{v}_\Gamma = \boldsymbol{0}$; consequently, the so-called \emph{convection} terms (second term on the left-hand side of Equations \ref{eq:massEqn}, \ref{eq:energyEqn} and \ref{eq:momentumEqn}) drop out in a Lagrangian formulation, but remain in an Eulerian approach.

\subsection{Constitutive Relations}
To close the governing equations, constitutive relations defining the Cauchy stress $\boldsymbol{\sigma}$ and the heat flux $\boldsymbol{q}$ need to be specified.

\paragraph{Fluid constitutive laws}
One of the main purposes of the \texttt{solids4foam} toolbox is to allow use of multiple different fluid and solid constitutive relations and solution methodologies; the multitude of fluid models available in the public release of OpenFOAM have been described elsewhere, for example, \citep{Weller1998, Jasak1996}; consequently, for illustrative purposes, only the case of a laminar incompressible isothermal Newtonian fluid is described here, where the Cauchy stress is given as:
\begin{eqnarray}
	\boldsymbol{\sigma} = 2 \eta \boldsymbol{D} - p \textbf{I}
\end{eqnarray}
The dynamic viscosity is $\eta > 0$, $p$ is the pressure field, and $\textbf{I}$ is the second order identity tensor.
The rate of deformation $\boldsymbol{D}$ is given as the symmetric component of the velocity gradient:
\begin{eqnarray}
	\boldsymbol{D}
	&=& \text{symm}[\boldsymbol{\nabla} \boldsymbol{v}] \notag \\
	&=& \frac{1}{2} \left(\boldsymbol{\nabla} \boldsymbol{v}
	+ \boldsymbol{\nabla} \boldsymbol{v}^T \right)
\end{eqnarray}
where the operator $\text{symm}[\cdot]$ refers to the symmetric component of a tensor.

\paragraph{Solid constitutive laws: \emph{linear geometry}}
In the case of small strains and small rotations, the linearised strain tensor is defined in terms of the displacement gradient:
\begin{eqnarray}
	\boldsymbol{\epsilon}
	&=& \text{symm}[\boldsymbol{\nabla} \boldsymbol{u}] \notag \\
	&=& \frac{1}{2} \left(\boldsymbol{\nabla} \boldsymbol{u}
	+ \boldsymbol{\nabla} \boldsymbol{u}^T \right)
\end{eqnarray}
The deviatoric component of the strain given as:
\begin{eqnarray}
	\boldsymbol{e} = \text{dev}[\boldsymbol{\epsilon}] = \boldsymbol{\epsilon} - \tfrac{1}{3} \text{tr} [\boldsymbol{\epsilon}] \textbf{I}
\end{eqnarray}
where the deviatoric operator is indicated by $\text{dev}[\cdot]$ and the trace operator is indicated by $\text{tr}[\cdot]$.

For the \emph{linear geometry} case, where small strains and small rotations are assumed, there are a large number of well-known constitutive laws; for illustrative purposes, a selection of laws, implemented within the \texttt{solids4foam} toolbox, are briefly described here.
As stated previously, it is intended that implementation of new laws should be a straight-forward process.
The definition of the engineering stress $\boldsymbol{\sigma}_s$ for a number of popular solid constitutive laws (\texttt{Mechanical Law} in the class structure) that are suitable for linear geometry (small strains and rotations) are given in Table \ref{table:MechanicalLawLinearGeometry};
the corresponding mechanical parameters are described in Appendix \ref{App:A}, Table \ref{table:MechanicalLawLinearGeometryParameters}.
\begin{table}[htb]
  \centering
	\ra{1.3}
	\begin{tabular}{@{}lll@{}}
	\toprule
	Linear (Hookean) elastic	&
	$\boldsymbol{\sigma}_s$
	& $=2\mu \boldsymbol{\epsilon} + \lambda \text{tr}[\boldsymbol{\epsilon}] \textbf{I}$ \\
	& &
	$=2\mu \boldsymbol{e} + \kappa \text{tr}[\boldsymbol{\epsilon}] \textbf{I}$
	\\
	Orthotropic linear elastic	&
	$\boldsymbol{\sigma}_s$ &$= \boldsymbol{C}_e : \boldsymbol{\epsilon} $
	\\
	Linear visco-elastic	&
	$\boldsymbol{\sigma}_s(t)$ &$=
	\int_{-\infty}^{t} g(t - s)
	\; 2 \mu \boldsymbol{e} \; \text{d}s
	+ \kappa \, \text{tr}[\boldsymbol{\epsilon}] \textbf{I}$
	\\
	& & \quad $g(t) = \gamma_{\infty} + \sum_{i=1}^N \gamma_i \exp[-\nicefrac{t}{\tau_i}]$
	\\
	Thermo-linear elastic &
	$\boldsymbol{\sigma}_s$
	&
	$=2\mu \boldsymbol{\epsilon} + \lambda \text{tr}[\boldsymbol{\epsilon}] \textbf{I}
	- (2\mu + 3\lambda) \alpha (T - T_0) \textbf{I} $ \\
	& &
	$=2\mu \boldsymbol{e}
		+ \kappa \text{tr}[\boldsymbol{\epsilon}] \textbf{I}
		- 3 \kappa \alpha (T - T_0) \textbf{I}$
	\\
	Poro-linear elastic  &
	$\boldsymbol{\sigma}_s$
	& $=2\mu \boldsymbol{\epsilon} + \lambda \text{tr}[\boldsymbol{\epsilon}] \textbf{I}
	- p \textbf{I} $ \\
	& &
	$=2\mu \boldsymbol{e} + \kappa \text{tr}[\boldsymbol{\epsilon}] \textbf{I}
	- p \textbf{I} $
	\\
	Linear elastic, Mises/J$_2$ plastic &
	$\boldsymbol{\sigma}_s$
	& $=2\mu \boldsymbol{\epsilon}_e + \lambda \text{tr}[\boldsymbol{\epsilon}_e] \textbf{I}$ \\
	& & \quad $\boldsymbol{\epsilon} = \boldsymbol{\epsilon}_e + \boldsymbol{\epsilon}_p$
	, \quad $\dot{\boldsymbol{\epsilon}}_p = \dot{\Lambda} \,
	(\nicefrac{\text{dev}[\boldsymbol{\sigma}_s]}{||\text{dev}[\boldsymbol{\sigma}_s]||})$  \\
	& & \quad $f(\boldsymbol{\sigma}_s, \epsilon^{eq}_p) = ||\text{dev}[\boldsymbol{\sigma}_s]||
	- \sqrt{\frac{2}{3}}\sigma_Y(\epsilon^{eq}_p)$
	\\
	& & \quad  $\Delta \Lambda \geq 0,
	\quad f(\boldsymbol{\sigma}_s, \epsilon^{eq}_p) \leq 0,
	\quad \Delta \Lambda f(\boldsymbol{\sigma}_s, \epsilon^{eq}_p) = 0$
	\\
	Linear elastic, Mohr-Coulomb plastic &
	$\boldsymbol{\sigma}_s$
	& $=2\mu \boldsymbol{\epsilon}_e + \lambda \text{tr}[\boldsymbol{\epsilon}_e] \textbf{I}$ \\
	& & \quad $\boldsymbol{\epsilon} = \boldsymbol{\epsilon}_e + \boldsymbol{\epsilon}_p$
	, \quad $\dot{\boldsymbol{\epsilon}}_p = \dot{\Lambda} \,
	(\nicefrac{\partial g}{\partial \boldsymbol{\sigma}_s})$ \\
	& &	\quad $f(\boldsymbol{\sigma}_s)
	= \sigma_1 - \sigma_3 + (\sigma_1 + \sigma_3) \sin{\phi} - 2 c \cos{\phi}$
	\\
	& &	\quad $g(\boldsymbol{\sigma}_s)
	= \sigma_1 - \sigma_3 + (\sigma_1 + \sigma_3) \sin{\psi} $
	\\
	& & \quad  $\Delta \Lambda \geq 0,
	\quad f(\boldsymbol{\sigma}_s) \leq 0,
	\quad \Delta \Lambda f(\boldsymbol{\sigma}_s) = 0$	\\
	\bottomrule
	\end{tabular}
\caption{Solid constitutive laws: linear geometry}
\label{table:MechanicalLawLinearGeometry}
\end{table}

\paragraph{Solid constitutive laws: \emph{nonlinear geometry}}
When considering finite strains, there are once again a large number of popular constitutive relations to define the Cauchy (true) stress.
Table \ref{table:MechanicalLawNonlinearGeometry} summarises a number of finite strain laws implemented within the \texttt{solids4foam} toolbox.
It should be noted that in the limit of small strains and rotations, these nonlinear relations reduce to their linearised counterparts, for example, neo-Hookean elastic reduces to linear (Hookean) elastic.

It is useful to define a number of additional kinematic quantities when dealing with finite strains:
the deformation gradient, $\boldsymbol{F}$, and relative deformation gradient, $\boldsymbol{f}$, may be defined with respect to the initial configuration:
\begin{eqnarray}
	&\boldsymbol{F} = \textbf{I} + (\boldsymbol{\nabla}_0 \boldsymbol{u})^T \\
	&\boldsymbol{f} = \boldsymbol{F} \cdot \boldsymbol{F}_{[m-1]}^{-1}
\end{eqnarray}
or equivalently with respect to the updated configuration as:
\begin{eqnarray} \label{eq:relativeDeformationGradient}
	&\boldsymbol{F} = \boldsymbol{f} \cdot \boldsymbol{F}_{[m-1]} \\
	&\boldsymbol{f} = \textbf{I} + (\boldsymbol{\nabla}_u \Delta \boldsymbol{u})^T
\end{eqnarray}
where $\boldsymbol{\nabla}_0$ represents the gradient with respect to the initial undeformed configuration, and $\boldsymbol{\nabla}_u$ represents the gradient with respect to the so-called updated configuration (configuration at the end of the previous time step);
$\boldsymbol{F}_{[m-1]}$ is the deformation gradient at the end of the previous time step;
the increment of displacement is defined as $\Delta \boldsymbol{u} = \boldsymbol{u} - \boldsymbol{u}_{[m-1]}$, with $\boldsymbol{u}_{[m-1]}$ being the displacement at the end of the previous time step.
\begin{table}[htb]
  \centering
	\ra{1.3}
	\begin{tabular}{@{}lll@{}}
	\toprule
	Neo-Hookean elastic &
	$\boldsymbol{\sigma}$
	& $= \mu \, \text{dev}[\bar{\boldsymbol{b}}] + \frac{\kappa}{2} (\frac{J^2 - 1}{J}) \textbf{I}$ \\
	& & \quad $\bar{\boldsymbol{b}} = J^{-\nicefrac{2}{3}} \boldsymbol{F} \cdot \boldsymbol{F}^T$ \\
	St. Venant Kirchhoff elastic &
	$\boldsymbol{\sigma}$
	& $= J^{-1} \boldsymbol{F} \cdot \boldsymbol{S} \cdot \boldsymbol{F}^T $ \\
	& & \quad $\boldsymbol{S} = 2\mu \boldsymbol{E}
	+ \lambda \text{tr}[\boldsymbol{E}] \textbf{I}$ \\
	& & \quad $\boldsymbol{E} = \frac{1}{2}(\boldsymbol{F}^T \cdot \boldsymbol{F} - \textbf{I})$ \\
	Orthotropic St. Venant Kirchhoff elastic &
	$\boldsymbol{\sigma}$
	& $= J^{-1} \boldsymbol{F} \cdot \boldsymbol{S} \cdot \boldsymbol{F}^T $ \\
	& & \quad $\boldsymbol{S} = \boldsymbol{C}_e : \boldsymbol{E}$ \\
	Neo-Hookean elastic, Mises/J$_2$ plastic &
	$\boldsymbol{\sigma}$
	& $= \mu \, \text{dev}[\bar{\boldsymbol{b}}_e] + \frac{\kappa}{2} (\frac{J^2 - 1}{J}) \textbf{I}$ \\
	& & \quad $\boldsymbol{F} = \boldsymbol{F}_e \cdot \boldsymbol{F}_p$,
	\quad $\boldsymbol{N} = \nicefrac{\text{dev}[\boldsymbol{\sigma}]}{||\text{dev}[\boldsymbol{\sigma}]||}$ \\
	& & \quad $\frac{\partial}{\partial t}(\boldsymbol{F}^{-1} \cdot \boldsymbol{b}_e \cdot \boldsymbol{F}^{-T})
	= \frac{2}{3} \dot{\Lambda} \, \text{tr}[\boldsymbol{b}_e]
	\boldsymbol{F}^{-1} \cdot \boldsymbol{N}\cdot \boldsymbol{F}^{-T}$  \\
	& & \quad $f(\boldsymbol{\sigma}, \epsilon^{eq}_p) = ||\text{dev}[\boldsymbol{\sigma}]||
	- \sqrt{\frac{2}{3}}\sigma_Y(\epsilon^{eq}_p)$
	\\
	& & \quad  $\Delta \Lambda \geq 0,
	\quad f(\boldsymbol{\sigma}, \epsilon^{eq}_p) \leq 0,
	\quad \Delta \Lambda f(\boldsymbol{\sigma}, \epsilon^{eq}_p) = 0$
	\\
	\bottomrule
	\end{tabular}
\caption{Solid constitutive laws: nonlinear geometry}
\label{table:MechanicalLawNonlinearGeometry}
\end{table}

\paragraph{Thermal constitutive law}
For both fluid and solid materials, the heat flux is typically given by Fourier's law of heat conduction:
\begin{eqnarray}
	\boldsymbol{q} = - k \boldsymbol{\nabla} T
\end{eqnarray}
where $k > 0$ is the thermal conductivity.

\subsection{General Form of Mathematical Models}
\paragraph{Fluid mathematical models}
The fluid mathematical models, as implemented in OpenFOAM, are described a number of times elsewhere, for example, \citep{Weller1998}, and are not repeated here in detail.
The standard \foam fluid model implementations have been taken and refactored into the \texttt{Fluid Model} class design described here to allow straight-forward interfacing with the \texttt{solids4foam} modules.
Initially, laminar and turbulent isothermal incompressible Navier-Stokes formulations have been included. It is intended to include additional established methods as future steps, for example, multi-phase, buoyancy, and compressible formulations.
In brief, for isothermal incompressible flows within a static domain, the governing equations take the form of the so-called Navier-Stokes equations \citep{Ferziger2002}:
\begin{eqnarray}
	\oint_\Gamma \boldsymbol{n} \cdot \boldsymbol{v} \; \text{d}\Gamma &=& 0
	\\
	\int_\Omega \frac{\partial \boldsymbol{v}}{\partial t}  \; \text{d}\Omega
	\;+\; \oint_\Gamma \boldsymbol{v} \left[\boldsymbol{n} \cdot \boldsymbol{v} \right]   \; \text{d}\Gamma
	&=&
	\frac{1}{\rho} \oint_\Gamma \boldsymbol{n} \cdot \left[ \eta \boldsymbol{\nabla} \boldsymbol{v} \right] \; \text{d}\Gamma
	\;-\; \frac{1}{\rho} \int_\Omega \, \boldsymbol{\nabla} p \; \text{d}\Omega
	\;+\; \int_\Omega \,\boldsymbol{b} \; \text{d}\Omega
\end{eqnarray} 
where the surface force terms in the momentum equation are decomposed into a viscous term - first term on the right-hand side of the second equation above - and a pressure term - second term on the right-hand side of the second equation above.

\paragraph{Solid mathematical models: \emph{linear geometry}}
In the case where the strains and rotations are small (\ie $\text{det}[\boldsymbol{\epsilon}] \ll 1$), the difference between the deformed and undeformed configurations, as well as the distinction between Eulerian and Lagrangian descriptions,  can be neglected; in addition, the Cauchy and Engineering stress tensors coincide \ie $\boldsymbol{\sigma}_s = \boldsymbol{\sigma}$.
In such cases, the conservation of linear momentum (Equation \ref{eq:momentumEqn}) may be formulated as:
\begin{eqnarray} \label{eq:smallStrainMomentum}
	\int_{\Omega_o} \rho_o \frac{\partial^2 \boldsymbol{u}}{\partial t^2}  \; \text{d}\Omega_o
	\;=\; \oint_{\Gamma_o} \boldsymbol{n}_o \cdot  \boldsymbol{\sigma}_s \; \text{d}\Gamma_o
	\;+\; \int_{\Omega_o} \rho_o \,\boldsymbol{b} \; \text{d}\Omega_o
\end{eqnarray}
where subscript $o$ indicates that a quantity refers to the original undeformed configuration for example, $\Omega_o$ is the initial undeformed volume, $\rho_o$ is the initial density field, and the velocity vector is expressed in terms of the time derivative of the displacement vector $\boldsymbol{v} = \frac{\partial \boldsymbol{u}}{\partial t}$.

Equivalently, the above equation may be expressed in terms of the increment of displacement, $\Delta \boldsymbol{u} = \boldsymbol{u} - \boldsymbol{u}_{[m-1]}$:
\begin{eqnarray} \label{eq:smallStrainMomentumIncremental}
	\int_{\Omega_o} \rho_o \frac{\partial^2}{\partial t^2} (\boldsymbol{u}_{[m-1]} + \Delta \boldsymbol{u})  \; \text{d}\Omega_o
	\;=\; \oint_{\Gamma_o} \boldsymbol{n}_o \cdot  \boldsymbol{\sigma}_s \; \text{d}\Gamma_o
	\;+\; \int_{\Omega_o} \rho_o \,\boldsymbol{b} \; \text{d}\Omega_o
\end{eqnarray}

In the above linear geometry mathematical model, the definition of engineering stress may be given by any constitutive laws consistent with the assumption of small strains, small rotations and linear geometry; for example, constitutive laws presented in Table \ref{table:MechanicalLawLinearGeometry}.
In the case where the definition of engineering stress $\boldsymbol{\sigma}_s$ is a linear function of the displacement vector, for example, Linear (Hookean) elastic in Table \ref{table:MechanicalLawLinearGeometry}, then the momentum equation (Equations \ref{eq:smallStrainMomentum} and \ref{eq:smallStrainMomentumIncremental}) becomes a linear function of displacement field $\boldsymbol{u}$;
however, for nonlinear constitutive laws, for example, Linear elastic, Mises/J$_2$ plastic in Table \ref{table:MechanicalLawLinearGeometry}, the momentum equation depends nonlinearly on the displacement field.
Additional nonlinearity may be introduced through boundary conditions that depend nonlinearly on the displacement vector, for example, contact or cohesive boundary conditions.

\paragraph{Solid mathematical models: \emph{nonlinear geometry}}
When the assumptions of small strains and small rotations are no longer valid, Nanson's relation \citep{Bathe1996, Maneeratana2000}, $\boldsymbol{\Gamma} = J \boldsymbol{F}^{-T} \cdot \boldsymbol{\Gamma}_o$ (or $\boldsymbol{n} = J \boldsymbol{F}^{-T} \cdot \boldsymbol{n}_o$), relating the deformed, $\boldsymbol{\Gamma}$, and initial, $\boldsymbol{\Gamma}_o$, area vectors, may be employed to reformulate the governing momentum equation (Equation \ref{eq:momentumEqn}) in terms of the initial undeformed configuration, indicated by subscript $o$:
\begin{eqnarray} \label{eq:totalLag}
	\int_{\Omega_o} \rho_o \frac{\partial^2 \boldsymbol{u}}{\partial t^2}  \; \text{d}\Omega_o
	\;=\; \oint_{\Gamma_o} (J \boldsymbol{F}^{-T} \cdot \boldsymbol{n}_o) \cdot  \boldsymbol{\sigma} \; \text{d}\Gamma_o
	\;+\; \int_{\Omega_o} \rho_o \,\boldsymbol{b} \; \text{d}\Omega_o
\end{eqnarray}
where the Jacobian of the deformation gradient is $J = \text{det}[\boldsymbol{F}]$, and $\text{det}[\cdot]$ indicates the determinant operator.
Equivalently the governing momentum equation may be expressed in terms of the updated configuration, indicated by subscript $u$:
\begin{eqnarray} \label{eq:updatedLag}
	\int_{\Omega_u} \frac{\partial}{\partial t} \left( \rho_u \frac{\partial \boldsymbol{u}}{\partial t} \right)  \; \text{d}\Omega_u
	\;=\; \oint_{\Gamma_u} (j \boldsymbol{f}^{-T} \cdot \boldsymbol{n}_u) \cdot  \boldsymbol{\sigma} \; \text{d}\Gamma_u
	\;+\; \int_{\Omega_u} \rho_u \,\boldsymbol{b} \; \text{d}\Omega_u
\end{eqnarray}
where, in this case, Nanson's relation relates the deformed, $\boldsymbol{\Gamma}$, and updated, $\boldsymbol{\Gamma}_u$, area vectors, $\boldsymbol{\Gamma} = j \boldsymbol{f}^{-T} \cdot \boldsymbol{\Gamma}_u$ (or $\boldsymbol{n} = J \boldsymbol{F}^{-T} \cdot \boldsymbol{n}_u$), via the relative deformation gradient (Equation \ref{eq:relativeDeformationGradient}) and the relative Jacobian, $j = \text{det}[\boldsymbol{f}]$.

When the momentum equation is integrated over the initial undeformed configuration, Equation \ref{eq:totalLag}, it is referred to as a total Lagrangian (TL) approach, whereas when integrated over the updated configuration, Equation \ref{eq:updatedLag}, it is referred to as an updated Lagrangian (UL) approach.
Both approaches are mathematically equivalent and, as discussed by \citet{Bathe1996}, the only advantage of one formulation over the other lies in its greater numerical efficiency.
Both TL and UL formulations have been discretised and implemented within \texttt{solids4foam}.
As both TL and UL formulations adopt Lagrangian forms of the equation, the conservation of mass is automatically satisfied (\ie convection term is zero).
In the case of the UL approach, the mesh must be moved at the end of each time-step; this mesh movement step is trivial in the case of the standard FE methods (or vertex-based FV methods) where the discrete displacement field resides at the mesh vertices. As the cell-centred FV method yields a discrete displacement field at the cell/element centroids, this means that the displacement field must be interpolated from the cell/element centroids to the vertices, using, for example, inverse distance weighted averaging or linear least squares \citep{Cardiff2014:orthotropicPaper, Cardiff2016:metalForming}.
An alternative approach, adopted by Maneeratana \etal \citep{Maneeratana1999:ACME, Maneeratana1999:ECCM, Maneeratana1999, Maneeratana2000}, is to directly update the geometry (volumes, area vectors, weights) using the assumed kinematic relations, for example, $\Omega_u = j \Omega_u^{[m-1]}$, $\boldsymbol{\Gamma}_u = j \boldsymbol{f}^{-T} \cdot \boldsymbol{\Gamma}_u^{[m-1]} $.

It should also be noted that the TL and UL formulations in Equations \ref{eq:totalLag} and \ref{eq:updatedLag} are expressed in terms of the total displacement vector $\boldsymbol{u}$, but can equally be expressed in terms of the displacement increment vector $\Delta \boldsymbol{u} = \boldsymbol{u} - \boldsymbol{u}_{[m-1]}$.
For the incremental TL formulation, Equation \ref{eq:totalLag} becomes:
\begin{eqnarray} \label{eq:totalLagIncr}
	\int_{\Omega_o} \rho_o \frac{\partial^2}{\partial t^2} (\boldsymbol{u}_{[m-1]} + \Delta\boldsymbol{u})  \; \text{d}\Omega_o
	\;=\; \oint_{\Gamma_o} (J \boldsymbol{F}^{-T} \cdot \boldsymbol{n}_o) \cdot  \boldsymbol{\sigma} \; \text{d}\Gamma_o
	\;+\; \int_{\Omega_o} \rho_o \,\boldsymbol{b} \; \text{d}\Omega_o
\end{eqnarray}
and the incremental UL formulation, Equation \ref{eq:updatedLag} becomes:
\begin{eqnarray} \label{eq:updatedLagIncr}
	\int_{\Omega_u} \frac{\partial}{\partial t} \left( \rho_u \frac{\partial (\boldsymbol{u}_{[m-1]} + \Delta\boldsymbol{u})}{\partial t}  \right) \; \text{d}\Omega_u
	\;=\; \oint_{\Gamma_u} (j \boldsymbol{f}^{-T} \cdot \boldsymbol{n}_u) \cdot  \boldsymbol{\sigma} \; \text{d}\Gamma_u
	\;+\; \int_{\Omega_u} \rho_u \,\boldsymbol{b} \; \text{d}\Omega_u
\end{eqnarray}
The total displacement and incremental displacement formulations are equivalent and the difference between them lies in their numerical efficiency for a given problem.
In the above TL and UL mathematical models, the definition of Cauchy stress may be given by any constitutive law valid for finite strains, for example, constitutive laws presented in Table \ref{table:MechanicalLawNonlinearGeometry}.

In the TL and UL approaches, there are \emph{at least} two forms of nonlinearity:
\begin{itemize}
\item \textbf{material}: the definitions of Cauchy stress $\boldsymbol{\sigma}$ in Table \ref{table:MechanicalLawNonlinearGeometry} are all nonlinear functions of the displacement vector $\boldsymbol{u}$;

\item \textbf{geometric}: as the geometry (for example, area vectors, $\boldsymbol{\Gamma}$, and volumes, $\Omega$) is a function of the displacement field, this introduces nonlinearity where the Cauchy stress is multiplied by scaling terms ($J \boldsymbol{F}^{-T}$ for TL and $j \boldsymbol{f}^{-T}$ for UL) that are a function of displacement.
\end{itemize}

Additional nonlinearity may be introduced through boundary conditions that depend on the displacement vector, for example, contact or cohesive boundary conditions.


\section{Numerical Methods}

\subsection{Discretisation}
The mathematical models presented in the previous section are here discretised using the cell-centred finite volume method.
Similar to FE methods, the cell-centred FV method provides a discrete approximation of the previously presented exact integral equations;
however, unlike standard FE methods, where the strong form of the governing equation is cast into its equivalent weak form before discretisation, the current FV method directly discretises the strong integral form of the governing equations.

\paragraph{Solution Domain Discretisation}
The solution domain is discretised in space and time:
the total simulation time is divided into a finite number of time steps, of fixed or varying size,
and the domain space is divided into a finite number of contiguous convex polyhedral cells, bounded by polygonal faces, that do not overlap and fill the space completely.
In contrast to standard FE methods, where shape functions are specific to the shape of the element, in the current cell-centred FV method no distinction is made between different cell volume shapes, as all general polyhedra (for example, tetrahedra, hexahedra, triangular prism, dodecahedra, \etc) are discretised in the same general fashion.
Further details of the structure of the mesh are given in, for example, \citep{Cardiff2016:blockCoupled, Cardiff2016:metalForming, Weller1998, Jasak2000:linearElasticity}.

Although the presented toolbox is currently based on the cell-centred FV method, the class structure has been designed to be independent of the implemented discretisation method, and as such, an FE implementation could be included in future developments as a solid or fluid model, where the same interface functions and classes are still used.

\paragraph{Equation Discretisation \& Solution Methodology}
Although the equation discretisation and solution procedure are distinct topics, they are jointly discussed here due to their interdependent relationship.

\paragraph{Fluid mathematical models}
The standard approaches implemented within OpenFOAM have been refactored into class-form (\texttt{Fluid Models} in Figure \ref{fig:classStructure1}) and linked with the \texttt{solidsfoam} toolbox.
For incompressible isothermal flow, there are a choice of algorithms, including Pressure-Implicit with Splitting of Operators (PISO), Semi-Implicit Method for Pressure-Linked Equations (SIMPLE), variants of PISO or SIMPLE, and block coupled (pressure-velocity) approaches; further details can be found in \citep{Jasak1996, Weller1998, Clifford2011}.

\paragraph{Solid mathematical models: \emph{linear geometry}}
Historically, implicit cell-centred FV methods have employed so-called segregated or staggered solution procedures, where the governing vector momentum equation is temporarily decoupled into three scalar component equations; these scalar component equations are solved sequentially using iterative linear solvers and outer fixed-point/Picard iterations are performed to re-couple the equations.
The reason for using such segregated solution procedures may be attributed to the genesis of FV solid mechanics methods from FV CFD methods.
In contrast, traditional FE methods typically use \emph{block coupled} solution procedures, where all three displacement components are simultaneously solved in a large block system, typically using a direct linear solver.
In recent years, similar block coupled solution procedures have been developed by \citet{Das2011} and \citet{Cardiff2016:blockCoupled}, echoing standard FE solution procedures.
Although the current article has primarily focussed on developments related to the OpenFOAM software, in particular \citep{Weller1998, Jasak2000:linearElasticity, Jasak2000:contact, Tukovic2007:fsi, Tukovic2007:updatedLagrangainFAMENA, Karac2009, Tukovic2010acp, Tukovic2012, Cardiff2012:contactPaper, Carolan2013, Cardiff2014:orthotropicPaper, Cardiff2014:hipPaper,  Cardiff2014:mappedMusclesPaper, Tukovic2014, Cardiff2015:hydraulicFractures, Cardiff2016:blockCoupled, Cardiff2016:metalForming, Tukovic2017:fsi}, it should be noted that many of the discussed FV discretisation and solution procedures stem from and relate closely to the work of Demird\v{z}i\'{c}, Muzaferija, Ivankovi\'{c} and co-workers \citep{Demirdzic1988, Ivankovic1994, Demirdzic1995, Demirdzic1997, Ivankovic1997, Ivankovic1998, Demirdzic2000, Maneeratana2000, Demirdzic2005, Basic2005, Bijelonja2005, Bijelonja2006, Demirdzic2015:fourthOrder}. In addition, there has been notable contributions to the field of FV solid mechanics by a number of other authors taking related but distinct approaches \citep{Beale1990, Trangenstein1991, Spalding1993, Hattel1995, Trangenstein1994, Fryer1991, Wheel1996, Slone2002, Pan2010, Das2011, Cavalcante2012:generalized, Nordbotten2014}.

In keeping with the OOP philosophy, the structure of the \texttt{Solid Models} class has been designed so that either segregated or coupled solution methodologies may be employed in combination with any of the mechanical constitutive laws. From the user perspective, this means that the choice between solution methodologies can be controlled via input files.

To allow sufficient flexibility in the choice of solution methodology and constitutive relation, the surface force term (divergence of stress) in the conservation of linear momentum is partitioned into \emph{implicit} and \emph{explicit} components;
for the case of linear geometry, the conservation of linear momentum becomes:
\begin{eqnarray} \label{eq:momentumPartitioned}
	\int_{\Omega_o} \rho_o \frac{\partial^2 \boldsymbol{u}}{\partial t^2}  \; \text{d}\Omega_o
	\;=\;
	\overbrace{\oint_{\Gamma_o} \boldsymbol{T}_{\sigma} \; \text{d}\Gamma_o}^{\text{implicit}}
	\;+\;
	\overbrace{ \oint_{\Gamma_o} \boldsymbol{n}_o \cdot  \boldsymbol{\sigma}_s \; \text{d}\Gamma_o
	\;-\; \oint_{\Gamma_o} \boldsymbol{T}_{\sigma} \; \text{d}\Gamma_o
	}^{\text{explicit}}
	\;+\; \int_{\Omega_o} \rho_o \,\boldsymbol{b} \; \text{d}\Omega_o
\end{eqnarray}
where $\boldsymbol{T}_{\sigma}$ is an approximation (or linearisation) of the traction field, $\boldsymbol{n}_o \cdot  \boldsymbol{\sigma}_s$, in terms of the displacement field.
In this context, \emph{implicit} indicates contribution to the matrix of the resulting discretised algebraic linear system and \emph{explicit} indicates contribution to the source vector of the linear system; however, it should be noted that the procedure is also implicit in the time marching sense.
Outer fixed-point/Picard iterations are performed over the equation until explicit terms change less than some predefined tolerance; in that case, the first and third terms on the right-hand side of Equation \ref{eq:momentumPartitioned} cancel out and the calculated displacement field will satisfy the governing equation.
Additional convergence checks may be required to ensure the constitutive relation has converged.
Assuming convergence is achieved, the choice of the implicit traction approximation, $\boldsymbol{T}_{\sigma}$, affects only the rate of convergence, and does not affect the final converged answer.
In the case of a linear (Hookean) elastic law and a segregated solution procedure, \citet{Jasak2000:linearElasticity} has shown the choice of implicit term that ensures optimal rate of convergence is:
\begin{eqnarray}
	\boldsymbol{T}_{\sigma}^{~\text{segregated}} &= (2\mu + \lambda) \; \boldsymbol{n}_o \cdot \boldsymbol{\nabla}_o \boldsymbol{u} \notag \\
	&= (\kappa + \frac{4}{3}\mu) \; \boldsymbol{n}_o \cdot \boldsymbol{\nabla}_o \boldsymbol{u}
\end{eqnarray}
whereas in the case of block coupled procedure the entire traction vector can be treated implicitly:
\begin{eqnarray}
	\boldsymbol{T}_{\sigma}^{~\text{coupled}} &=& \boldsymbol{n}_o \cdot \boldsymbol{\sigma} \notag \\
	&=& \mu \; \boldsymbol{n}_o \cdot \boldsymbol{\nabla}_o \boldsymbol{u}
	+ \mu \; \boldsymbol{n}_o \cdot (\boldsymbol{\nabla}_o \boldsymbol{u})^T
	+ \lambda \; \text{tr}(\boldsymbol{\nabla}_o \boldsymbol{u}) \boldsymbol{n}_o \notag \\
	&=& \mu \; \boldsymbol{n}_o \cdot \boldsymbol{\nabla}_o \boldsymbol{u}
	+ \mu \; \boldsymbol{n}_o \cdot (\boldsymbol{\nabla}_o \boldsymbol{u})^T
	+ (k - \tfrac{2}{3}\mu) \; \text{tr}(\boldsymbol{\nabla}_o \boldsymbol{u}) \boldsymbol{n}_o
\end{eqnarray}
For the case of a linear (Hookean) elastic material and a block coupled approach, the explicit surface force terms in Equation \ref{eq:momentumPartitioned} are zero \ie $\boldsymbol{T}_\sigma = \boldsymbol{n}_o \cdot  \boldsymbol{\sigma}_s$, and the no outer fixed-point/Picard iterations are required \citep{Cardiff2016:blockCoupled}, assuming linear boundary conditions.

In the case of a linear elastic law and segregated or coupled approaches, the optimal choice of implicit term, $\boldsymbol{T}_{\sigma}$, is clear; however, when considering nonlinear constitutive laws (for example, linear elastic with Mises/J$_2$ plasticity) the most appropriate choice is less trivial.
Setting the implicit term, $\boldsymbol{T}_{\sigma}$, based on the linear elastic approximation has been found to provide convergence for implemented constitutive laws, for example, finite strain elastoplasticity \citep{Cardiff2016:metalForming}; however, such a choice may not be optimal in terms of convergence rate or stability.
Insight can be provided by considering the approach taken in standard nonlinear FE methods, where Newton-Raphson iterative schemes are employed.
To ensure the characteristic quadratic rate of asymptotic convergence is achieved for the Newton iterations, the so-called \emph{consistent tangent matrix} must be employed, where the linearisation of the integration algorithm within the mechanical constitutive law is consistent with the the outer Newton tangent \citep{Simo1998, Belytschko2000}.
The tangent matrix for small strains and isotropic linear elasticity, expanded in Voigt notation, is:
\begin{eqnarray}
\boldsymbol{C}^{\text{elastic}} =
\kappa \textbf{I} \otimes \textbf{I} + 2\mu(\textbf{1} - \tfrac{1}{3} \textbf{I} \otimes \textbf{I}) =
\begin{pmatrix}
\kappa + \tfrac{4}{3}\mu & \kappa - \tfrac{4}{3}\mu & \kappa - \tfrac{4}{3}\mu & 0 & 0 & 0 \\
\kappa - \tfrac{4}{3}\mu & \kappa + \tfrac{4}{3}\mu & \kappa - \tfrac{4}{3}\mu & 0 & 0 & 0 \\
\kappa - \tfrac{4}{3}\mu & \kappa - \tfrac{4}{3}\mu & \kappa + \tfrac{4}{3}\mu & 0 & 0 & 0 \\
0 & 0 & 0 & \mu & 0 & 0 \\
0 & 0 & 0 & 0 & \mu & 0 \\
0 & 0 & 0 & 0 & 0 & \mu
\end{pmatrix}
\end{eqnarray}
where $\kappa + \tfrac{4}{3} \mu = 2\mu + \lambda$, $\textbf{1}$ is the fourth-order identity tensor, and $\otimes$ indicates an outer tensor product.
The segregated algorithm typically used in cell-centred FV methods, as described by \citet{Jasak2000:linearElasticity}, effectively approximates the divergence of stress by a simple decoupled Laplacian/diffusion term:
\begin{eqnarray}
\boldsymbol{\nabla} \cdot \boldsymbol{\sigma} &=&
\boldsymbol{\nabla} \cdot \left[ \boldsymbol{C} : \boldsymbol{\nabla} \boldsymbol{u} \right]  \\ \notag
&=&
\boldsymbol{\nabla} \cdot \left[ K_{\text{imp}} \boldsymbol{\nabla} \boldsymbol{u} \right] + [\text{explicit correction}]
\end{eqnarray}
where the stiffness for this implicit Laplacian term, $K_{\text{imp}}$, is set to be the upper left diagonal of the tangent matrix \ie $K_{\text{imp}} = \kappa + \tfrac{4}{3} \mu$.
As described above, the off-diagonal components, providing coupling between the components of displacement, are then treated explicitly via outer deferred corrections.

Inspired by the consistent tangent matrix used with a Newton method, it is possible to consider other choices for the $K_{\text{imp}}$ implicit stiffness term when using a segregated algorithm and fixed-point/Picard outer iterations.
For example, consider the consistent tangent matrix for J$_2$ perfect plasticity with isotropic linear elasticity \citep{Simo1998}:
\begin{eqnarray}
\boldsymbol{C}^{\text{J$_2$}} &=&
\kappa \textbf{I} \otimes \textbf{I}
+ 2\mu \Theta (\textbf{1} - \tfrac{1}{3} \textbf{I} \otimes \textbf{I})
- 2\mu \Theta \boldsymbol{N} \otimes \boldsymbol{N} \\ \notag
&=&
\begin{pmatrix}
\kappa + \mu\Theta(\tfrac{4}{3} - 2 N_{11}^2)
& \kappa - \mu\Theta(\tfrac{4}{3} + 2 N_{11} N_{22})
& \kappa - \mu\Theta(\tfrac{4}{3} + 2 N_{11} N_{33}) & 0 & 0 & 0 \\
\kappa - \mu\Theta(\tfrac{4}{3} + 2 N_{22} N_{11})
& \kappa + \mu\Theta(\tfrac{4}{3} - 2 N_{22}^2)
& \kappa - \mu\Theta(\tfrac{4}{3} + 2 N_{22} N_{33}) & 0 & 0 & 0 \\
\kappa - \mu\Theta(\tfrac{4}{3} + 2 N_{33} N_{11})
& \kappa - \mu\Theta(\tfrac{4}{3} + 2 N_{33} N_{22})
& \kappa + \mu\Theta(\tfrac{4}{3} - 2 N_{33}^2) & 0 & 0 & 0 \\
0 & 0 & 0 & \mu\Theta & 0 & 0 \\
0 & 0 & 0 & 0 & \mu\Theta & 0 \\
0 & 0 & 0 & 0 & 0 & \mu\Theta
\end{pmatrix}
\end{eqnarray}
where $\Theta = 1 - \tfrac{2\mu \Delta\gamma}{||\boldsymbol{s}_{n+1}^{\text{trial}}||}$, and note that $\Theta \rightarrow 1$ as $\Delta \gamma \rightarrow 0$, and $\Theta \leq 1 $.
By once again taking the top left diagonal coefficients, this suggest the following vector coefficient for the implicit Laplacian term to be used with the segregated FV procedure:
\begin{eqnarray}
\boldsymbol{K}^{\text{imp}} &=&
\begin{pmatrix}
\kappa + \mu\Theta(\tfrac{4}{3} - 2 N_{11}^2) \\
\kappa + \mu\Theta(\tfrac{4}{3} - 2 N_{22}^2) \\
\kappa + \mu\Theta(\tfrac{4}{3} - 2 N_{33}^2)
\end{pmatrix}
\end{eqnarray}
It should, however, be recognised that the asymptotic convergence rate of fixed-point/Picard iterations is linear, in contrast to the superior quadratic convergence of the Newton methods.
Consequently, a larger number of total outer iterations would be expected in the current segregated FV implementation; however, through the use of efficient iterative linear solvers, each outer iteration is much less expensive, resulting in a competitive method for certain classes of problems.
For the segregated approach, the use of the Laplacian coefficient indicated by the above analysis may reduce the overall number of other iterations; however, limited efficiency increases would be expected given the use of fixed-point iterations and efficient iterative solvers for the inner system.
In contrast, given the close similarly of the coupled approach \citep{Cardiff2016:blockCoupled} to standard FE implementations, the use of full/quasi Newton iteration with consistent tangent matrices would be a promising direction to explore and will be the focus of future developments.

\paragraph{Solid mathematical models: \emph{nonlinear geometry}}
Similar to linear geometry case above, the TL and UL geometrically nonlinear momentum equations can be re-formulated, respectively, as:
\begin{eqnarray} \label{eq:momentumPartitionedTL}
	\int_{\Omega_o} \rho_o \frac{\partial^2 \boldsymbol{u}}{\partial t^2}  \; \text{d}\Omega_o
	\;=\;
	\overbrace{\oint_{\Gamma_o} \boldsymbol{T}_{\sigma} \; \text{d}\Gamma_o}^{\text{implicit}}
	\;+\;
	\overbrace{
	\oint_{\Gamma_o} (J \boldsymbol{F}^{-T} \cdot \boldsymbol{n}_o) \cdot  \boldsymbol{\sigma} \; \text{d}\Gamma_o
	\;-\; \oint_{\Gamma_o} \boldsymbol{T}_{\sigma} \; \text{d}\Gamma_o
	}^{\text{explicit}}
	\;+\; \int_{\Omega_o} \rho_o \,\boldsymbol{b} \; \text{d}\Omega_o
\end{eqnarray}
\begin{eqnarray} \label{eq:momentumPartitionedUL}
	\int_{\Omega_u} \rho_u \frac{\partial^2 \boldsymbol{u}}{\partial t^2}  \; \text{d}\Omega_u
	\;=\;
	\overbrace{\oint_{\Gamma_u} \boldsymbol{T}_{\sigma} \; \text{d}\Gamma_u}^{\text{implicit}}
	\;+\;
	\overbrace{
	\oint_{\Gamma_u} (j \boldsymbol{f}^{-T} \cdot \boldsymbol{n}_u) \cdot  \boldsymbol{\sigma} \; \text{d}\Gamma_u
	\;-\; \oint_{\Gamma_u} \boldsymbol{T}_{\sigma}  \; \text{d}\Gamma_u
	}^{\text{explicit}}
	\;+\; \int_{\Omega_u} \rho_u \,\boldsymbol{b} \; \text{d}\Omega_u
\end{eqnarray}

The relations above are given for the case where the total displacement, $\boldsymbol{u}$, is the primary unknown; in the case that the increment of displacement, $\Delta \boldsymbol{u}$, is the primary unknown, then the total displacement, $\boldsymbol{u}$, in the inertia term on the left-hand side of the equations above is replaced by $\boldsymbol{u}_{[m-1]} + \Delta\boldsymbol{u}$.

\paragraph{System of linear equations}
The above mathematical models are subsequently discretised using the cell-centred FV method, where linear variations of the displacement field are assumed across each cell; the resulting discretisation is second-order accurate in space, where the discretisation error reduces at a second-order rate as the cell size is reduced.
Specific details of the discretisation of each of the terms (for example, inertia, Laplacian, divergence) can be found in \citep{Cardiff2016:metalForming, Cardiff2016:blockCoupled, Cardiff2014:orthotropicPaper, Jasak2000:linearElasticity, Weller1998}.
The final discretised form of the linear momentum equation results in a system of linear algebraic equations of the form:
\begin{eqnarray}
[\boldsymbol{A}][\boldsymbol{\phi}] = [\boldsymbol{b}]
\end{eqnarray}
where $[\boldsymbol{A}]$ is a $N \times N$ sparse matrix with weak diagonal dominance, for the segregated algorithm, and $N$ is the number of cells in the mesh.
For the case of a block coupled solution procedure, the matrix is no longer diagonally dominant and $N$ represents the sum of the number of cells and boundary faces.
The solution vector $[\boldsymbol{\phi}]$ contains the unknown cell-centre displacements, $\boldsymbol{u}$, or displacement increments, $\Delta \boldsymbol{u}$, and $\boldsymbol{b}$ is the source vector containing the explicit discretised terms, body force terms and boundary condition contributions.
For the coupled approach, the solution vector also contains the boundary face unknown displacements or displacement increments.
The inner linear sparse system can be solved using an iterative solver, such as the incomplete Cholesky pre-conditioned conjugate gradient (ICCG) method \citep{Jacobs1980}, or directly using Gaussian elimination or LU decomposition.
It should be noted that the standard OpenFOAM library contains only iterative linear solvers implementations; however, it is possible to use direct solvers by linking with external libraries, such as Eigen \citep{eigenweb}, MUMPS \citep{MUMPS}, PETSc \citep{petsc-web-page} or Trilinos \citep{Trilinos}.
Algebraic multi-grid methods can be used to accelerate the iterative methods.
The iterative linear solvers employed here, within the \foam library, allow for efficient parallelisation on distributed memory computer cluster via the method of domain decomposition;
this involves partitioning the mesh into sub-regions, each of which are solved on a separate CPU core.
Inter-CPU-core communication, via the Message Passing Interface (MPI) protocol, provides the necessary information to solve the system of linear equations.

For the segregated solution procedure, the use of iterative linear solvers has been found to be more efficient; whereas, for the block coupled procedure, direct linear solvers may be more efficient for smaller system and iterative solvers for larger ones.
Although not implemented here, the use of geometric (as opposed to algebraic) multi-grid methods would be expected to provide significant speedups for both linear geometry and nonlinear cases, as demonstrated by a number of authors, for example, \citep{Muzaferija1994, Demirdzic1997, Fainberg1996, Ivankovic1997}; however, this approach comes with the additional overhead of generating meshes, a task that is not trivial to automate.

In order to improve stability, it is occasionally necessary to employ under-relaxation of the solution variable ($\boldsymbol{u}$ or $\Delta \boldsymbol{u}$) via field under-relaxation or linear equation under-relaxation.
As outer iterations are performed within each time-step, the use of under-relaxation does not affect the predicted transient response.
In addition, to avoid the appearance of checkerboard oscillations, a so-called Rhie-Chow diffusive term \citep{Rhie1983} is added within the momentum equation; this correction term reduces at a third-order rate as the mesh is refined and consequently does not affect the second-order accuracy of the spatial discretisation.

It should be noted that the solid model formulations presented here all employ a continuum element approach, and as yet, reduced order approaches like beam and plate/shell formulations have not been implemented; this point is discussed further in the future outlook section at the end of article.


\paragraph{Boundary Conditions}
In addition to implementation of standard (constant or time-varying, uniform or non-uniform) displacement, traction and symmetry conditions, a number of nonlinear boundary conditions have been implemented:
for example, normal and frictional contact models based on a penalty formulation, for example, see \citep{Cardiff2012:contactPaper, Cardiff2016:metalForming}, as well as cohesive zone model fracture conditions, for example, \citep{Cardiff2015:hydraulicFractures, Lee2015}.
As in standard FE methods, Dirichlet conditions, for example, displacement, are strongly enforced; however, in contrast to FE methods, Neumann conditions, for example, traction conditions, are also strongly enforced within the current cell-centred FV discretisation.
The implications of this difference, for solid mechanics problems, may be of little consequence, assuming mesh independence is achieved.

\subsection{Implementation}
For illustrative purposes, the algorithm for the segregated linear geometry (small strains) \texttt{Solid Model} class is shown in Algorithm \ref{algorithm:linGeomSolid} below.
Notice that the general solution methodology is independent of the details of the stress field definition, allowing straight-forward inclusion of new constitutive laws without having to re-implement, test and verify the outer solution procedure.
For the interested reader, the corresponding \foam C$++$ code is given in Appendix \ref{app:linGeomSolid}, where it is clear to see a distinguishing feature of \foam language: the coding syntax closely resembles the mathematical equations being solved.
\begin{algorithm}
\caption{Solution procedure: linear geometry segregated algorithm}
\label{algorithm:linGeomSolid}
\begin{algorithmic}[1]
\ForAll{time steps}:
	\While {explicit terms are not converged}
		\State Momentum equation: assemble and solve in terms of $\boldsymbol{u}$ (or $\Delta \boldsymbol{u}$)
		\Comment{Equation \ref{eq:momentumPartitioned}}
		\State Update kinematics for example, $\boldsymbol{\nabla} \boldsymbol{u}$, $\boldsymbol{\epsilon}$
		\State Update stress
		\Comment{Calculation of stress by run-time selectable mechanical law}
	\EndWhile \textbf{end while}
\EndFor \State \textbf{end for}
\end{algorithmic}
\end{algorithm}

\section{Example Problems}
The current section presents the analysis of a number of contrasting physical problems, to highlight the features and capabilities of the presented \texttt{solids4foam} toolbox.
In addition to the cases presented here, a number of problems presented in previous publications are also included within the \texttt{solids4foam} toolbox; for example, contact mechanics \citep{Cardiff2012:contactPaper}, isotropic and orthotropic linear elasticity \citep{Cardiff2016:blockCoupled, Cardiff2014:orthotropicPaper, Tukovic2012}, poro-elasto-plasticity \citep{Tang2015}, metal forming \citep{Cardiff2016:metalForming}, and FSI benchmarks \citep{Tukovic2014, Tukovic2017:fsi}.


\subsection{Heated spherical pressure vessel}
This case consists of a hollow spherical vessel that is internally heated and pressurised, providing an example of a transient coupled temperature-displacement analysis;
the case has been analysed a number of times previously, for example, \citep{SimScaleDoc:sphericalPresureVessel, Afkar2014}.
Although the solution is 1-D axisymmetric, the problem is analysed here as 3-D with symmetries for demonstrative purposes as well as allowing assessment of parallel efficiency.
\begin{figure}[htb]
	\centering
	\subfigure[Geometry and loading conditions]
	{
		\includegraphics[width=0.48\textwidth]{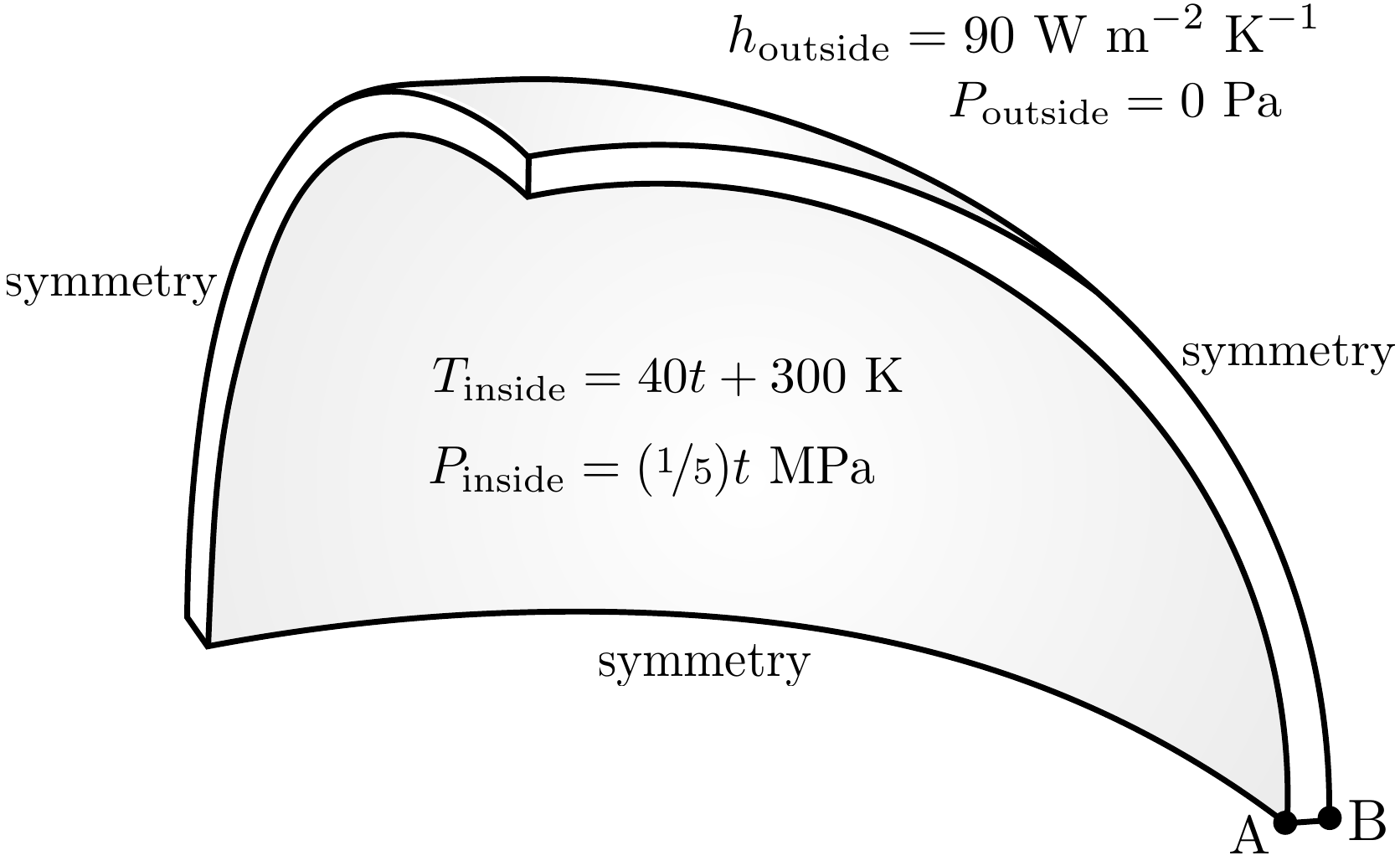}
		\label{fig:hotSphereGeom}
	}
	\subfigure[Mesh containing $9~375$ cells]
	{
		\includegraphics[width=0.4\textwidth]{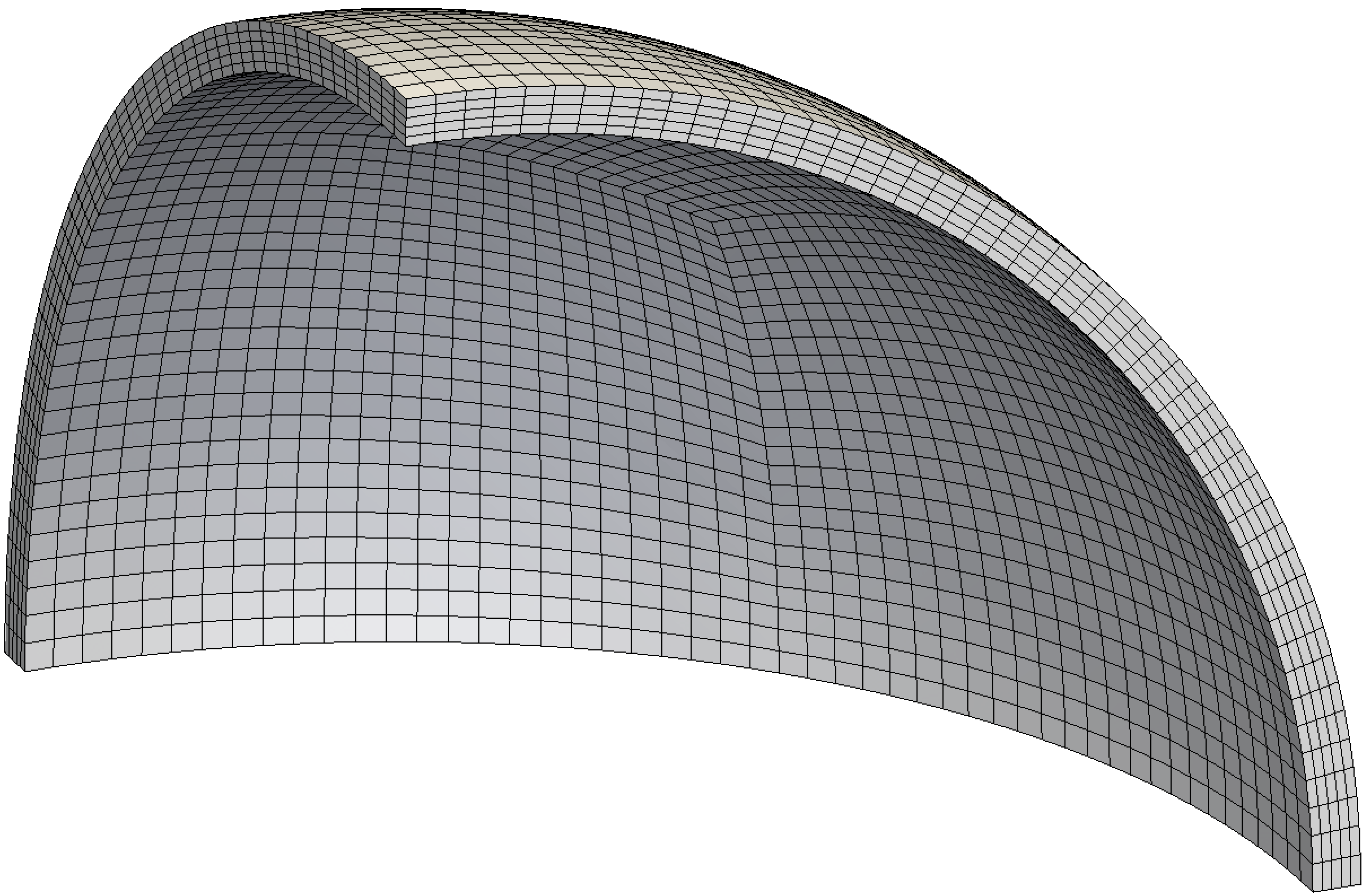}
		\label{fig:hotSphereMesh}
	}
	\caption{Heated spherical pressure vessel: problem geometry, loading conditions and mesh}
\end{figure}

The one-eighth spherical geometry is shown in Figure \ref{fig:hotSphereGeom}, where the inner radius is 190 \milli\meter\ and the outer radius is 200 \milli\meter.
A structured hexahedral mesh has been generated using the OpenFOAM \texttt{blockMesh} mesh utility; in this case, \texttt{blockMesh} from the {OpenFOAM-dev} fork \citep{OpenFOAM.org} of OpenFOAM has been used as it allows meshing of spherical surfaces.
Four successively refined hexahedral meshes have been employed containing $9~375$, $75~000$, $600~000$ and $4~800~000$ cells respectively; the coarsest mesh is shown in Figure \ref{fig:hotSphereMesh}.

As indicated in Figure \ref{fig:hotSphereGeom}, the loading conditions consist of a time-varying pressure and temperature on the inside of the sphere:
\begin{eqnarray}
	P_{\text{inside}} = \left(\frac{1}{5}\right)t \quad\text{\mega\pascal}\notag \\
	T_{\text{inside}} = 40t + 300  \quad\text{\kelvin}
\end{eqnarray}
where $t$ is time, varying from 0 to 5 \second.
The outer surface is specified as traction free, with a surface convective heat coefficient of $h = 90$ \watt\per\meter\squared\kelvin; the reference atmospheric temperature is assumed to be $T_\infty = 300$ \kelvin.
As one-eighth of the geometry is modelled, three symmetry conditions are employed.

The assumed thermo-material properties are given in Table \ref{table:hotSphericalPressureVessel_properties}, where the stress-free reference temperature is $300$ \kelvin.
\begin{table}[htb]
  \centering
	  \ra{1.3}
		\begin{tabular}{@{}lll@{}}
		\toprule
		Young's modulus & $E$ & 190 \giga\pascal \\
		Poisson's ratio & $\nu$ & 0.305 \\
		Density & $\rho$ & 7750 \kilo\gram\per\meter\cubed \\
		Coefficient of thermal expansion & $\alpha$ & $9.7\times10^{-6}$ \per\kelvin \\
		Specific heat capacity & $C$ & 486 \joule\per\kelvin \\
		Thermal conductivity & $k$ & 20 \watt\per\metre\kelvin \\
		\bottomrule
		\end{tabular}
\caption{Heated spherical pressure vessel: thermo-mechanical properties}
\label{table:hotSphericalPressureVessel_properties}
\end{table}

At the initial time, $t = 0$ \second, the temperature of the domain is 300 \kelvin, the displacement is $(0~0~0)$ \meter, and the initial velocity is $(0~0~0)$ \meter\per\second.
A constant time-step of 1 \second\ is employed, simulating five time-steps in total, from 0 to 5 \second; the time-step size of 1 \second\ was selected as it was found to produce time-step independent results on all meshes.

\begin{figure}[htb]
	\centering
	\includegraphics[width=0.6\textwidth]{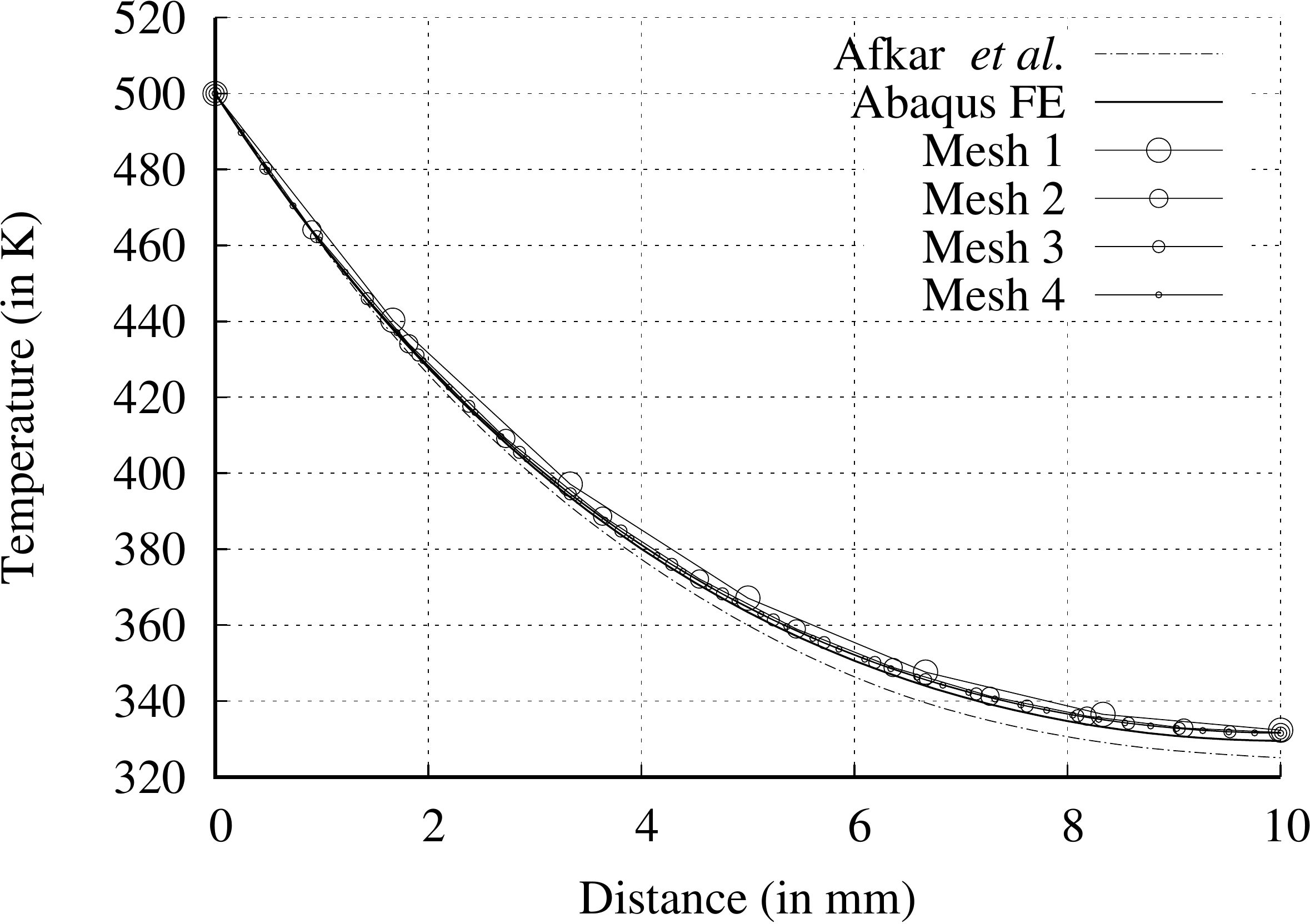}
	\caption{Heated spherical pressure vessel: temperature distribution across the wall thickness at $t = 5$ \second}
	\label{fig:hotSphericalPressureVessel_temperature}
\end{figure}

Figure \ref{fig:hotSphericalPressureVessel_temperature} shows the temperature distribution across the wall thickness from point A to point B at $t = 5$ \second\, where the position of points A and B are indicated in Figure \ref{fig:hotSphereGeom}.
The predictions are shown for the four meshes; for reference, results are given from \citet{Afkar2014} and predictions generated with commercial finite element software Abaqus are shown using the finest mesh; the FE model uses axisymmetric full integration bi-linear elements (Abaqus element code: CAX4T).
Little difference can be seen between the predictions on the different meshes, indicating the discretisation error is small in all cases; this is an expected result given the temperature distribution is close to linear.

\begin{figure}[htb]
	\centering
	\subfigure[Von Mises stress distribution]
	{
		\includegraphics[width=0.48\textwidth]{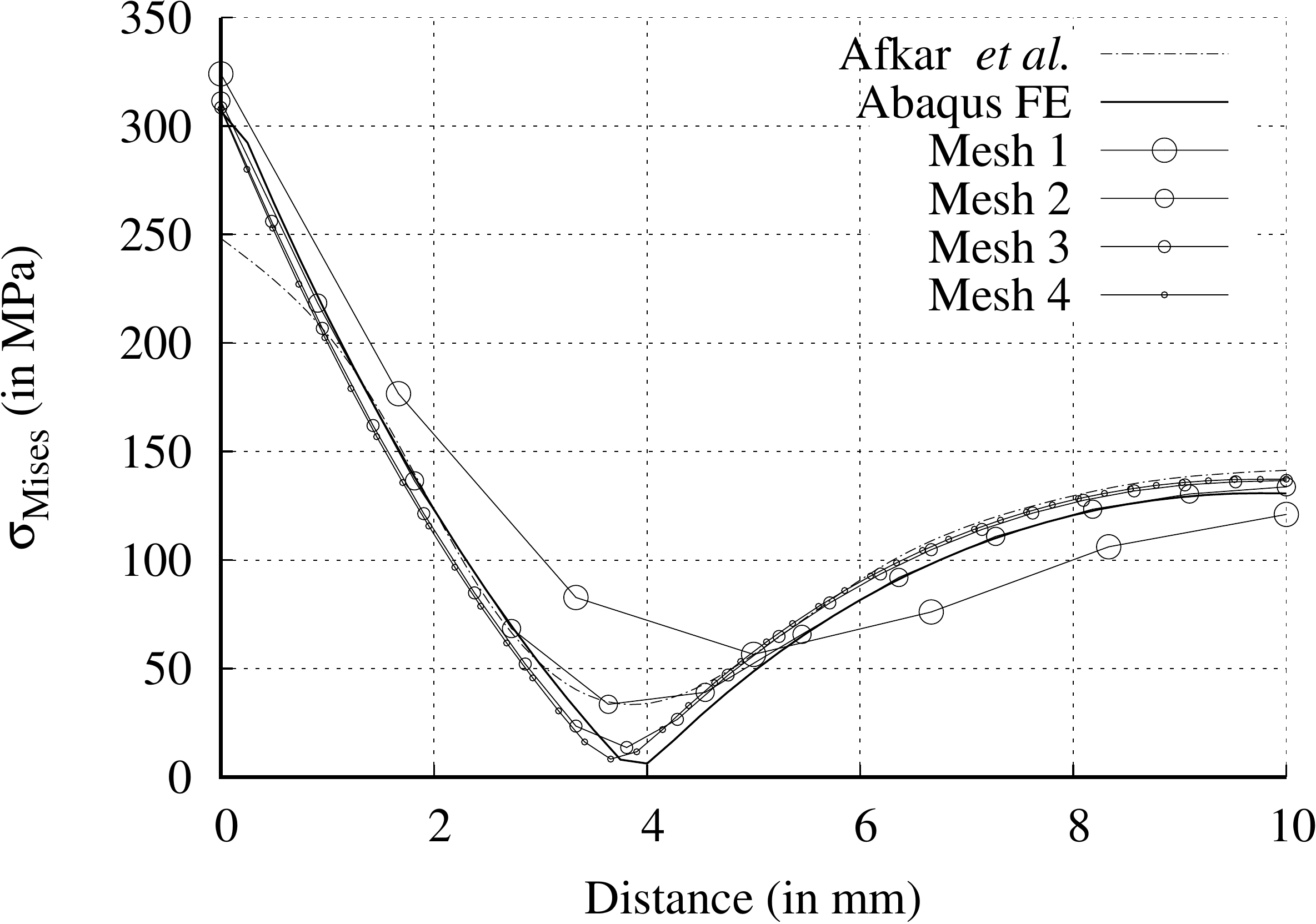}
		\label{fig:hotSphericalPressureVessel_mises}
	}
	\subfigure[Reduction of discretisation error in the von Mises with increasing mesh refinement]
	{
		\includegraphics[width=0.48\textwidth]{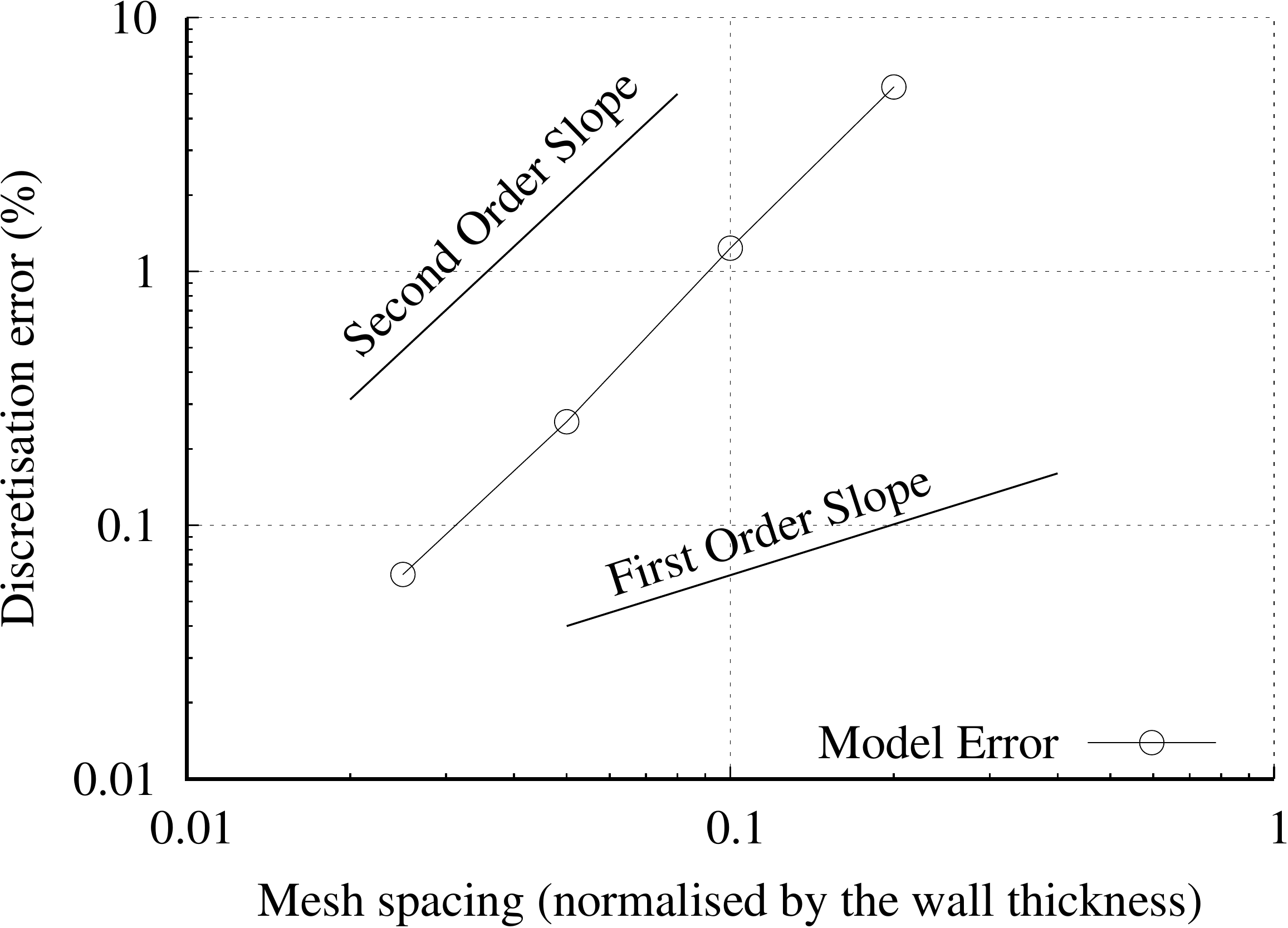}
		\label{fig:hotSphereError}
	}
	\caption{Heated spherical pressure vessel: von Mises stress distribution across the wall thickness at $t = 5$\second\ and the reduction in discretisation error as the cell size is reduced}
\end{figure}

Figure \ref{fig:hotSphericalPressureVessel_mises} shows the von Mises distribution across the wall thickness, from point A to point B, for the four meshes; once again for comparison results are given from \citet{Afkar2014} and using the finest mesh with FE software Abaqus.
As the mesh density is increased, the predictions can be seen to become mesh independent, with even the coarsest mesh capturing the trend.
By examining the error in the von Mises stress prediction at the inner surface, the discretisation error, shown in Figure \ref{fig:hotSphereError}, can be seen to reduce at an approximately second error rate \ie discretisation error reduces by a factor of four when the mesh spacing is halved.
The reference solution is approximated using Richardson's extrapolation \citep{Roache1997}.

\begin{figure}[htb]
	\centering
	\subfigure[Simulation time]
	{
		\includegraphics[width=0.48\textwidth]{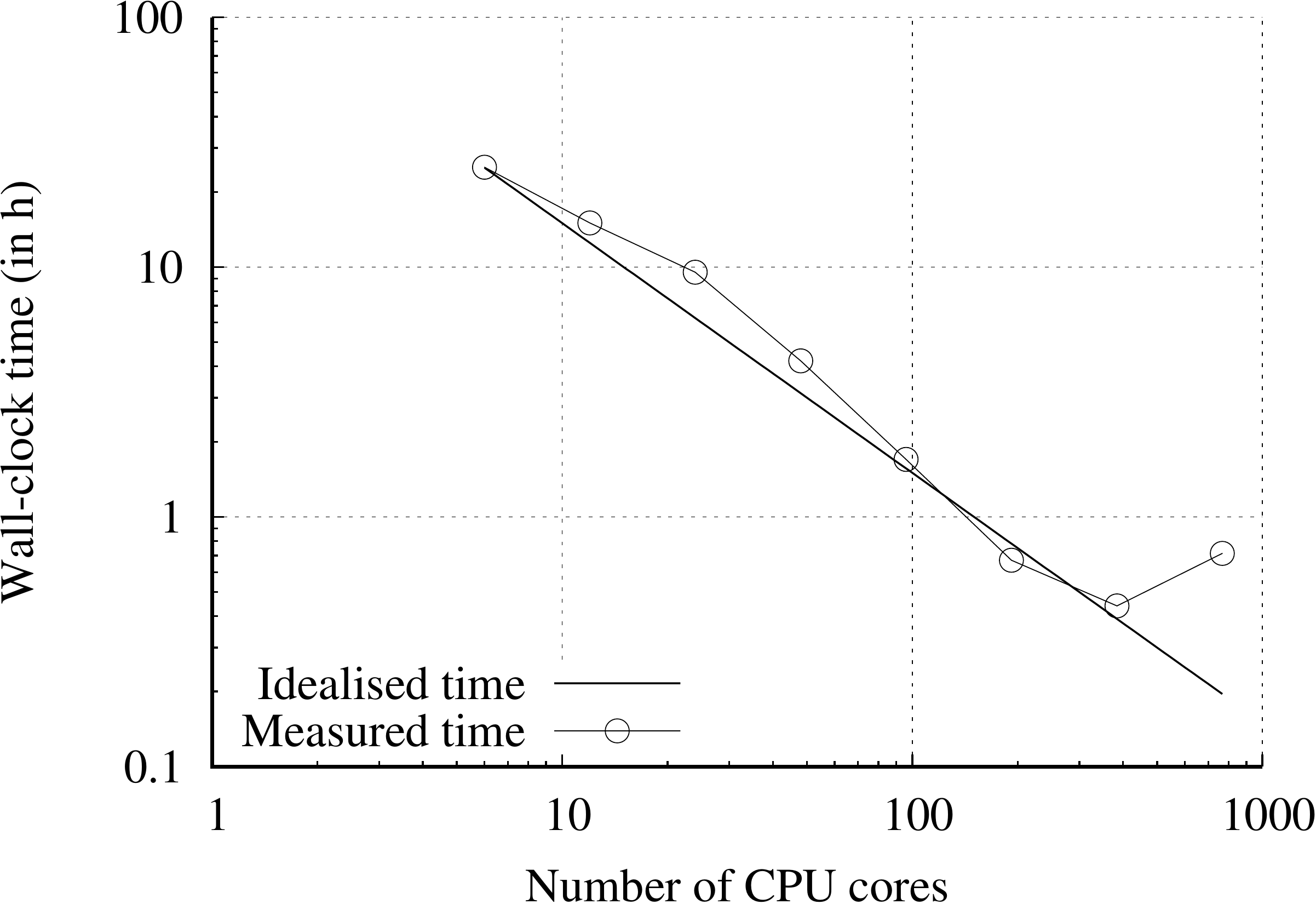}
		\label{fig:hotSphereParallelTimes}
	}
	\subfigure[Parallel speedup]
	{
		\includegraphics[width=0.48\textwidth]{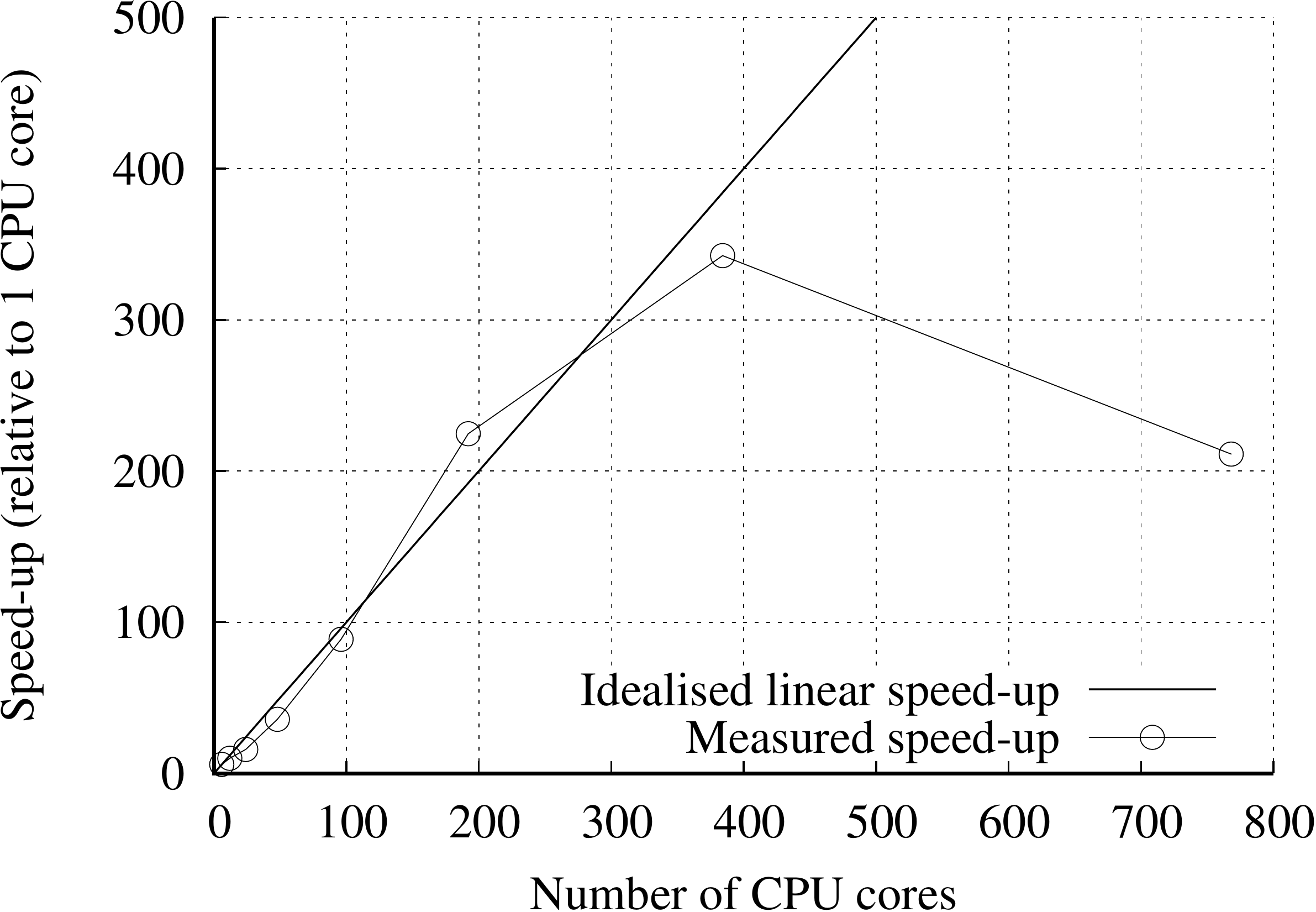}
		\label{fig:hotSphereParallelSpeedup}
	}
	\caption{Heated spherical pressure vessel: parallel efficiency for the mesh containing 4.8 million cells}
\end{figure}

To illustrate the parallel efficiency of the method, the finest mesh, containing 4.8 million cells, is solved using increasing numbers of CPU cores; the METIS decomposition method \citep{Karypis1999} has been used to decompose the domain.
In the ideal case, the method should scale linearly \ie relative to the time taken to simulate the case using one CPU core, it should take one tenth of this time using ten CPU cores, one hundredth of the time using 100 CPU cores, and one thousandth of the time using 1000 CPU cores, etc.
The time taken to analyse the case using increasing numbers of CPU cores is shown in Figure \ref{fig:hotSphereParallelTimes}; the parallel speedup, defined as $\nicefrac{T_{\text{core=1}}}{T_{\text{core=N}}}$, is shown in Figure \ref{fig:hotSphereParallelSpeedup}, where $T_{\text{core=1}}$ is the time taken to simulate the problem using 1 CPU core and $T_{\text{core=N}}$ is the time taken using $N$ CPU cores; in this case $N$ was varied from 6 to 768 CPU cores.
As it would have taken an excessive amount of time to solve the problem on one CPU core, the time taken to solve the problem on one CPU core has been approximated as $T_{\text{core=1}} \approx 6 \times T_{\text{core=6}}$.
The simulation wall-clock times are given for reference in Table \ref{table:hotSphereParallelTimes}, where the approximate number of cells per CPU core are indicated.
The \emph{Fionn} supercomputer from the Irish Centre for High-End Computing has been used for all calculations, where each compute node contains two Intel$^{\tiny{\textregistered}}$ Xeon$^{\tiny{\textregistered}}$ (E5-2660 v2 @ 2.20GHz, 25.6 \mega\byte\ of cache) CPUs and 64 GiB of RAM.

In the current parallel efficiency analysis, the total number of cells in the mesh remains constant and the question becomes \emph{how many cores should we use to get answer in the shortest time?}
This corresponds to a so-called \emph{strong} parallel scalability test, where the total amount of computational work remains constant and the number of CPU cores is increased to try reduce the time required.
In contrast, a \emph{weak} parallel scalability aims to keep the amount of work per CPU core constant, while increasing the number of CPU cores used \ie \emph{for a given amount of time and assuming we are not restricted in the number of CPU cores we can use, how large a problem can we solve?}
From Figures \ref{fig:hotSphereParallelTimes} and \ref{fig:hotSphereParallelSpeedup}, it can be seen that the method scales in an approximately linear fashion up until 384 CPU cores \ie the solution is attained approximately 384 times faster using 384 CPU cores than using 1 CPU cores.
This corresponds to approximately $12~500$ cells per CPU core.
As the number of CPU cores are increased further to 768, the parallel efficiency drops and in this case attaining the solution using 768 CPU cores is in fact slower than using 384 CPU cores;
this can be explained by noting that parallel performance is a competition between computational work done per core versus the amount of information communicated between the cores.
As the number of the cores is increased, the amount of work done per core is decreasing while the amount of communication between the cores is increasing; at a certain point, the time taken for communication between the cores becomes significant and increasing the number of cores further can be expected to increase the simulation time.

As a rule of thumb for FV methods in CFD using standard iterative solution methodologies, achieving good parallel efficiency requires approximately $20 \times10^3$ to $50 \times10^3$ cells per CPU core; in the current case, the noted behaviour was slightly more efficient requiring only $12.5 \times10^3$ cells per core, where the solution for the finest mesh was achieved in under half an hour using 384 CPU cores, in comparison to over 25 hours when using 1 CPU core; this is in keeping with previous findings for FV solid mechanics, for example, \citep{Cardiff2014:orthotropicPaper}.
\begin{table}[htb]
  \centering
	  \ra{1.3}
		\begin{tabular}{@{}llll@{}}
		\toprule
		Number of CPU cores & Wall-clock time (in \hour) & Speed-up & Cells per CPU core \\
		\toprule
		1 & 150.70 (\emph{estimated}) & - & $4800\times10^3$ \\
		6 & 25.12 & 6 & $800\times10^3$ \\
		12 & 15.01 & 10.04 & $400\times10^3$ \\
		24 & 9.54 & 15.80 & $200\times10^3$ \\
		48 & 4.21 & 35.77 & $100\times10^3$ \\
		96 & 1.70 & 88.74 & $50\times10^3$ \\
		192 & 0.67 & 224.65 & $25\times10^3$ \\
		384 & 0.44 & 342.50 & $12.5\times10^3$ \\
		768 & 0.71 & 211.18 & $6.25\times10^3$ \\
		\bottomrule
		\end{tabular}
\caption{Heated spherical pressure vessel: parallel simulation times for the mesh 4.8 million cells}
\label{table:hotSphereParallelTimes}
\end{table}

\subsection{A 3-D Elastic Punch}
In this case, a flat punch with rounded edges is pressed into a flat-topped cylinder.
The case has been proposed as a contact mechanics benchmark by the National Agency for Finite Element Methods and Standards (NAFEMS) \citep{NAFEMS, NAFEMS:contact}, and has been examined previously a number of times, for example, \citep{CodeAster:documentation}.

\begin{figure}[htb]
	\centering
	\subfigure[Geometry]
	{
		\includegraphics[width=0.48\textwidth]{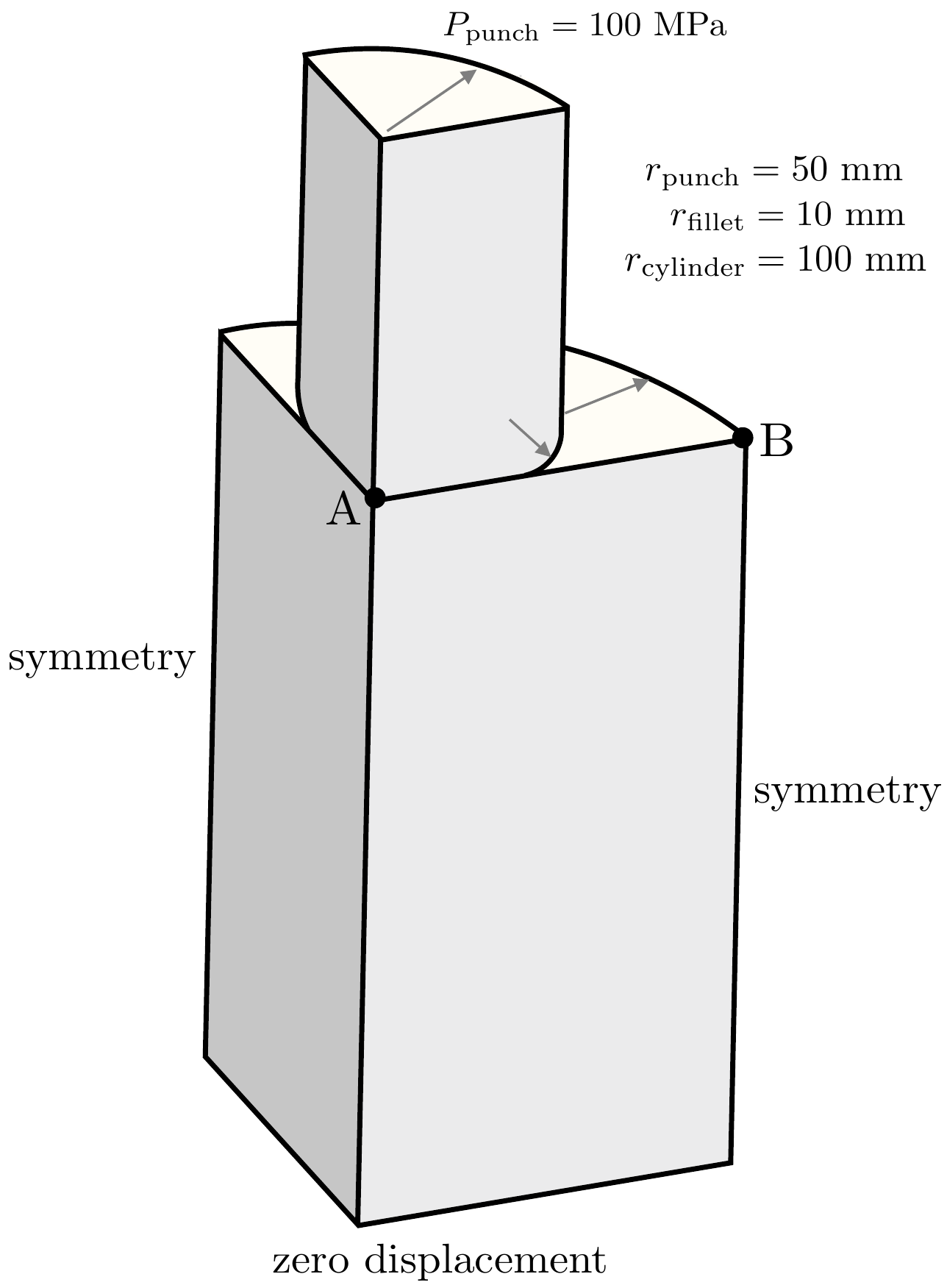}
		\label{fig:punchGeom}
	} \quad \quad
	\subfigure[Mesh containing 4~946 cells]
	{
		\includegraphics[width=0.29\textwidth]{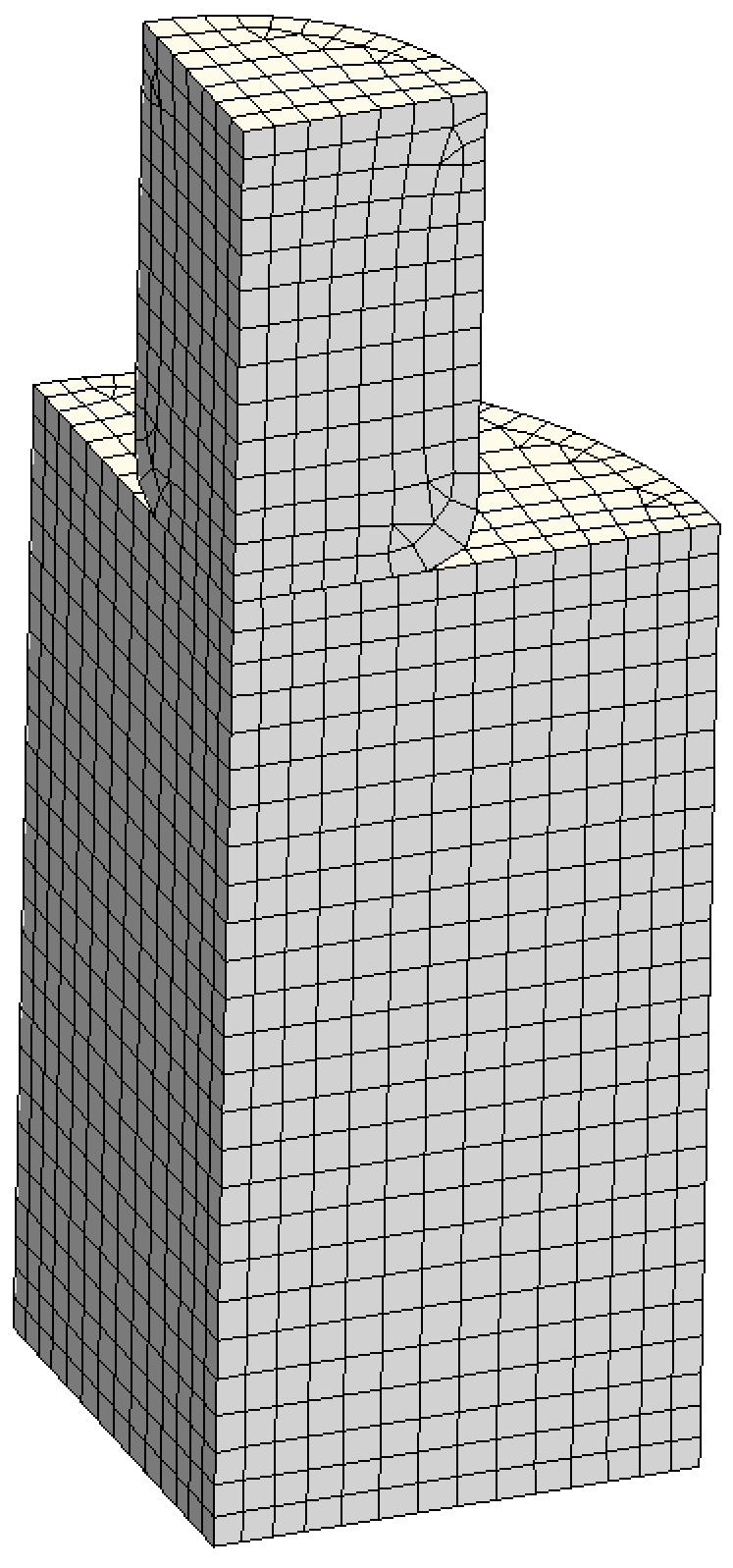}
		\label{fig:punchMesh}
	}
	\caption{A 3-D Elastic Punch: problem geometry and mesh}
\end{figure}

The case demonstrates a 3-D quasi-static contact analysis of differing materials; the case and solution is expected to be 2-D axisymmetric but is simulated here as 3-D with symmetries for demonstrative purposes.
The problem geometry, shown in Figure \ref{fig:punchGeom}, consists of a cylindrical punch with radius 50 \milli\meter\ and a 10 \milli\meter\ filleted edge, in contact with a flat-topped cylinder of 100 \milli\meter\ radius.
The geometry has been created using open-source software FreeCAD \citep{FreeCAD} and exported in the IGES format to the open-source Gmsh software \citep{Geuzaine2009} to create a facetted STL surface file; subsequently, a Cartesian mesh has been created using the OpenFOAM meshing utility \texttt{cfMesh} \citep{cfmeshweb}.
Four approximately uniform Cartesian-based meshes of increasing refinement \-- 4~946, 31~513, 228~026 and 1~789~600 cells \-- have been employed, where the coarsest mesh is shown in Figure \ref{fig:punchMesh}.

The punch is assumed to be steel and the cylinder is assumed to be aluminium, where the mechanical properties are given in Table \ref{table:punch_properties}; it should be noted that the Poisson's ratios need not be set the same, and have only been done so here to allow comparison with the NAFEMS results.
\begin{table}[htb]
  \centering
	  \ra{1.3}
		\begin{tabular}{@{}lll@{}}
		\toprule
		\emph{Punch} & & \\
		Young's modulus & $E$ & 210 \giga\pascal \\
		Poisson's ratio & $\nu$ & 0.3 \\
		\bottomrule
		\emph{Cylinder} & & \\
		Young's modulus & $E$ & 70 \giga\pascal \\
		Poisson's ratio & $\nu$ & 0.3 \\
		\bottomrule
		\end{tabular}
\caption{A 3-D Elastic Punch: mechanical properties}
\label{table:punch_properties}
\end{table}

A pressure of 100 MPa is applied to the upper surface of the punch, and the lower surface of the cylinder is fixed \ie displacement is zero; in addition, two symmetry plane conditions are employed as indicated in Figure \ref{fig:punchGeom}.
To enforce the contact constraint conditions between the punch and cylinder surfaces, a penalty formulation is employed for both the normal and friction models; further details are given previously \citep{Cardiff2012:contactPaper, Cardiff2016:metalForming}.
The case is solved in one static time-step where inertia and gravity terms are neglected.
A linear equation under-relaxation factor of 0.999 is used to ensure convergence; the use of equation (as opposed to field) under-relaxation has been found to be important in contact analyses.

\begin{figure}[htb]
	\centering
	\subfigure[Axial displacement]
	{
		\includegraphics[width=0.48\textwidth]{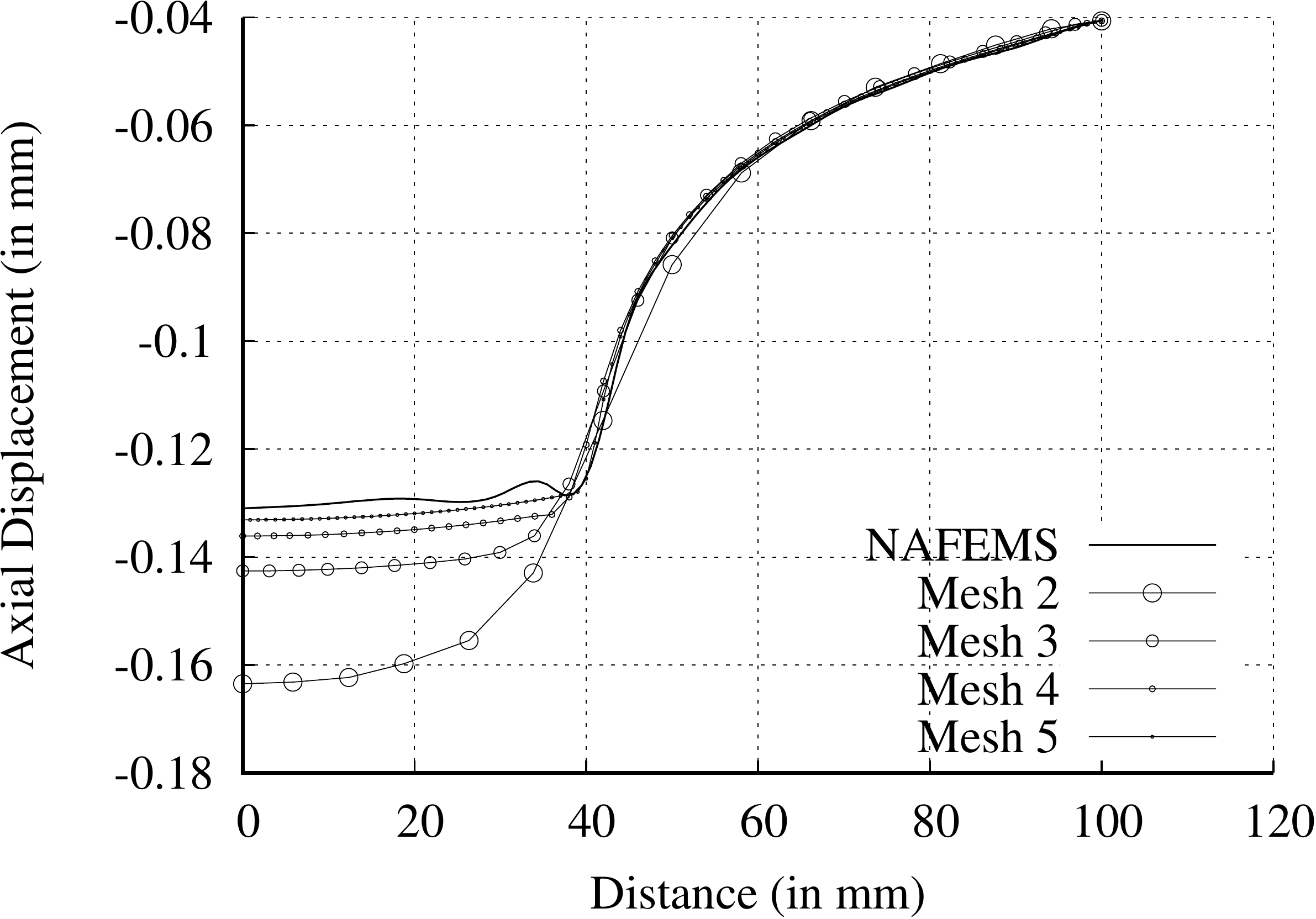}
		\label{fig:punch_axialDisp}
	}
	\subfigure[Radial displacement]
	{
		\includegraphics[width=0.48\textwidth]{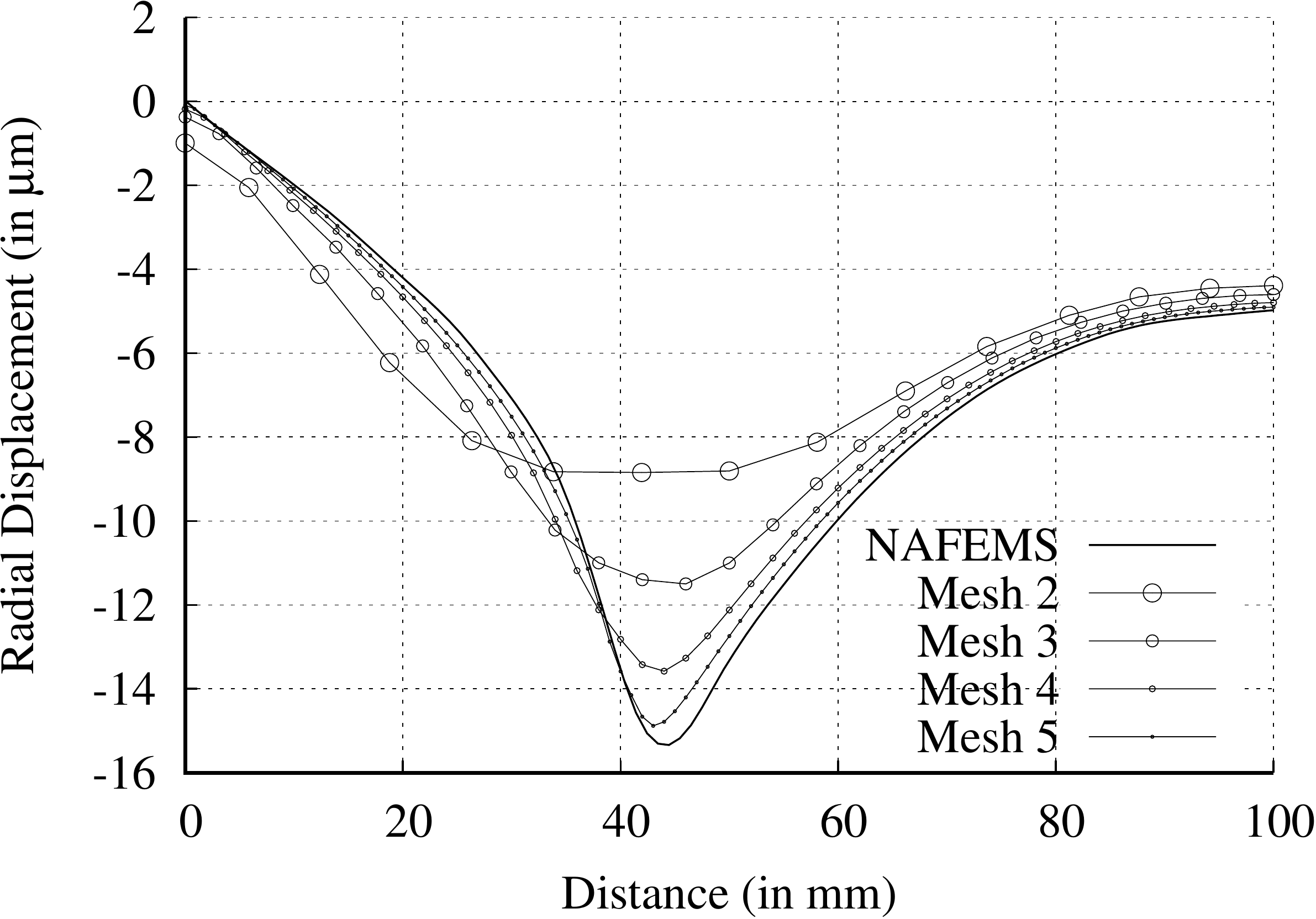}
		\label{fig:punch_radialDisp}
	}
	\caption{A 3-D Elastic Punch: displacement predictions on the cylinder surface}
\end{figure}

\begin{figure}[htb]
	\centering
	\subfigure[Von Mises stress distribution]
	{
		\includegraphics[width=0.48\textwidth]{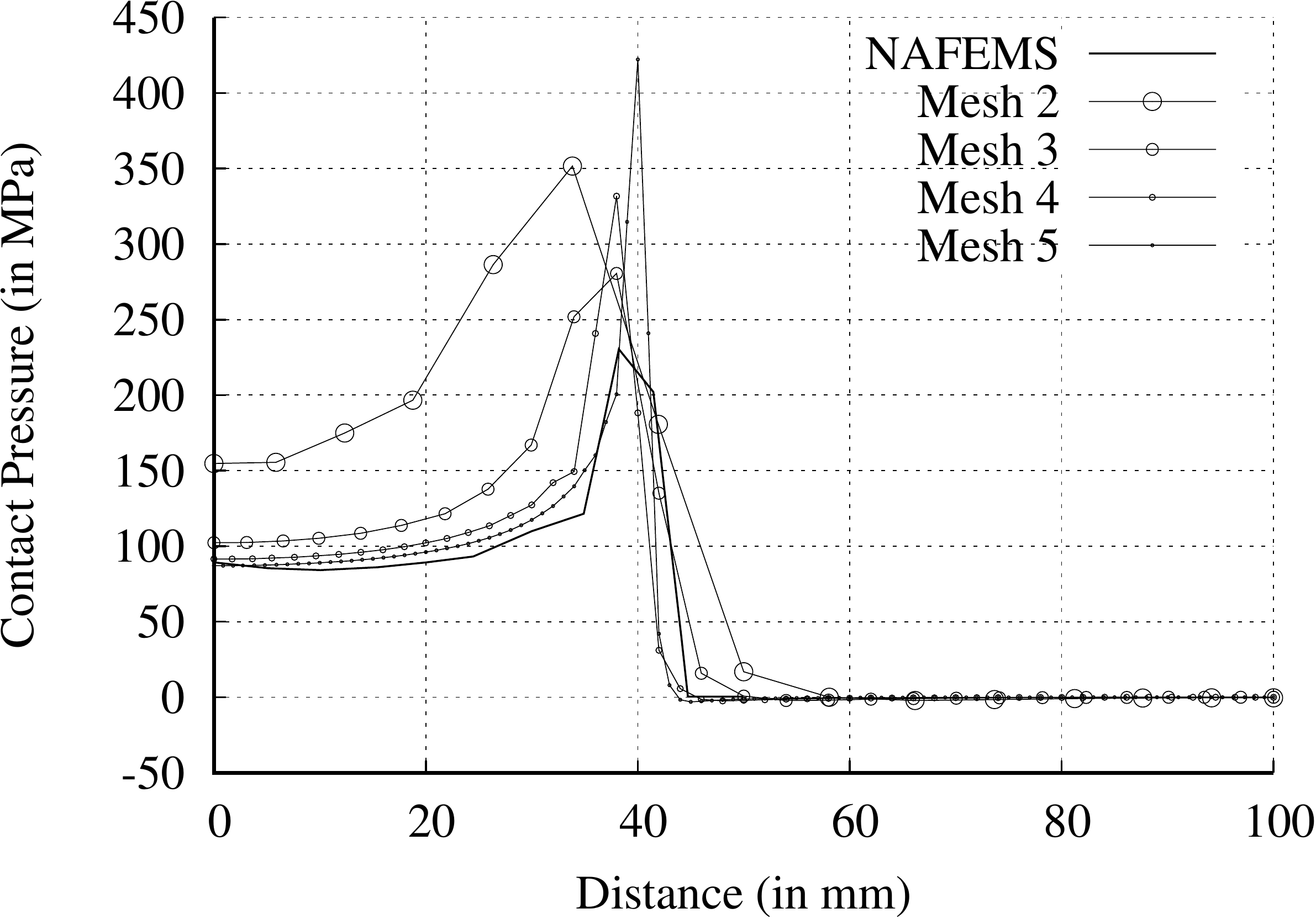}
		\label{fig:punch_contactPressure}
	}
	\subfigure[Reduction of discretisation error in the axial displacement at point A with increasing mesh refinement]
	{
		\includegraphics[width=0.48\textwidth]{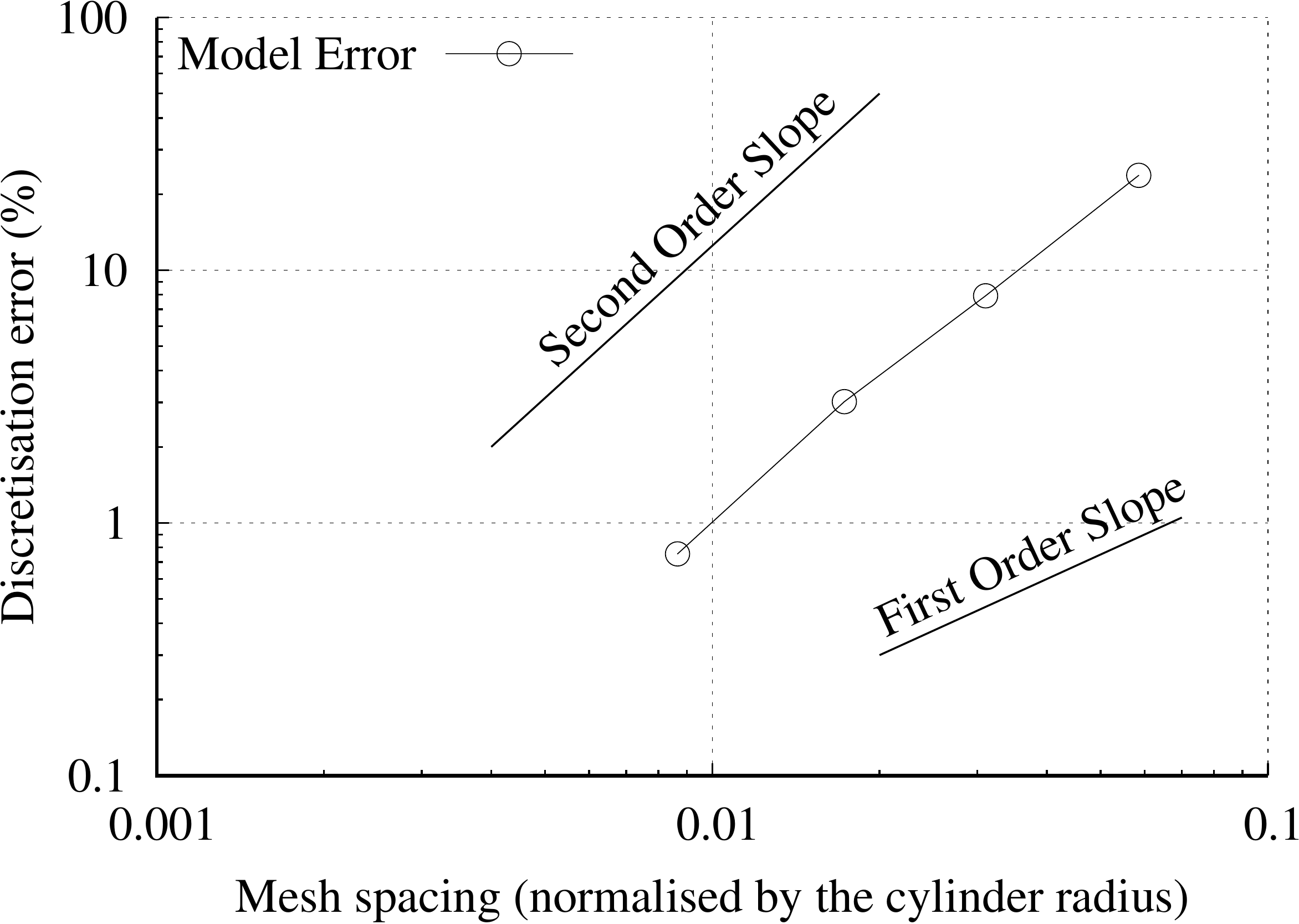}
		\label{fig:punchError}
	}
	\caption{A 3-D Elastic Punch: contact pressure stress distribution on the cylinder surface and the reduction in discretisation error}
\end{figure}

The predicted displacements and stresses are examined on the cylinder surface along a radial line from point A to point B, as shown in Figure \ref{fig:punchGeom};
the predicted axial displacement distribution is shown in Figure \ref{fig:punch_axialDisp}, where the results for the four meshes of increasing refinement are compared with predictions from the NAFEMS benchmark \citep{NAFEMS:contact};
the reference plots have been digitised using the {WebPlotDigitizer} software \citep{WebPlotDigitizer}.
Figure \ref{fig:punch_radialDisp} shows the predicted radial displacement distribution, and Figure \ref{fig:punch_contactPressure} shows the contact pressure distribution.
Both displacement and stress predictions can be seen to approach the reference solutions as the mesh is refined.

To examine discretisation error, the axial displacement at point A on the surface of the cylinder is examined for the four meshes;
Figure \ref{fig:punchError} shows the discretisation error to reduce from over 20\% to less than 1\% as the mesh is refined, where the reference displacement is approximated using Richardson's extrapolation \citep{Roache1997}.

\subsection{Stress relaxation of a viscoelastic tube}
This case consists of a tube constructed from a viscoelastic material, where the inner surface is quickly displaced in the radial direction; the tube relaxes towards steady-state after initially experiencing high wall stresses.
A similar case has been examined in the documentation of the commercial finite element software Comsol \citep{Comsol:documentation}.

\begin{figure}[htb]
	\centering
	\subfigure[Geometry]
	{
		\includegraphics[width=0.48\textwidth]{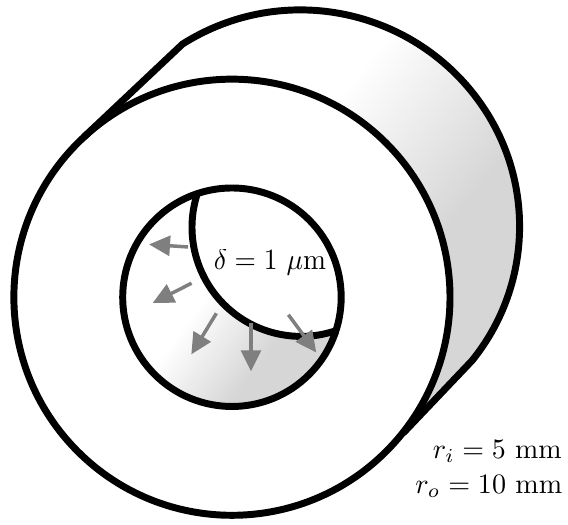}
		\label{fig:viscoTubeGeom}
	}
	\subfigure[2-D mesh containing 1~152 cells]
	{
		\includegraphics[width=0.48\textwidth]{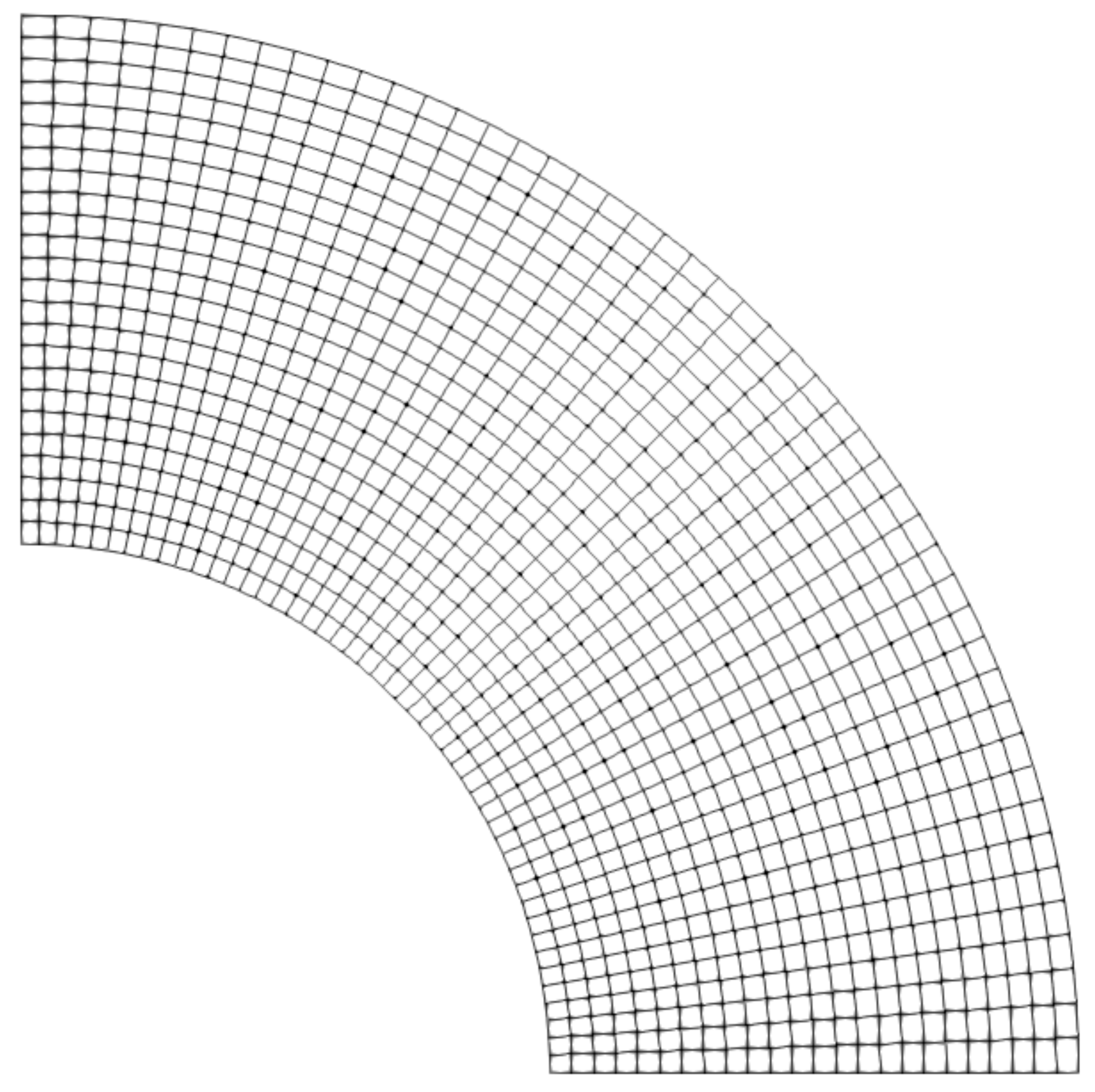}
		\label{fig:viscoTubeMesh}
	}
	\caption{Stress relaxation of a viscoelastic tube: problem geometry and mesh}
\end{figure}

The case demonstrates a transient viscoelastic analysis;
the case is modelled as 2-D with quarter symmetries, although the solution in this case is expected to be 1-D axisymmetric.
Figure \ref{fig:viscoTubeGeom} shows a schematic of the problem geometry, consisting of a long cylindrical tube with inner radius 5 \milli\meter\ and outer radius 10 \milli\meter.
A uniform quadrilateral mesh containing $1~152$ cells has been generated using OpenFOAM meshing utility \texttt{blockMesh}.

The tube material is assumed to be a viscoelastic polymer and its deviatoric response is represented by a Prony series, whereas its bulk/volumetric response is assumed to be elastic;
the mechanical properties are given in Table \ref{table:viscoTube_properties}.
\begin{table}[htb]
  \centering
	  \ra{1.3}
		\begin{tabular}{@{}lll@{}}
		\toprule
		\emph{Relaxed} Young's modulus & $E_{\infty}$ & 39.58 \giga\pascal \\
		Poisson's ratio & $\nu$ & 0.33 \\
		\emph{Maxwell Models} & & \\
		Young's modulus 1 & $E_1$ & 2.93 \giga\pascal \\
		Young's modulus 2 & $E_2$ & 5.86 \giga\pascal \\
		Young's modulus 3 & $E_3$ & 6.60 \giga\pascal \\
		Young's modulus 4 & $E_4$ & 18.32 \giga\pascal \\
		Relaxation time 1 & $\tau_1$ & 30 \second \\
		Relaxation time 2 & $\tau_2$ & 300 \second \\
		Relaxation time 3 & $\tau_3$ & 3~000 \second \\
		Relaxation time 4 & $\tau_4$ & 12~000 \second \\
		\bottomrule
		\end{tabular}
\caption{Stress relaxation of a viscoelastic tube: viscoelastic mechanical properties given in terms of a Prony series}
\label{table:viscoTube_properties}
\end{table}

The inner surface of the tube is displaced 1 \micro\meter\ in the radial direction; in order to mimic the displacement occurring instantaneously, it is applied within the first time step.
The outer surface of the tube is traction free and two symmetry conditions are applied, as indicated in Figure \ref{fig:viscoTubeGeom}.

The case is solved over a period of 7~000 \second, where inertia and gravity terms are neglected.
Four separate time step sizes are examined (700, 350, 175, 87.5 \second) to illustrate the effect of loading increment size on the employed material law integration, where  the current method employs the recursive implementation presented by \citet{Simo1998}.

\begin{figure}[htb]
	\centering
	\subfigure[Radial stress]
	{
		\includegraphics[width=0.48\textwidth]{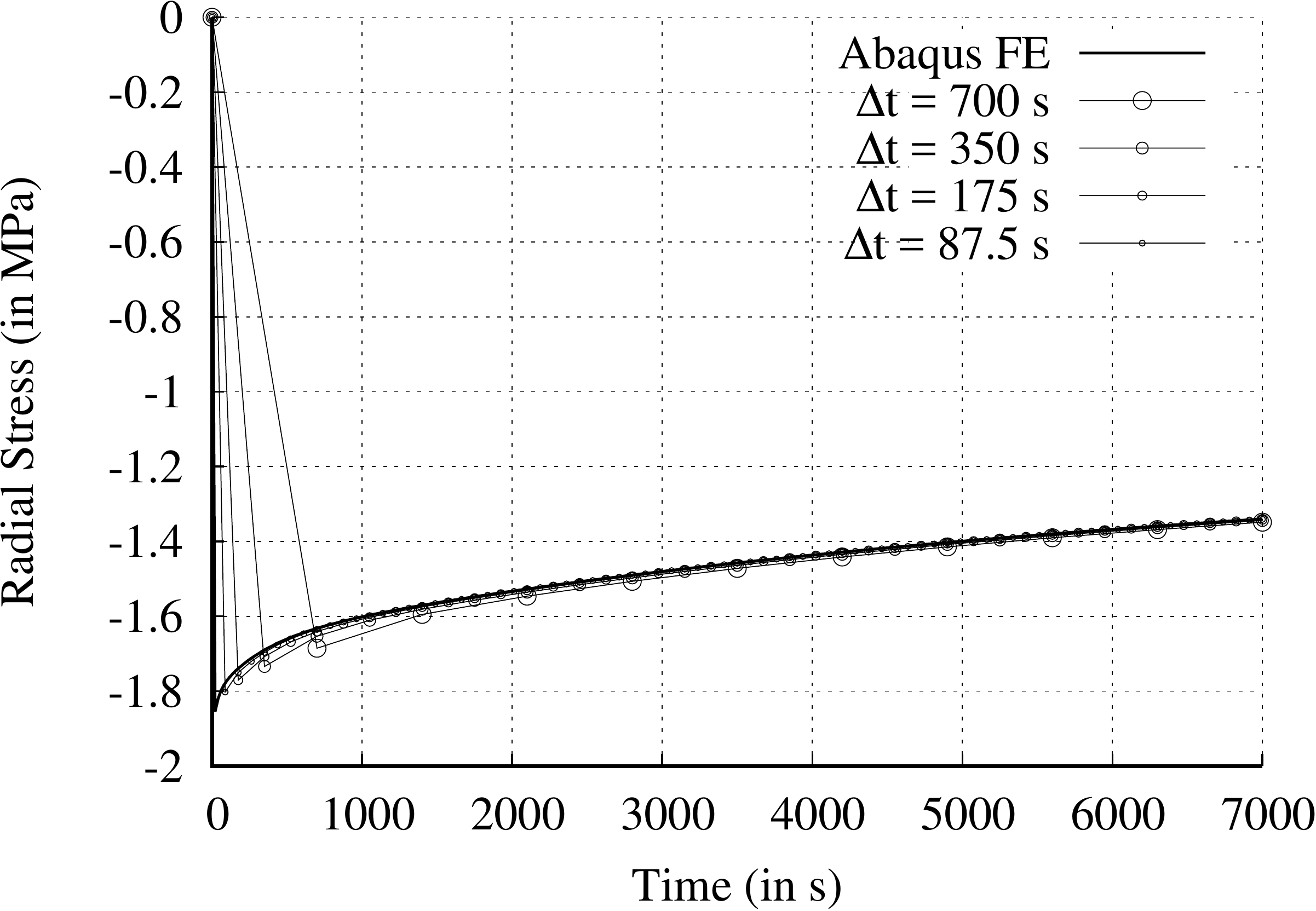}
		\label{fig:viscoTube_radialStress}
	}
	\subfigure[Radial displacement]
	{
		\includegraphics[width=0.48\textwidth]{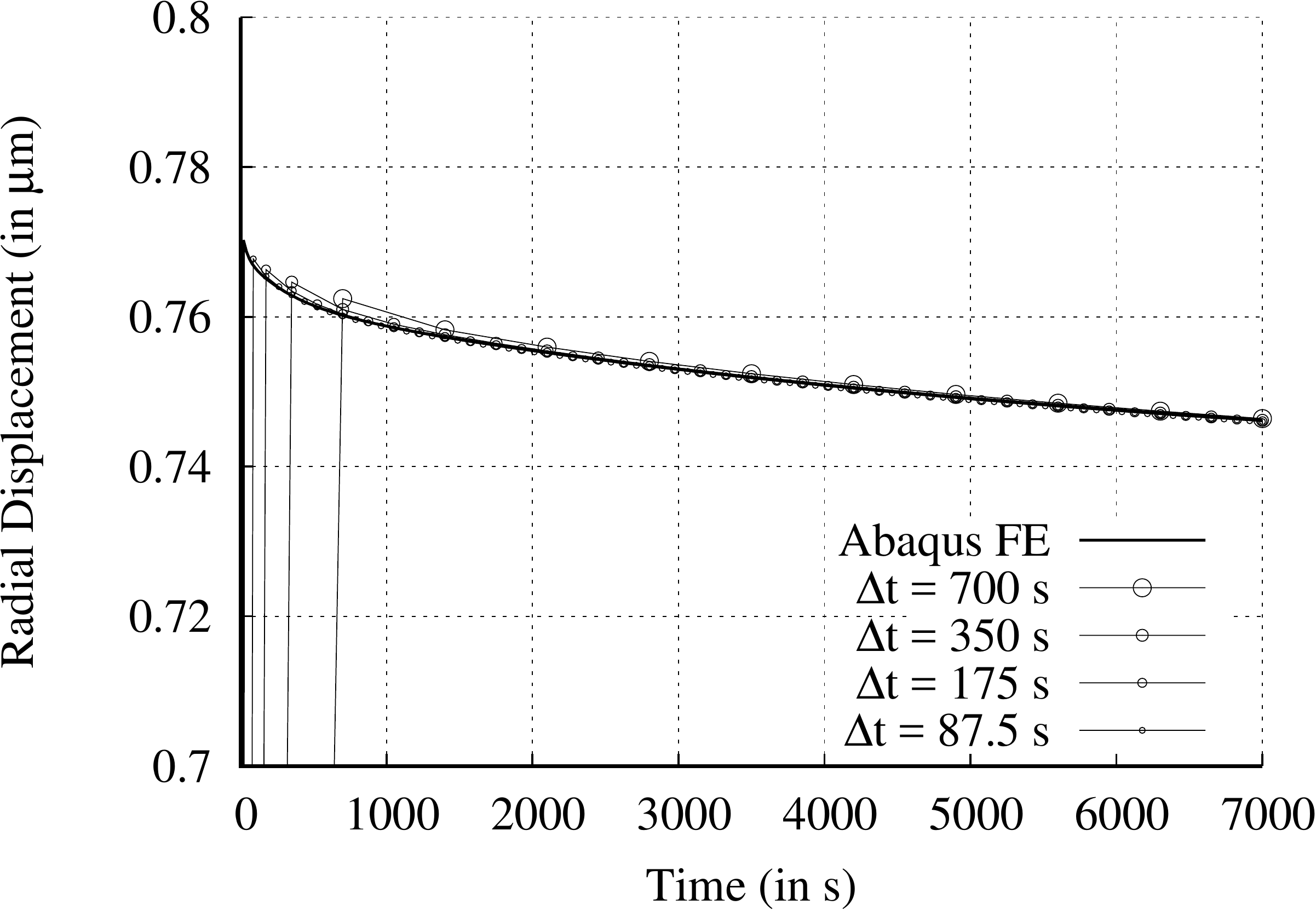}
		\label{fig:viscoTube_radialDisp}
	}
	\caption{Stress relaxation of a viscoelastic tube: radial stress and radial displacement predictions across the wall thickness, showing the effect of loading step size}
\end{figure}

The predicted radial stress on the inner tube surface versus time is examined;
the results for the four different time step sizes are shown in Figure \ref{fig:viscoTube_radialStress}, where results from FE software Abaqus are given for comparison;
the Abaqus results have been generated using a time step size of 20 \second\ to ensure the material integration error is small.
Figure \ref{fig:viscoTube_radialDisp} shows the predicted radial displacement on the inner tube versus time, where the \emph{y}-axis (displacement) scale does not start from zero to allow the differences to be more clearly seen.

It can be seen that the predicted stresses and displacements are relatively insensitive to the time-step, and as the step size is reduced the predictions are seen to approach the reference solution.


\subsection{Perforated Elastoplastic Plate}
\begin{figure}[htb]
	\centering
	\subfigure[Geometry]
	{
		\includegraphics[width=0.51\textwidth]{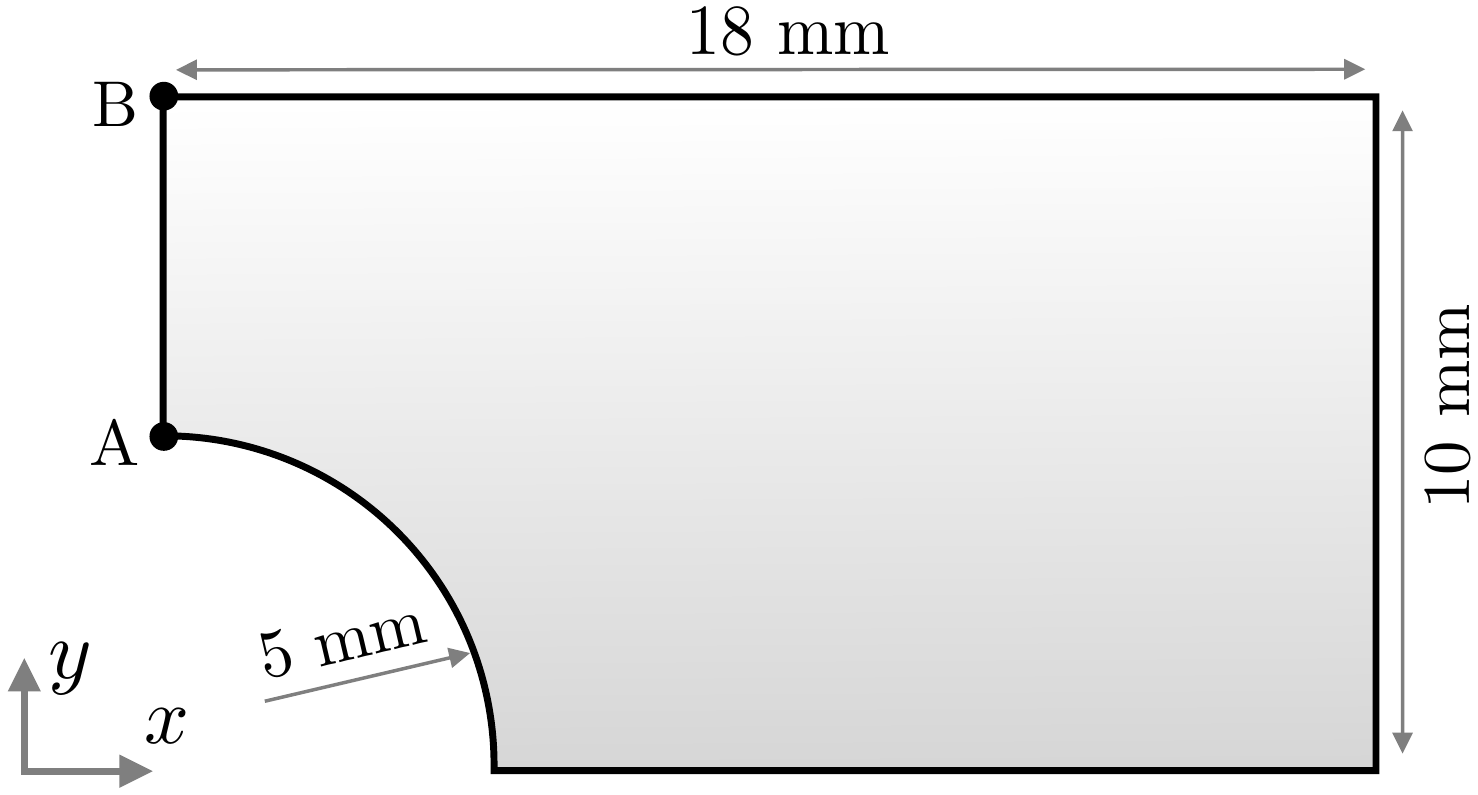}
		\label{fig:plateHoleGeom}
	}
	\subfigure[2-D mesh containing 300 cells]
	{
		\includegraphics[width=0.45\textwidth]{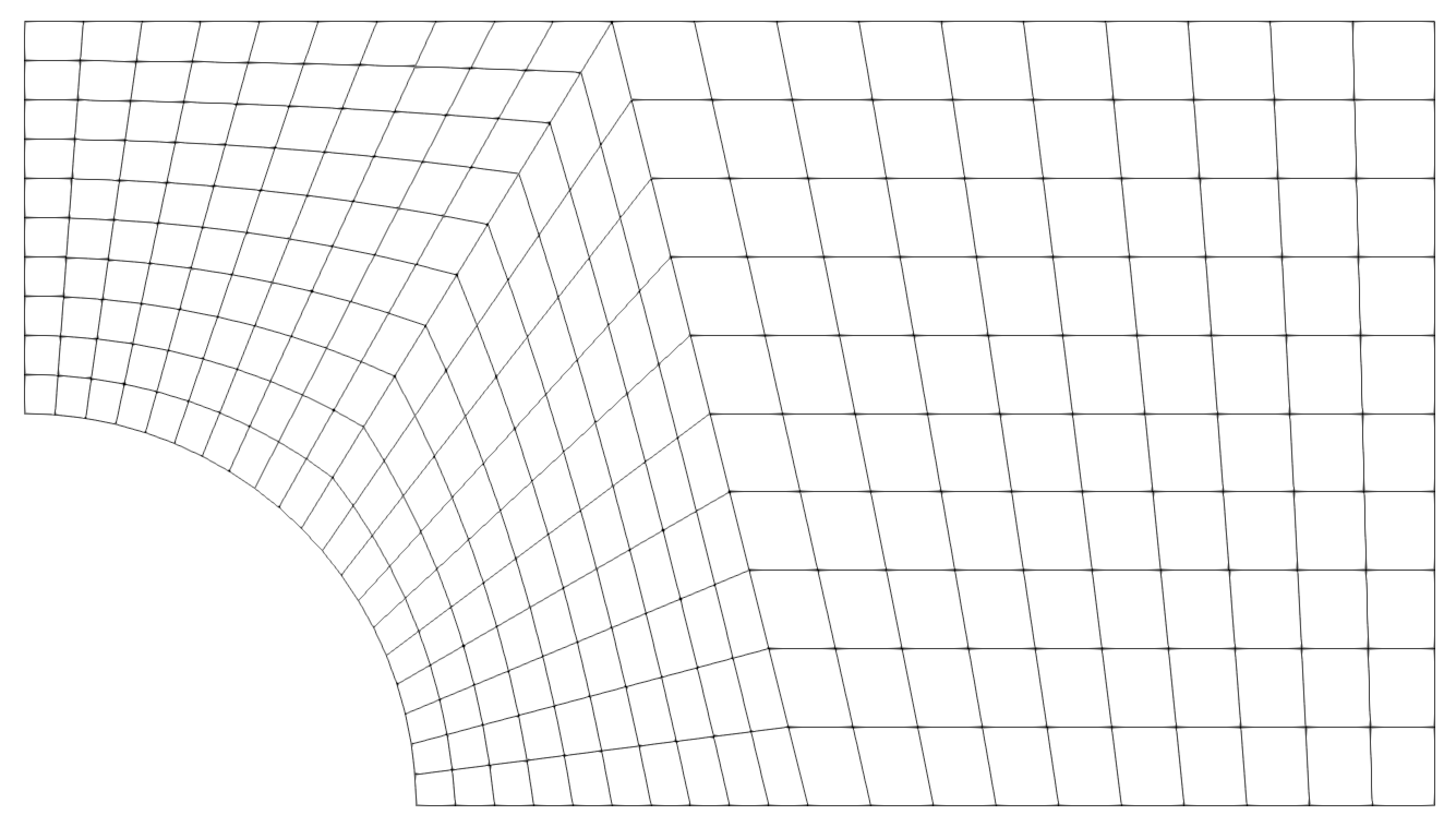}
		\label{fig:plateHoleMesh}
	}
	\caption{Perforated Elastoplastic Plate: problem geometry and mesh}
\end{figure}

A perforated plate subjected to tension has been examined a number of times in different forms, for example, \citep{Zienkiewicz2000, Comsol:documentation, Jasak2000:linearElasticity, Demirdzic1997}.
This classic case consists of a thin plate with a circular hole subjected to tension, where in this case, plastic deformation is considered.
Owing to the dual symmetries of the problem, one quarter of the plate is modelled and is represented here as plane strain 2-D.
The plate has dimensions of $20 \times 36$ \milli\meter\ with a hole of diameter 10 \milli\meter, shown schematically in Figure \ref{fig:plateHoleGeom}.
Four successively refined structured quadrilateral meshes containing cells 300, 1~200, 4~800 and 19~200 cells have been created using OpenFOAM meshing utility \texttt{blockMesh}; the coarsest mesh, containing 300 cells, is shown in in Figure \ref{fig:plateHoleMesh}.

The aluminium plate material is assumed to be isotropic Hookean elastic with von J$_2$/Mises yield criterion and isotropic hardening; the mechanical properties are given in Table \ref{table:plateHole_properties}.
\begin{table}[htb]
  \centering
	  \ra{1.3}
		\begin{tabular}{@{}lll@{}}
		\toprule
		Young's modulus & $E$ & 70 \giga\pascal \\
		Poisson's ratio & $\nu$ & 0.3 \\
		Initial yield stress & $\sigma_{Y0}$ & 243 \mega\pascal \\
		Isotropic hardening modulus & $H$ & 2.25 \giga\pascal \\
		\bottomrule
		\end{tabular}
\caption{Perforated Elastoplastic Plate: elastoplastic mechanical properties}
\label{table:plateHole_properties}
\end{table}

The right boundary of the plate is loaded with a time-varying normal traction, $T^n_{\text{right}}$, that linearly ramps up over the first 10 \second\ and linearly ramps down from 10 \second\ to 20 \second:
\begin{eqnarray}
	T^n_{\text{right}} = 
		\begin{cases}
			13.365 t \text{ Pa}		&	\text{for}\quad  0<t \leq 10 \text{ \second} \\ 
			13.365 (20 - t) \text{ Pa}	&	\text{for}\quad  10<t \leq 20 \text{ \second}
		\end{cases}
\end{eqnarray}
where $t$ is the current time in seconds.
Following the approach in \citep{Comsol:documentation}, the peak value, occurring at $t = 10$ \second, has been selected so that the mean stress over the section through the hole is 10\% above the initial yield stress \ie $=1.1 \times 243 \text{ Pa} \times \nicefrac{(20 \text{ mm} - 10 \text{ mm})}{20\text{ mm}}$.
The shear traction on the right boundary is zero.
The hole and upper plate boundaries are traction free; symmetry conditions are applied to the left and bottom plate boundaries;
2-D plane strain conditions are assumed.
The case is solved over 20 \second\ in time steps of 0.1 \second, where inertia and gravity are neglected.

\begin{figure}[htb]
	\centering
	\subfigure[$y$ component of displacement at the upper rim of the hole]
	{
		\includegraphics[width=0.48\textwidth]{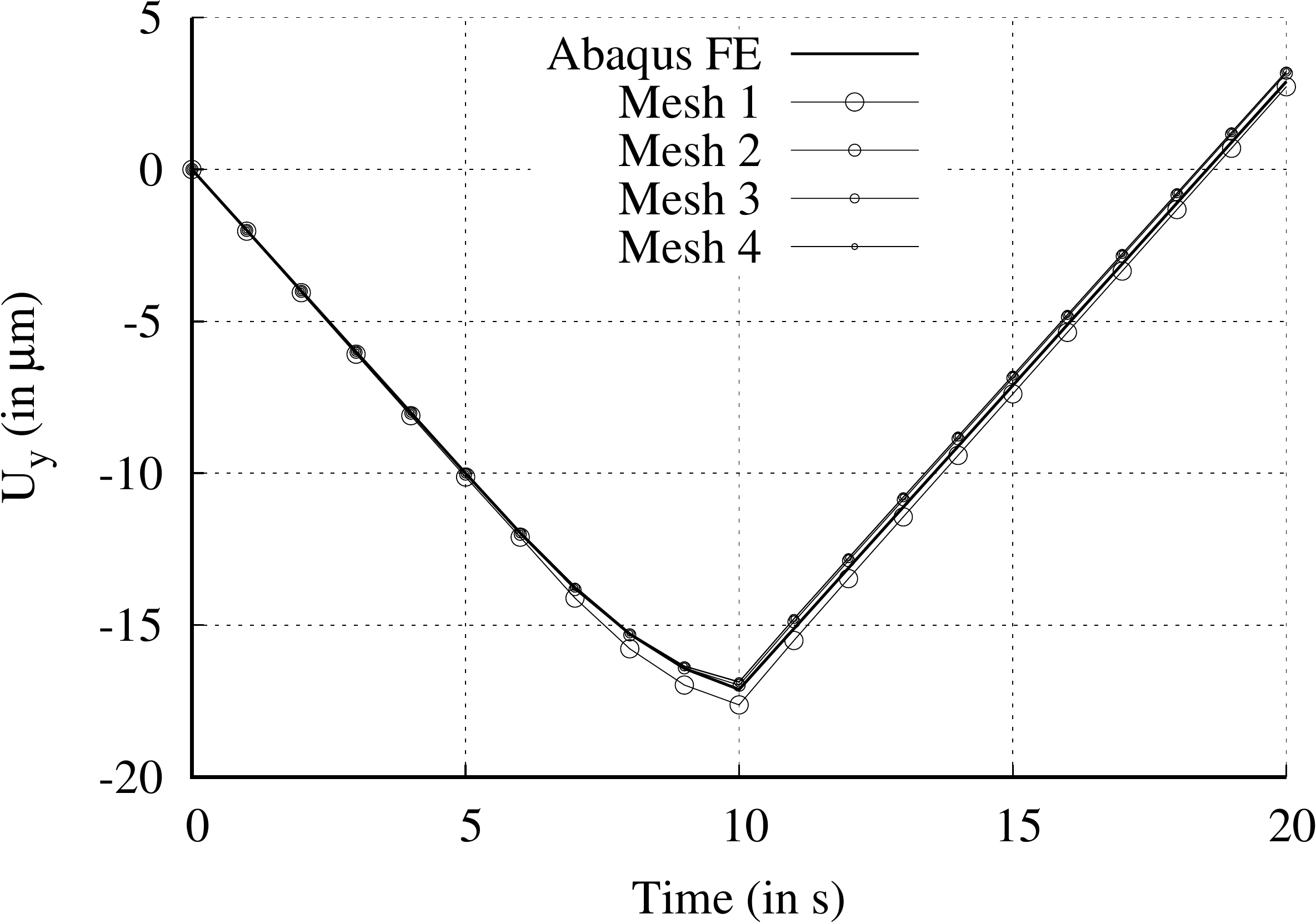}
		\label{fig:plateHolePointDisp}
	}
	\subfigure[Residual stress, $\boldsymbol{\sigma}_{xx}$, along the path A--B at $t = 20$ \second]
	{
		\includegraphics[width=0.48\textwidth]{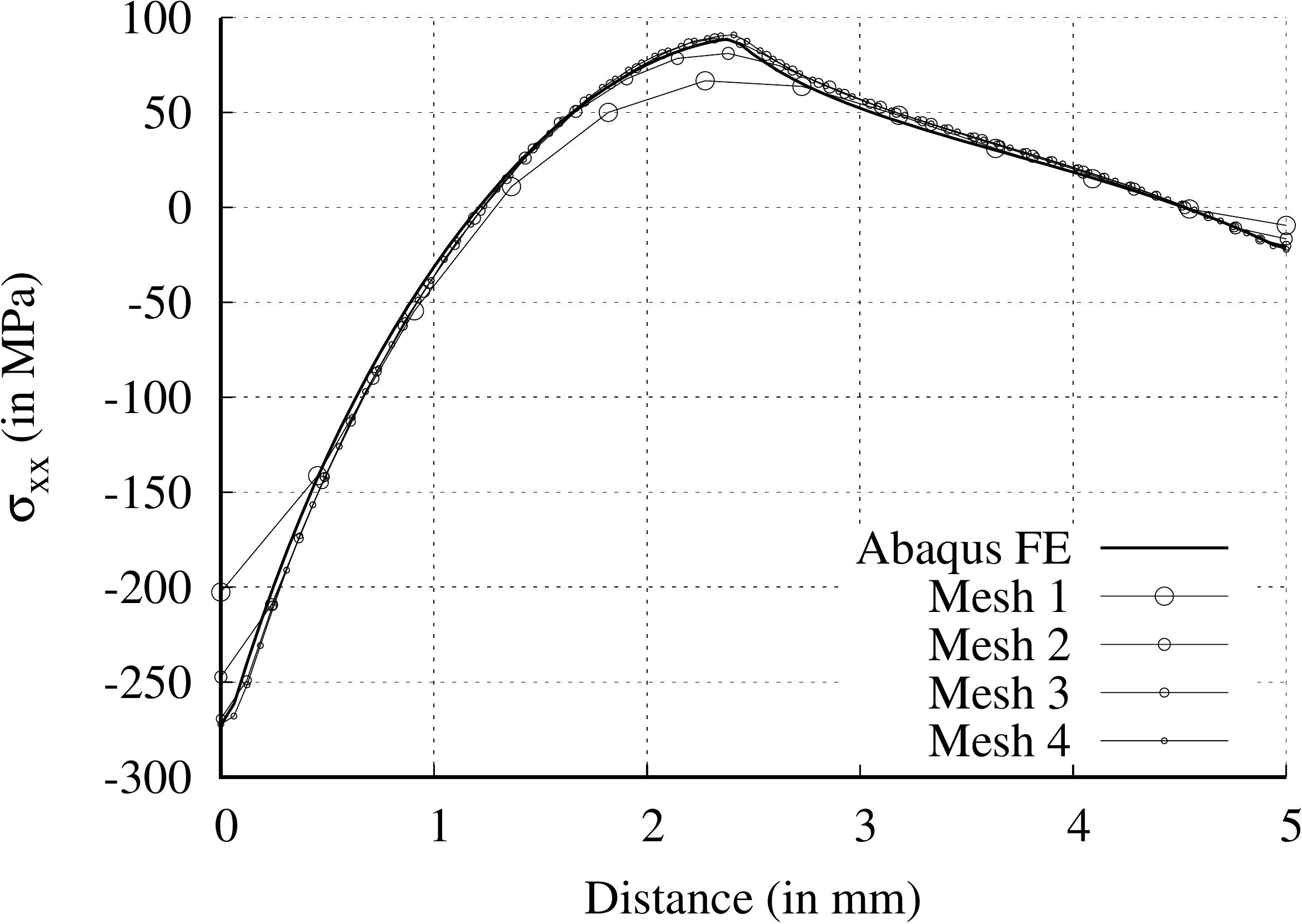}
		\label{fig:plateHoleStressTime20}
	}
	\caption{Perforated Elastoplastic Plate: displacement and residual stresses predictions}
\end{figure}

\begin{figure}[htb]
	\centering
	\subfigure[OpenFOAM: $t = 7$ \second]
	{
		\includegraphics[height=0.2\textwidth]{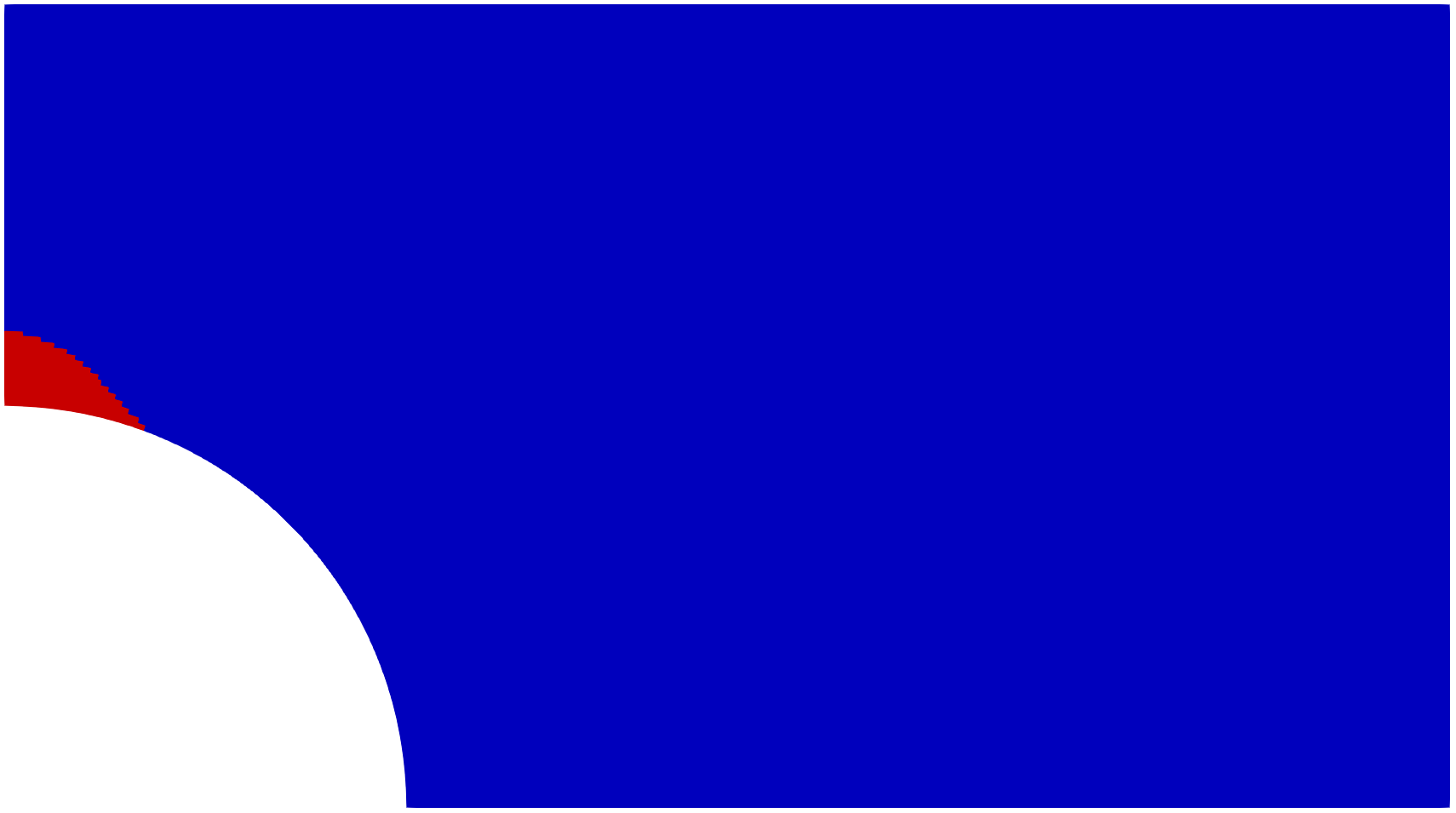}
	}
	\subfigure[Abaqus FE: $t = 7$ \second]
	{
		\includegraphics[height=0.2\textwidth, width=0.358\textwidth]{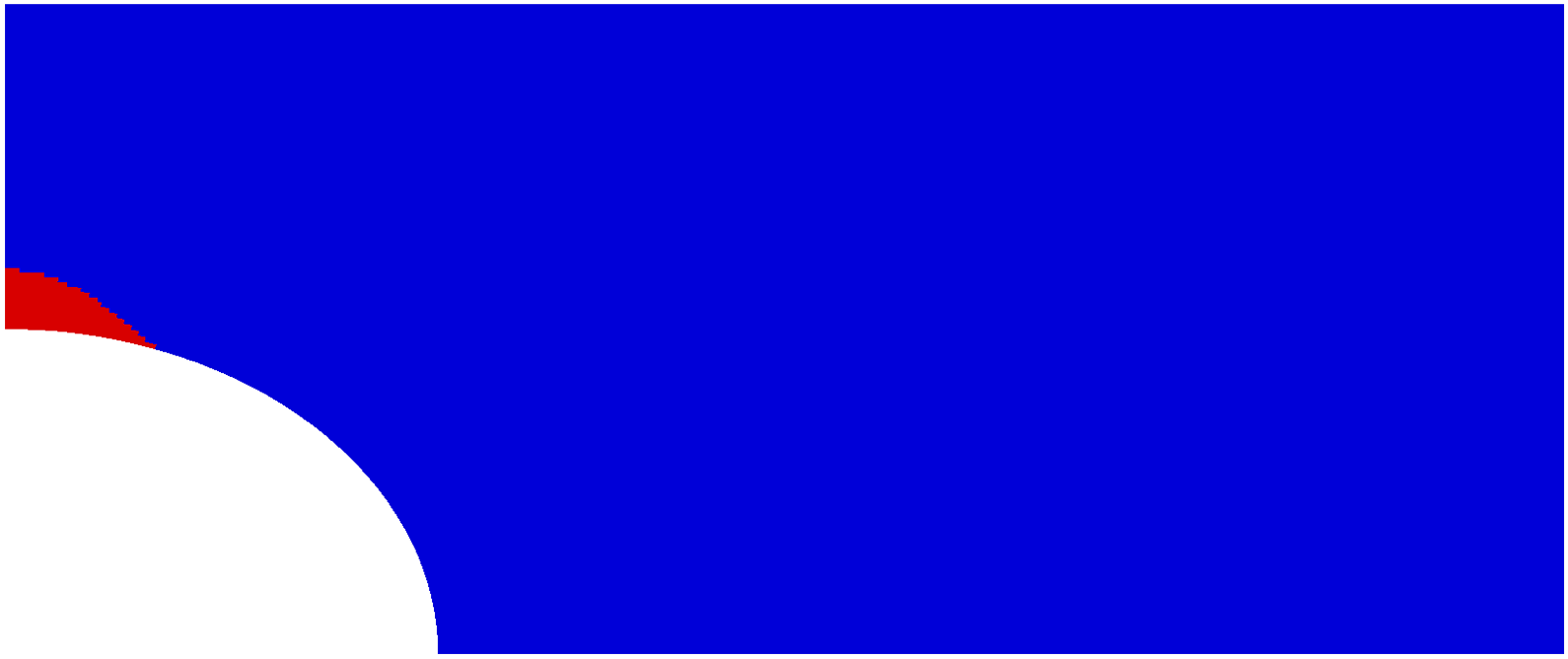}
	}
	\subfigure[OpenFOAM: $t = 8$ \second]
	{
		\includegraphics[height=0.2\textwidth]{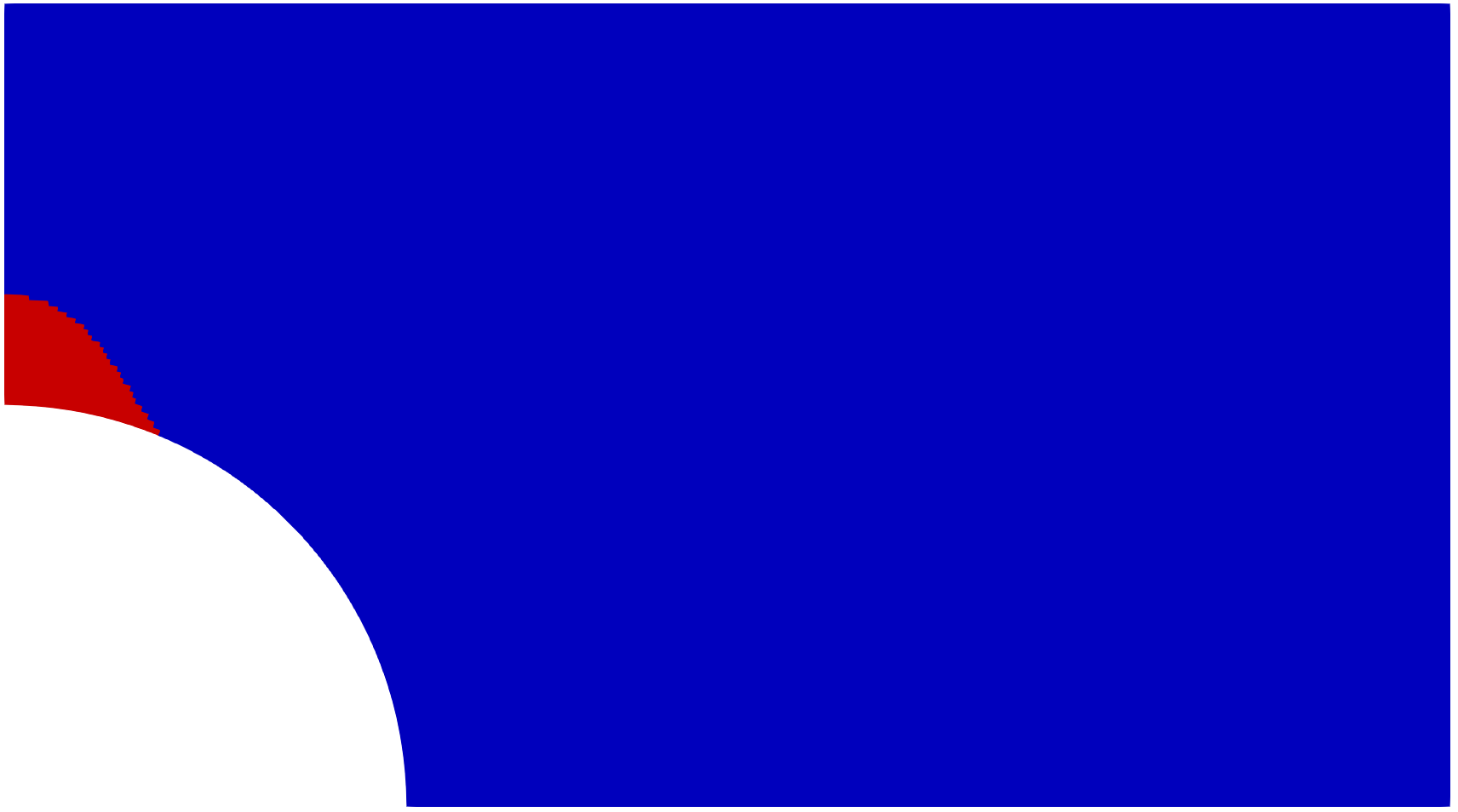}
	}
	\subfigure[Abaqus FE: $t = 8$ \second]
	{
		\includegraphics[height=0.2\textwidth, width=0.358\textwidth]{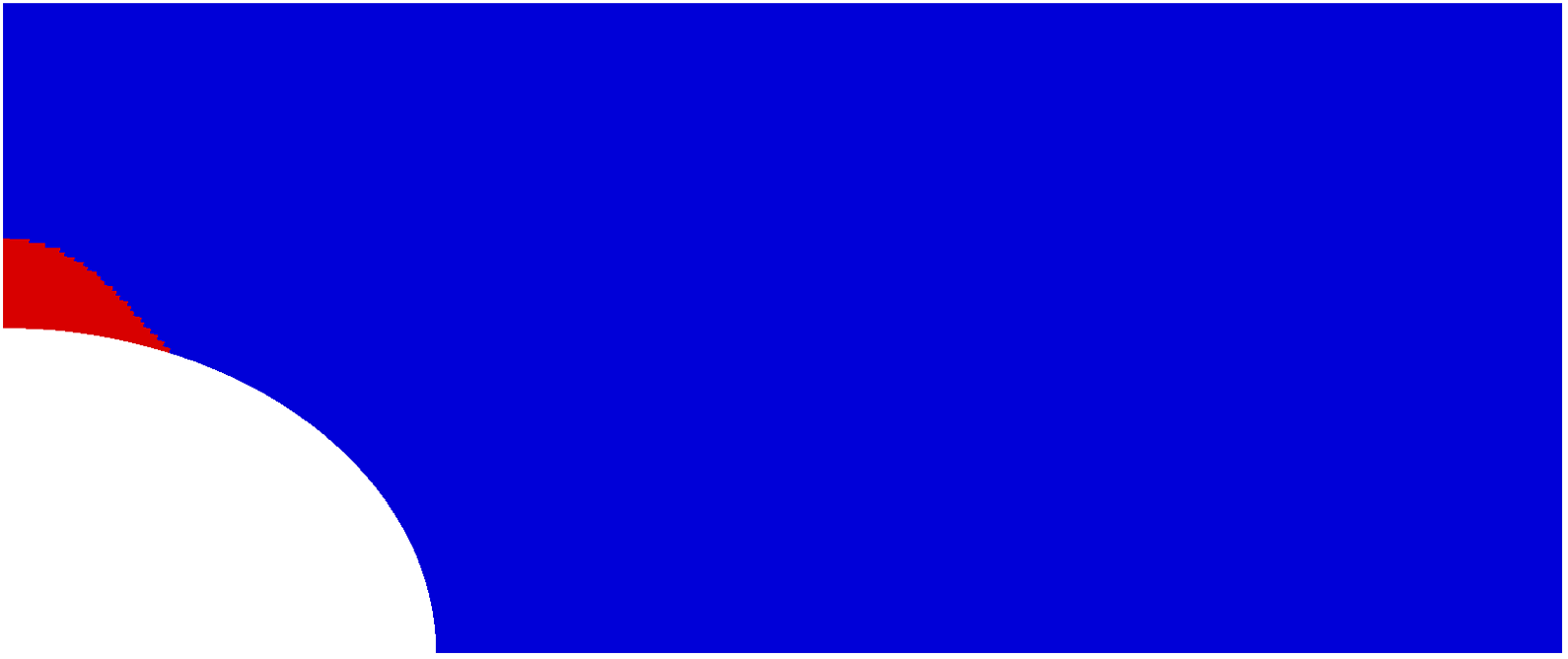}
	}
	\subfigure[OpenFOAM: $t = 9$ \second]
	{
		\includegraphics[height=0.2\textwidth]{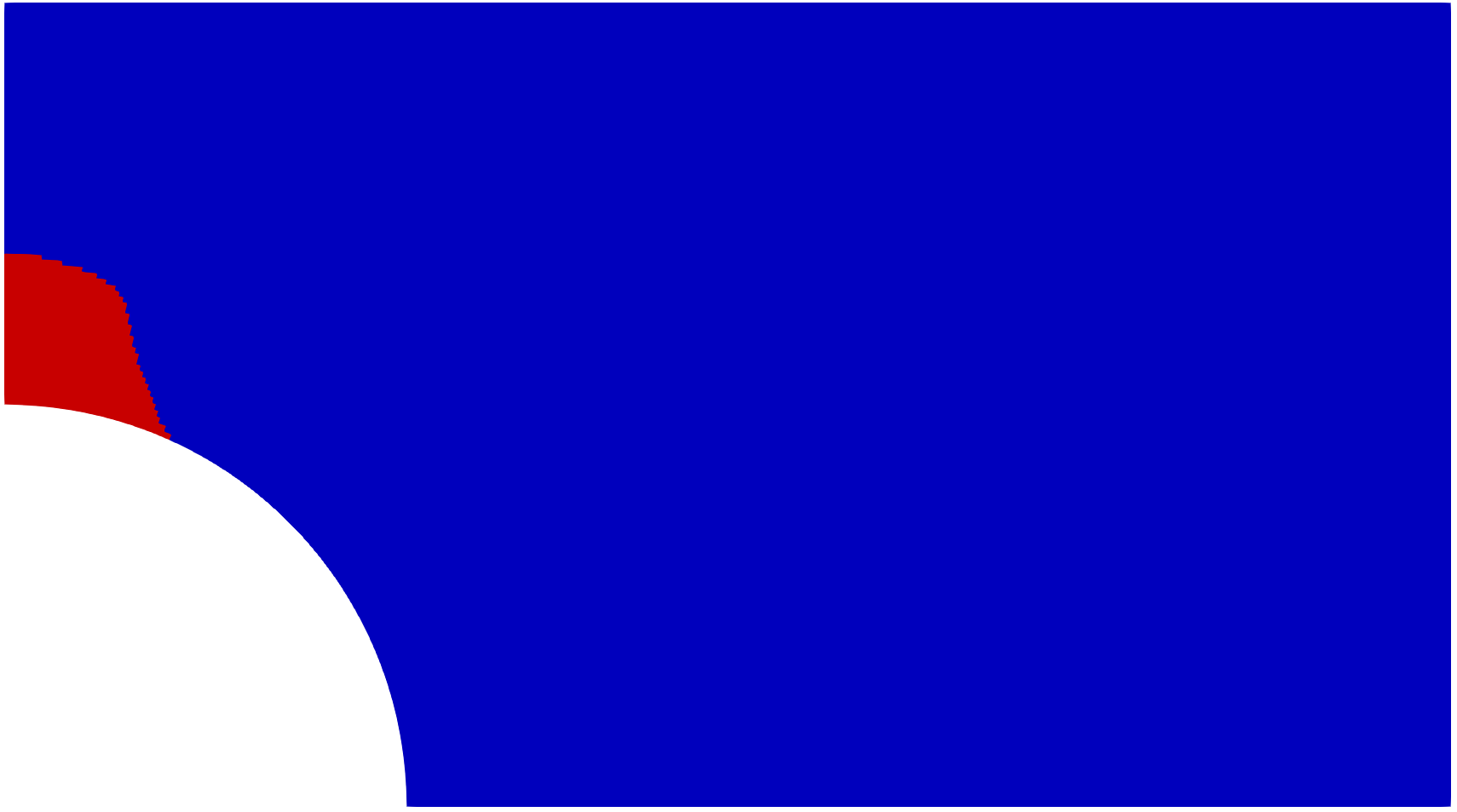}
	}
	\subfigure[Abaqus FE: $t = 9$ \second]
	{
		\includegraphics[height=0.2\textwidth, width=0.358\textwidth]{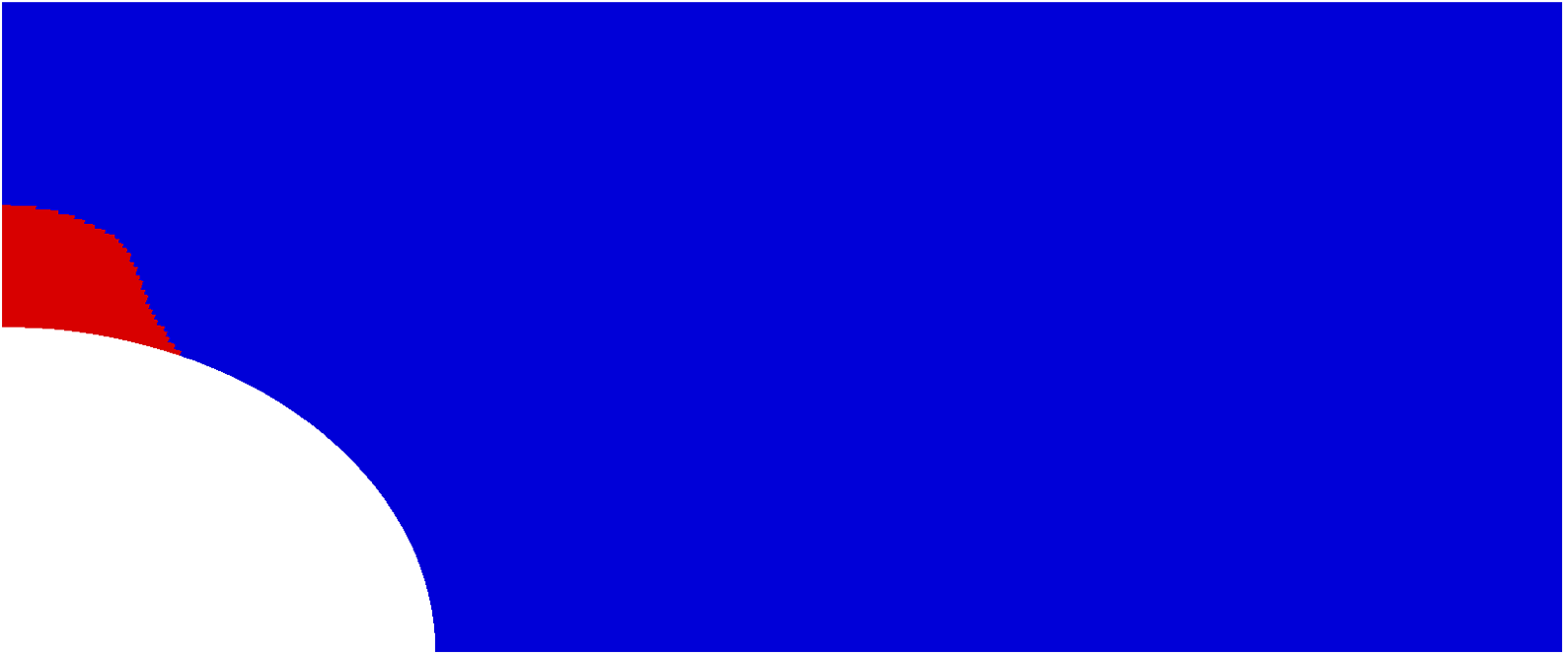}
	}
	\subfigure[OpenFOAM: $t = 10$ \second]
	{
		\includegraphics[height=0.2\textwidth]{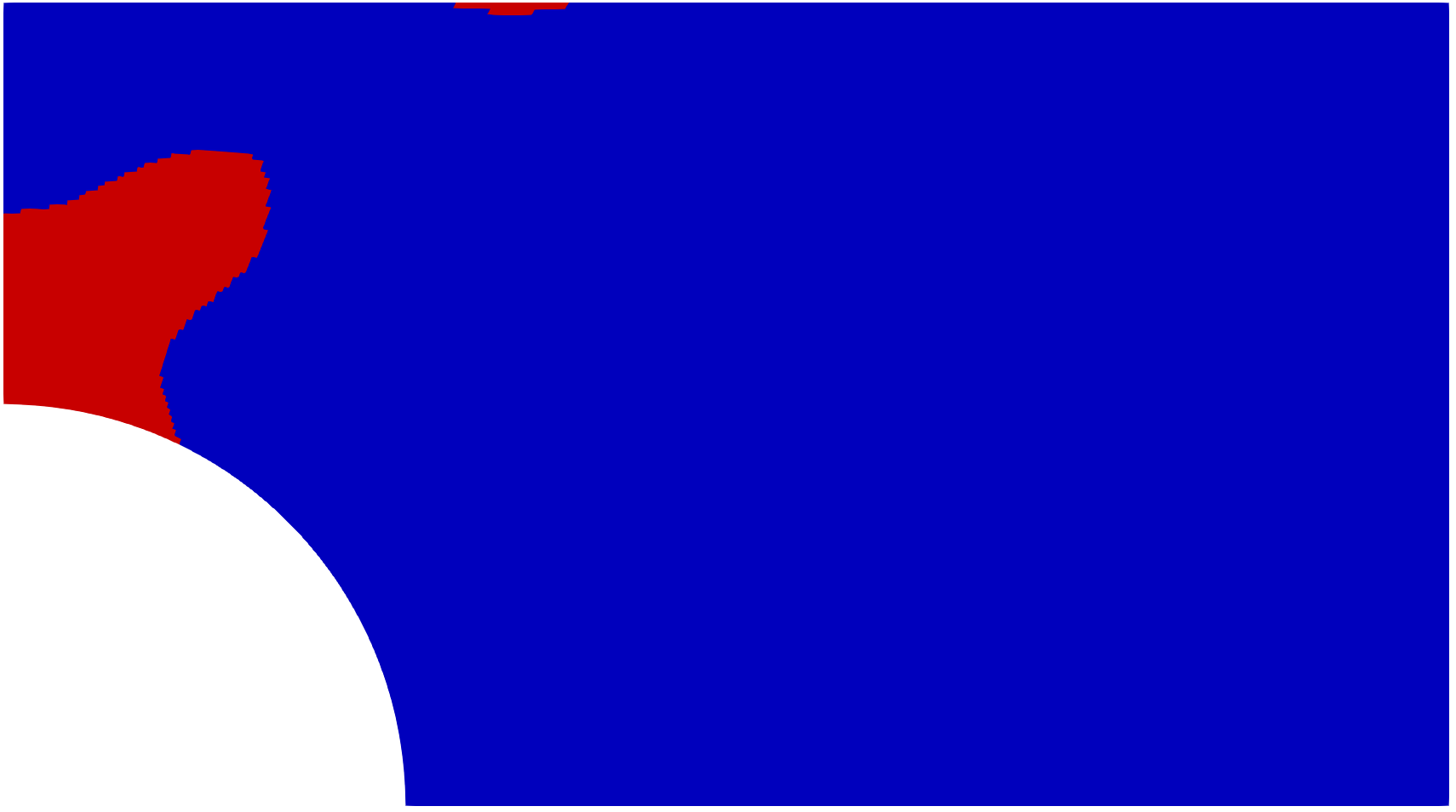}
	}
	\subfigure[Abaqus FE: $t = 10$ \second]
	{
		\includegraphics[height=0.2\textwidth, width=0.358\textwidth]{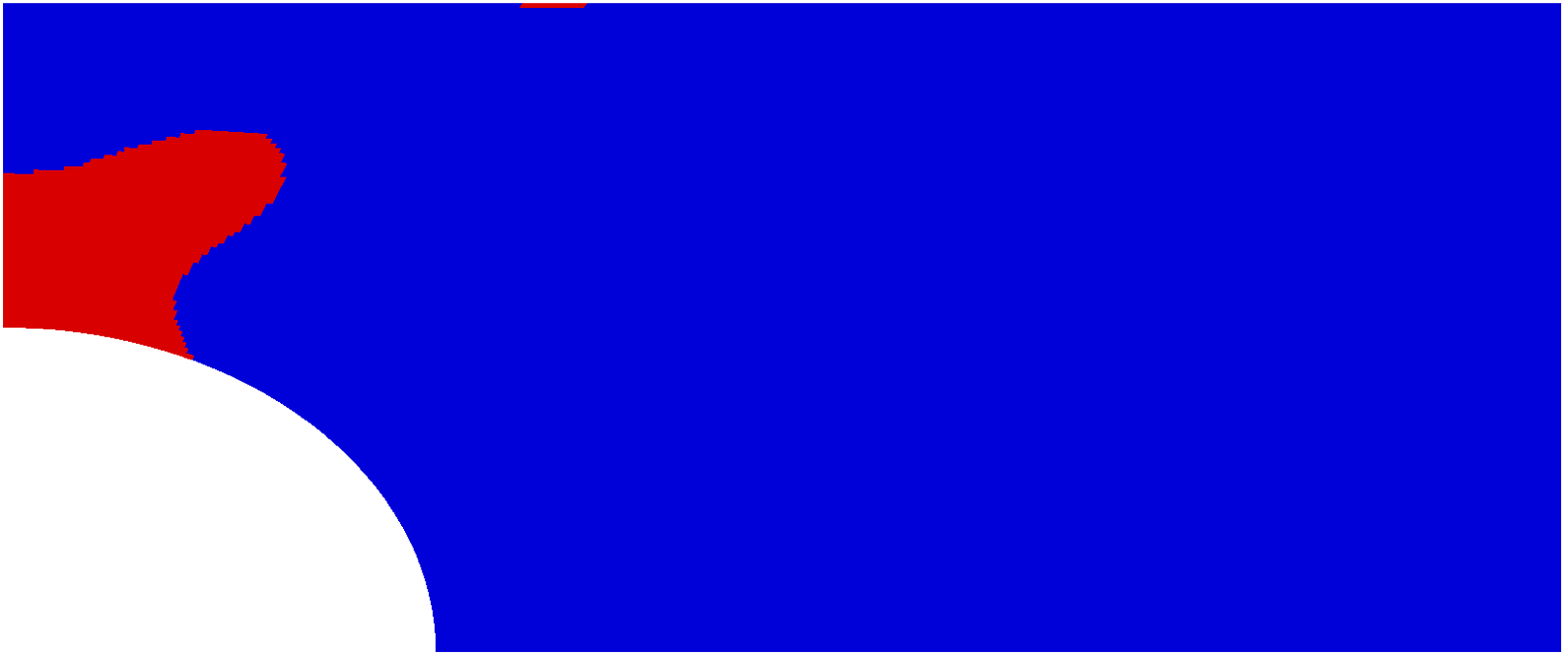}
	}
	\caption{Perforated Elastoplastic Plate: regions of predicted active yielding at succesive times}
	\label{fig:plateHoleActiveYield}
\end{figure}

The predicted $y$ component of displacement at point A on the top rim of the hole is shown in Figure \ref{fig:plateHolePointDisp}, where results from FE software Abaqus on the finest mesh are given for comparison.
As the mesh is refined, the predictions are seen to become mesh independent, in agreement with the FE results.
From Figure \ref{fig:plateHoleStressTime20}, the residual stress component $\sigma_{xx}$ can be seen along the path A--B at $t = 20$ \second.
As is expected, the region of material near the hole that previously underwent plastic deformation returns to a state of compression upon unloading.
In Figure \ref{fig:plateHoleActiveYield}, the regions of predicted active yielding are shown for four time instances;
the predictions can be seen to agree closely with the FE results.

\subsection{Micro-Beam in Channel Flow}
\begin{figure}[htb]
	\centering
	\subfigure[Geometry]
	{
		\includegraphics[width=0.95\textwidth]{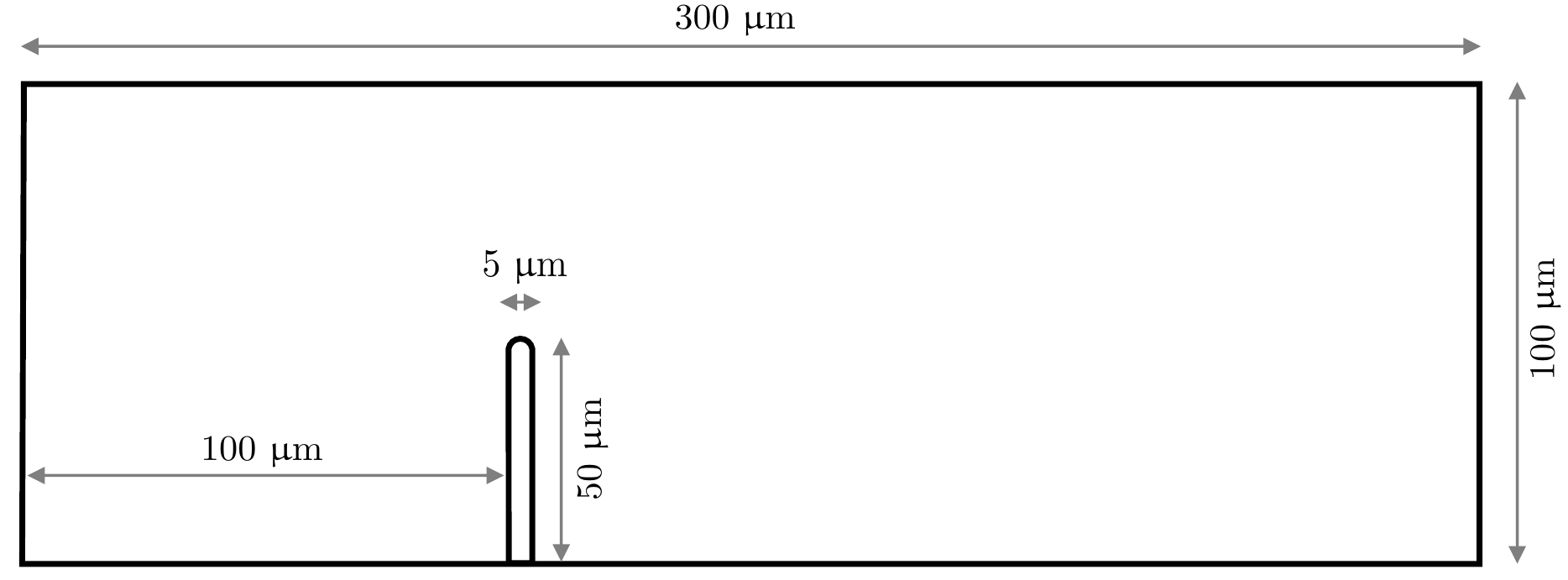}
		\label{fig:beamInFlowGeom}
	} \\
	\subfigure[2-D mesh containing 300 cells]
	{
		\includegraphics[width=0.9\textwidth]{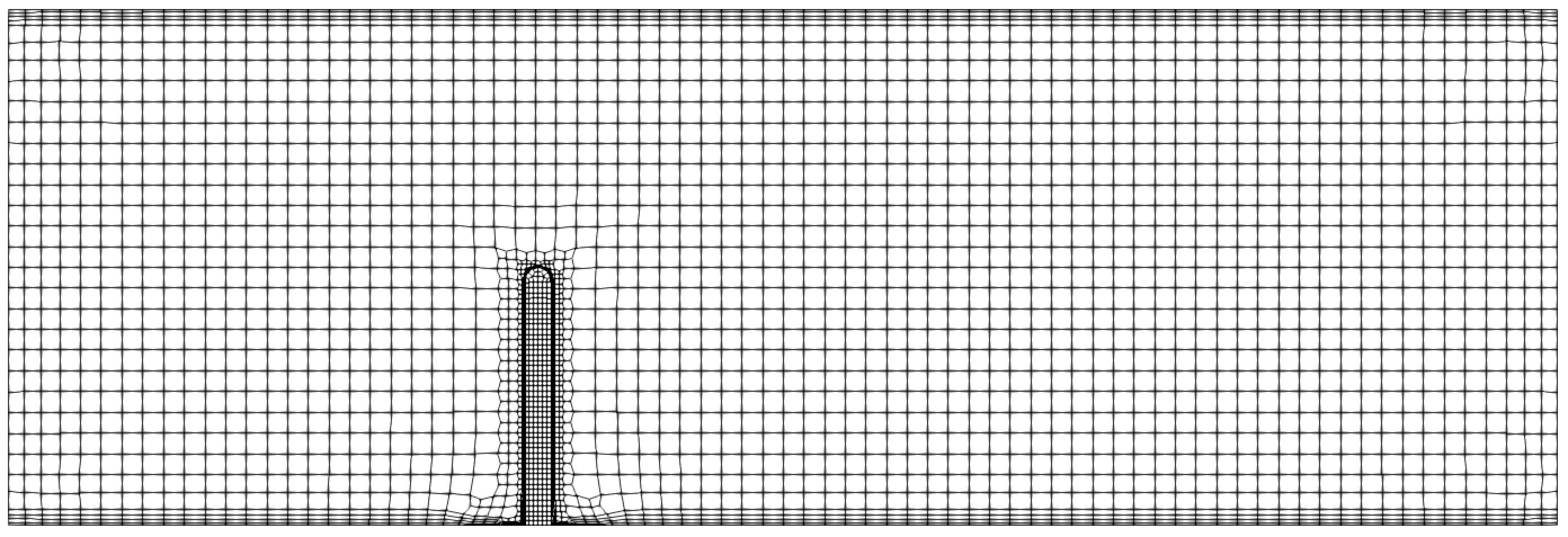}
		\label{fig:beamInFlowMesh}
	}
	\caption{Micro-Beam in Channel Flow: problem geometry and mesh}
\end{figure}

This case consists of an elastic beam subjected to laminar flow in a channel, and is based on a problem presented in the documentation for the commercial FE software \citet{Comsol:documentation};
it is possible to find a number of similar FSI cases in literature, for example, \citep{Turek2006, Tukovic2017:fsi}.
The case geometry, shown schematically in Figure \ref{fig:beamInFlowGeom}, consists of a solid plate, $50$ \micro\meter\ tall and $5$ \micro\meter\ wide with a semi-circular top, located $100$ \micro\meter\ from the left inlet boundary of a channel, where the channel is $100$ \micro\meter\ high and $300$ \micro\meter\ long.
The geometry for both fluid and solid domain has been created using open-source software Gmsh \citep{Geuzaine2009}; the geometry was then exported as a facetted STL surface and subsequently meshed using the OpenFOAM Cartesian meshing utility \texttt{cfMesh} \citep{cfMesh}.
Figure \ref{fig:beamInFlowMesh} shows the 2-D Cartesian mesh.
The fluid domain mesh consists of 3~126 cells, and the solid domain mesh consists of 296 cells.

The physical and mechanical properties of the solid and fluid are given in Table \ref{table:beamInFlow_properties};
the solid plate is assumed to be a neo-Hookean flexible polymer, where the out-of-plane length is long and hence plane strain conditions are assumed;
the fluid is assumed to be incompressible, Newtonian and isothermal.
\begin{table}[htb]
  \centering
	  \ra{1.3}
		\begin{tabular}{@{}lll@{}}
		\toprule
		\emph{solid} & & \\
		Density & $\rho_{\text{solid}}$ & 7850 \kilo\gram\per\metre\cubed \\
		Young's modulus & $E$ & 700 \kilo\pascal \\
		Poisson's ratio & $\nu$ & 0.3 \\
		\midrule
		\emph{fluid} & & \\
		Density & $\rho_{\text{fluid}}$ & 1000 \kilo\gram\per\metre\cubed \\
		Dynamic viscosity & $\eta$ & 0.001 \pascal$\cdot$\second \\
		\bottomrule
		\end{tabular}
\caption{Micro-Beam in Channel Flow: physical and mechanical properties}
\label{table:beamInFlow_properties}
\end{table}

At the left boundary of the channel, a parabolic inlet velocity is prescribed, where the average flow velocity initially increases rapidly, reaching a peak at 0.215 \second\, before gradually slowing to a steady-state value of 0.05 \meter\per\second.
The maximum inlet velocity normal to the left boundary, $U$, is given as:
\begin{eqnarray}
U = \frac{U_{\text{steady}} t^2}{\sqrt{(0.04 - t^2)^2 + (0.1t)^2}}
\end{eqnarray}
where $t$ is the current time in seconds, and $U_{\text{steady}}$ is the maximum inlet velocity at steady-state \ie at $t \gg 1$ \second.
A zero normal gradient of pressure is prescribed at the inlet.
The velocity at the upper and lower walls of the channel is set to zero \ie no slip condition, with a zero normal gradient of pressure.
The gauge pressure at the outlet is set to zero, with the normal gradient of velocity also set to zero.
The bottom surface of the beam is fixed (zero displacement) and gravity is neglected.
At the interface between the fluid and solid regions, continuity of force/traction, displacement and velocity is enforced via Dirichlet-Neumann coupling using Aitken adaptive under-relaxation; more details are found in Tukovi\v{c} \etal \citep{Tukovic2014, Tukovic2007:fsi, Tukovic2017:fsi}.

The model is solved for 1 \second\ of simulated time, where the time-step is adaptively set to ensure a Courant number of approximately 20.0 in the fluid domain; for the current mesh, this results in an average time-step of approximately 0.2 \milli\second.
The fluid is assumed to be water and is solved using the so-called PIMPLE (transient SIMPLE) segregated algorithm and the solid model employs the segregated updated Lagrangian solution procedure and neo-Hookean elastic constitutive law.
The material properties given in Table \ref{table:beamInFlowDeformedProperties}.
\begin{table}[htb]
  \centering
	  \ra{1.3}
		\begin{tabular}{@{}lll@{}}
		\toprule
		\emph{Fluid} & & \\
		Density & $\rho_f$ & 1000 \kilo\gram\per\meter\cubed \\
		Dynamic viscosity & $\eta$ & 0.001 \pascal\second \\
		\bottomrule
		\emph{Solid} & & \\
		Density & $\rho_s$ & 7850 \kilo\gram\per\meter\cubed \\
		Young's modulus & $E$ & 200 \kilo\pascal \\
		Poisson's ratio & $\nu$ & 0.33 \\
		\bottomrule
		\end{tabular}
\caption{Micro-Beam in Channel Flow: material properties}
\label{table:beamInFlowDeformedProperties}
\end{table}

\begin{figure}[htb]
	\centering
	\subfigure[solids4foam]
	{
		\includegraphics[width=0.48\textwidth]{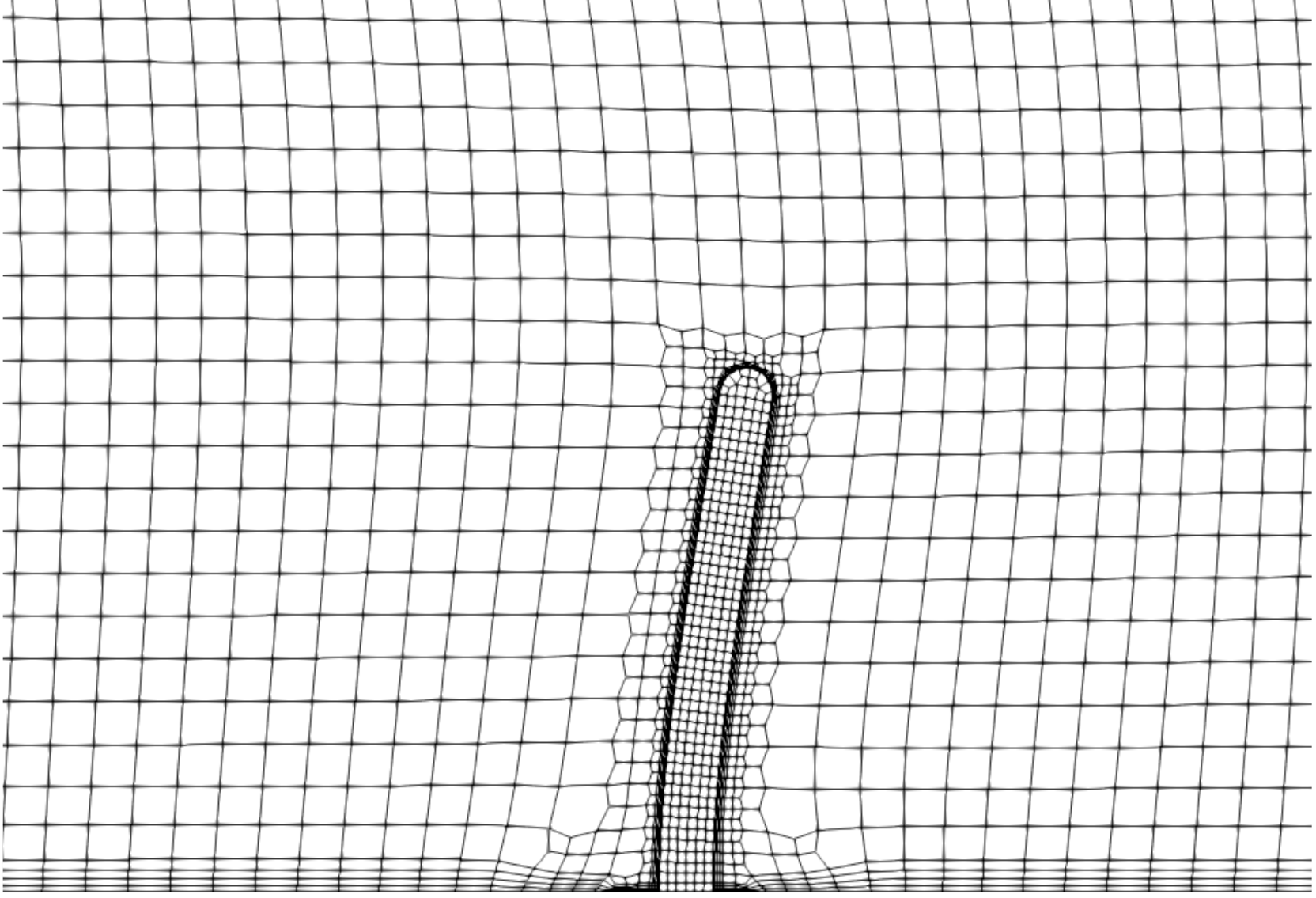}
		\label{fig:beamInFlowDeformedMesh}
	}
	\subfigure[Predictions from Comsol \citep{Comsol:documentation}]
	{
		\includegraphics[width=0.48\textwidth]{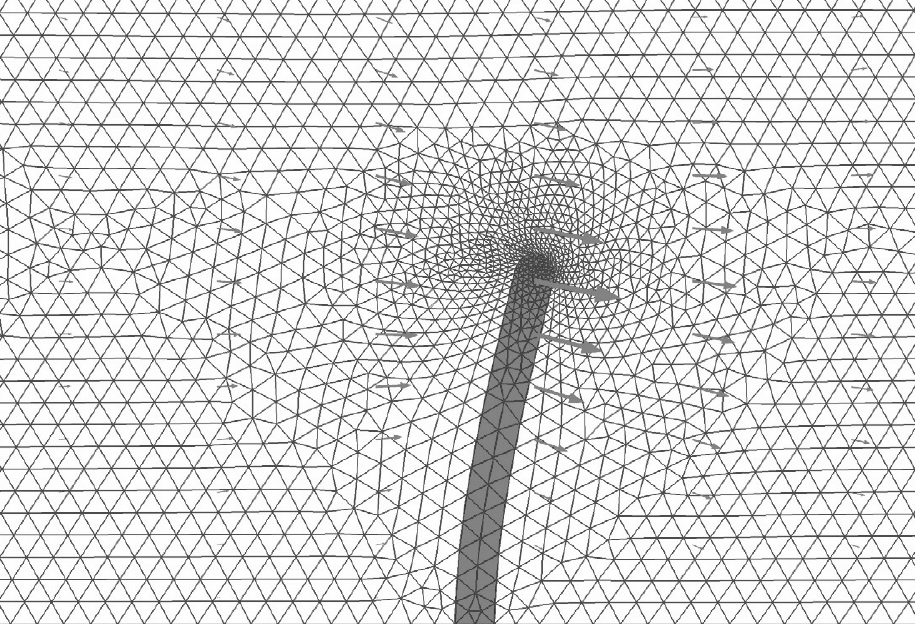}
		\label{fig:beamInFlowDeformedMeshComsol}
	}
	\caption{Micro-Beam in Channel Flow: deformed beam geometry at 0.15 \second}
\end{figure}

\begin{figure}[htb]
	\centering
	\includegraphics[width=0.48\textwidth]{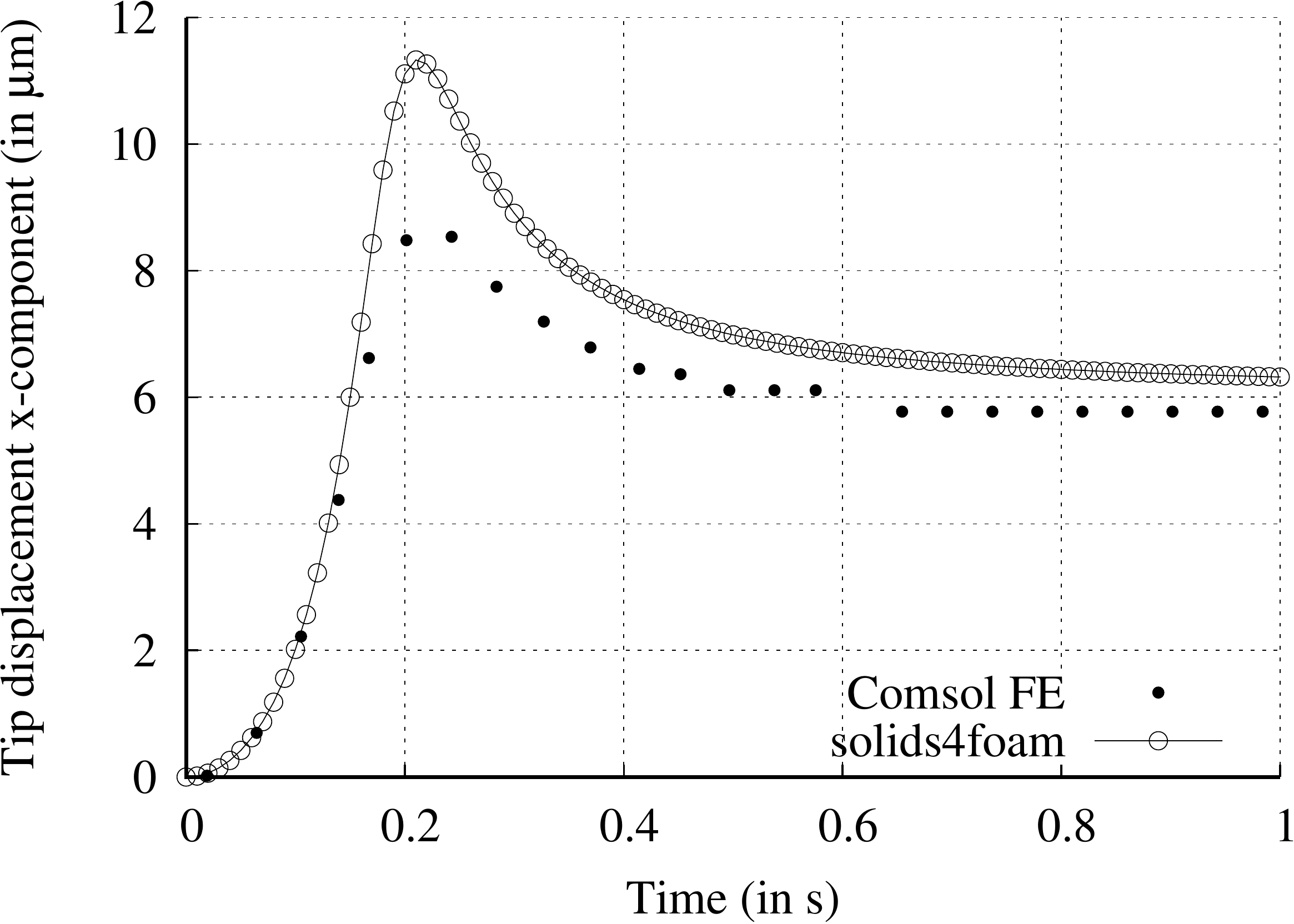}
	\caption{Micro-Beam in Channel Flow: time-varying \emph{x}-component of the beam tip displacement}
	\label{fig:beamInFlowTipDisplacement}
\end{figure}

The predicted deformed geometry of the beam at 0.15 \second\ is shown in Figure \ref{fig:beamInFlowDeformedMeshComsol}, where the predicted geometry from \citet{Comsol:documentation} is given for comparison.
The predicted time-varying \emph{x} component of displacement of the tip of the beam is shown in Figure \ref{fig:beamInFlowTipDisplacement}, with results from \citet{Comsol:documentation} given for comparison.
The displacement versus time trace from the Comsol documentation has been extracted using the {WebPlotDigitizer} software \citep{WebPlotDigitizer} and may introduce some errors.
Given that a relatively coarse mesh has been used in the Comsol model - there are two triangular elements across the beam thickness - the predicted displacements are generally in agreement.
Rather than provide a strict verification benchmark, this test case serves more for demonstration purposes.
Stringent benchmarking of the implemented FSI coupling procedure can be found in \citet{Tukovic2017:fsi}.

\begin{figure}[htb]
	\centering
	\subfigure[solids4foam]
	{
		\includegraphics[width=0.9\textwidth]{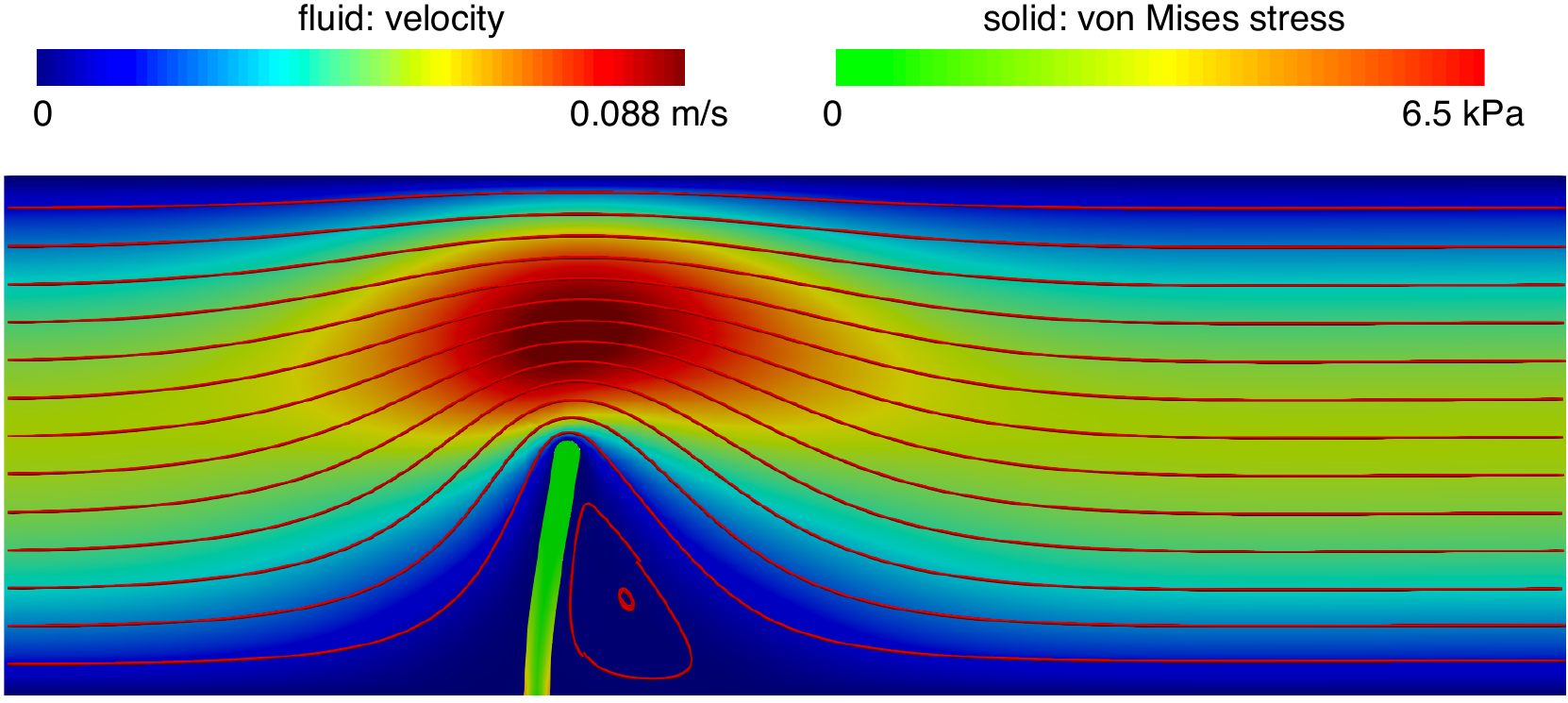}
		\label{fig:beamInFlowVelocity}
	} \\
	\subfigure[Comsol \citep{Comsol:documentation}]
	{
		\includegraphics[width=0.9\textwidth]{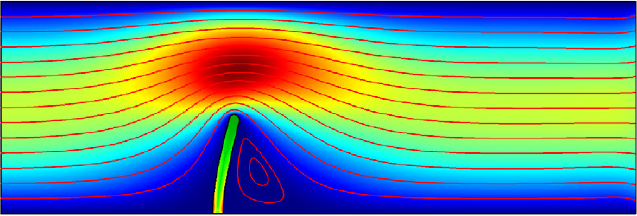}
		\label{fig:beamInFlowVelocityComsol}
	}
	\caption{Micro-Beam in Channel Flow: velocity distribution in the channel at 4.0 \second}
\end{figure}

The velocity distribution in the channel at 4 \second\ is shown in Figure \ref{fig:beamInFlowVelocity}, where once again results from \citet{Comsol:documentation} are given for comparison; contours of equivalent von Mises stress are shown on the solid beam. The results can be seen to be generally in agreement.



\section{Summary, Discussion \& Conclusions}

The current article has presented an open-source toolbox for solid mechanics and FSI simulations based on the FV method, where attention has been given to the design of the use of OOP design principles and code use, comprehension, maintenance, and extension.
The toolbox is built on the open-source library OpenFOAM and takes advantage of parallelisation using the method of domain decomposition.

\paragraph{Remark on discretisation in the context of solid mechanics: \emph{finite volume} vs. \emph{finite element}}
Since the 1970s, the FE method has dominated the field of computational solid mechanics, so much so that the phrase \emph{finite element analysis (FEA)} has become synonymous with computational solid mechanics.
As is evident from the current article alone, alternative methods, such as the finite volume method, have shown similar suitability for tackling a range of CSM problems.
The suitability of the FV method for such analyses is not surprising given the close relationship between FE methods and FV methods.
Over 20 years ago, \citet{Idelsohn1994} noted for 1-D diffusion problems that standard continuous Galerkin FE methods are exactly equivalent to the cell-centred FV method, producing the same matrix coefficients. The exact equivalence no longer holds in 3-D and for discretisation of convection/advection terms, but the close relationship is undeniable.
A key distinguishing feature of typical FV methods is that they produce strongly conservative discretisations, both locally and globally; this is in contrast to typical FE methods that are weakly conservative.
Assuming rigorous mesh sensitivity analyses are performed, the implication of this for CSM problems may be small; however, the development of novel FV discretisation may aid in overcoming shear and volume locking problems, for example, as explored by \citet{Haider2017}.

As the FV method has developed along the lines of an engineering approach that balances fluxes, this leads to a relatively simple formulation and derivation; in contrast, FE methods take a more mathematical approach, building on the framework of variational and weighted residuals.
The FV method applied to CSM comes in many forms, where the presented approach stems from the seminal work of \citet{Demirdzic1988} and co-workers in the late 1980s.
A large number of alternative forms are found in literature, including fully explicit approaches, vertex-centred methods and so-called parametric formulations \citep{Trangenstein1991, Fryer1991, Wheel1996, Slone2002, Pan2010, Kluth2010, Cavalcante2012:generalized, Lee2013, Aguirre2015}.

One indisputable disadvantage of the FV method for solid mechanics, relative to the FE method, is that the field of FV CSM is relatively small and somewhat disjointed.
This does not, however, take away from the potential of FV to not only equal but to outperform FE methods for some CSM applications; in particular, challenging nonlinear multi-physics problems seem to be prime candidates.

\paragraph{Remark on solution methodologies: segregated vs. coupled}
FE methods have typically employed coupled solution methodologies combined with direct linear solvers, where all three displacement components are solved for simultaneously; in contrast, implicit cell-centred FV methods, for example, \citep{Demirdzic1995, Cardiff2016:metalForming}, have used segregated/staggered approaches combined with iterative solvers, where each momentum component equations are temporarily decoupled and solved with outer fixed-point/Picard iterations providing the coupling. Recently, coupled solution approaches have also appeared for FV methods \citep{Das2011, Cardiff2016:blockCoupled}.
An advantage of using a segregated approach and iterative solver, for example, the preconditioned conjugate gradient method, is the relative low computer memory requirements, particularly in comparison to using coupled methods with direct solvers, for example, LU decomposition.
In the early days of CFD, when computer memory was at a premium, such iterative solution methodologies made large complex 3-D simulations possible.
In recent years with the cheap availability of computer memory, this is less of a concern; however, with the movement towards distributed memory computer clusters, iterative solvers once again are becoming attractive with their superior ability to scale in parallel.
When combined with multi-grid methods, both algebraic and geometric, the use of iterative solvers for FV (and FE) simulations has the potential to result in highly efficient parallelised solution procedures.

\paragraph{Remark on parallelisation}
The procedures presented here directly build on the domain decomposition method of the OpenFOAM library, where the model is partitioned into sub-domains, each of which is solved on a separate CPU core; inter-CPU-core communication provides the necessary coupling information to solve the underlying system of algebraic linear equations.
Clearly, efficient parallelisation is predicated on keeping the amount of communication between CPU cores small in comparison to the amount of work done by each CPU core.
This means in practice that it would be expected that  parallel efficiency drops as the number of cells per CPU core decreases and the number of inter-CPU core cell faces (proportional to inter-core communication) increases.
In the presented hot sphere case, this behaviour is evident, where a parallel efficiency of approximately one has been achieved up to 384 CPU cores with 12~500 cells per core, after which the efficiency significantly drops; in fact, the model required more time on 768 CPU cores than on 384.

The parallel efficiency of the presented solvers directly depends on the structure of the OpenFOAM linear equation solvers; however, there are a number of cases where the parallel efficiency is less than that demonstrated here; for example, in contact analyses a neighbour search must be performed to calculate the distances between the master and slave surfaces; this is currently achieved in parallel using so-called \emph{global face zones} where copies of all the faces and points on the contact boundaries are kept on all processors. In addition, these \emph{global face zones} are also currently used when passing information between the solid and fluid domains within FSI analyses.
In addition, the use of dynamic meshes can result in the number of cells per core becoming unbalanced and consequently load balancing is required to preserve parallel efficiency.

\paragraph{Future directions}
The desire for capturing increasingly varied phenomena within multi-physics will lead the future development directions.
In particular, this will include the implementation of important solid formulations not currently available within the toolbox, for example, Eigenproblem procedures for buckling and frequency analyses \citep{Bathe1996}, and beam/shell formulations for thin structures.
As regards the FV method for beams, plates and shells, suitable approaches have been developed by a number of authors \citep{Demirdzic1997:plates, Wheel1997, Fallah2014, Golubovic2017}.
The so-called finite area method, as developed by \citep{Tukovic2008, Tukovic2012:finiteArea}, for example, would provide a suitable framework for implementation.
Separately, the inclusion and exploration of higher-order FV discretisations, for example, approaches recently presented by \citet{Demirdzic2015:fourthOrder}, will also be considered.

In addition to the partitioned FSI approaches presented here, the development of so-called monolithic FSI procedures, where the solid and fluid regions are solved simultaneously, will be considered, and may provide superior stability and efficiency in certain cases. The refactoring of standard fluid models implemented within OpenFOAM will also be performed to link the methods with the \texttt{solids4foam} toolbox, for example, volume-of-fluid multi-phase procedures and block-coupled fluid methods.

\paragraph{Conclusions}
An open-source solid mechanics and FSI toolbox is presented, where emphasis has been given to employing effective coding design paradigms.
The applicability of the FV-based solvers have been demonstrated on a number of cases, and the methods have been shown to provide a comparable approach to standard FE methods.
The toolbox is released as open-source under the GNU Public License \citep{GNUGPL}.

\section{Acknowledgements}
Financial support is gratefully acknowledged from Bekaert through the University Technology Centre (UTC), and from the Irish Centre for Composites Research (IComp). In addition, the authors wish to acknowledge the DJEI/DES/SFI/HEA Irish Centre for High-End Computing (ICHEC) for the provision of computational facilities and support.

\appendix

\section{Appendix: Mechanical Parameters} \label{App:A}

\begin{table}[htb]
  \centering
	\ra{1.3}
	\begin{tabular}{@{}lll@{}}
	\toprule
	& $\lambda > 0$ & first Lam\'{e} parameter  \\
	& $\mu > 0$ & second Lam\'{e} parameter (shear modulus)  \\
	& $\kappa > 0$ & bulk modulus \\
	& $\boldsymbol{C}_e$ & fourth order elastic constitutive tensor \\
	& $g(t)$ & normalised relaxation function \\
	& $\gamma_\infty \geq 0$ & long-term relative modulus  \\
	& $\gamma_i \geq 0$ & relative modulus of the i$^{th}$ Maxwell model spring \\
	& $N > 0$ & bulk modulus \\
	& $\tau_i > 0$ & relaxation time of the i$^{th}$ Maxwell model spring  \\
	& $T_0$ & stress-free reference temperature \\
	& $p$ & pore-pressure \\
	& $\boldsymbol{\epsilon}_e$ & elastic strain \\
	& $\boldsymbol{\epsilon}_p$ & plastic strain \\
	& $\Lambda$ & plastic multiplier \\
	& $f$ & yield function \\
	& $\sigma_Y$ & yield stress/strength \\
	& $\epsilon^{eq}_p$ & equivalent plastic strain \\
	& $g$ & plastic potential \\
	& $\sigma_1 \geq \sigma_2 \geq \sigma_3$ & principal stresses \\
	& $\phi$ & friction angle \\
	& $c$ & cohesion strength \\
	& $\psi$ & dilation angle \\
	\bottomrule
	\end{tabular}
\caption{Mechanical parameters for the linear geometry constitutive laws}
\label{table:MechanicalLawLinearGeometryParameters}
\end{table}

\section{Appendix: \foam Implementation of the Linear Geometry Segregated Algorithm} \label{app:linGeomSolid}
The \foam implementation of the segregated linear geometry \texttt{Solid Model} class is shown in Listing \ref{listing:linGeomSolid}; implicit discretisation operators (finite volume method in OpenFOAM nomenclature) are indicted by \lstinline[language=C++, basicstyle=\ttfamily, breaklines=true]{fvm::} and explicit operators (finite volume calculus in OpenFOAM nomenclature) are indicated by \lstinline[language=C++, basicstyle=\ttfamily, breaklines=true]{fvc::}.
\footnotesize
\begin{lstlisting}[caption={Code excerpt from the segregated linear geometry class: \texttt{linGeomSolid}}, label={listing:linGeomSolid}]
// Momentum equation outer loop: fixed-point/Picard iterations
do
{
    // Store displacement field for under-relaxation and residual calculation
    D.storePrevIter();

    // Discretise the geometrically linear, linear momentum equation in terms
    // of the unknown total displacement field
    fvVectorMatrix DEqn
    (
        rho*fvm::d2dt2(D)
     == fvm::laplacian(impKf, D)
      - fvc::laplacian(impKf, D)
      + fvc::div(sigma)
      + rho*g
      + mechanical().RhieChowCorrection(D, gradD)
    );

    // Apply under-relaxation to the linear system
    DEqn.relax(DEqnRelaxFactor);

    // Solve the linear system using a run-time selectable linear
    // solver, for example, the preconditioned conjugate gradient method
    solverPerfD = DEqn.solve();

    // Apply under-relaxation to the displacement field
    D.relax();

    // Update the gradient of displacement field
    mechanical().grad(D, gradD);

    // Calculate the stress field using a run-time selectable mechanical law
    mechanical().correct(sigma);
}
while (!converged);
\end{lstlisting} \normalsize

\section{References}
\bibliographystyle{unsrtnat}  
\bibliography{Bibliography}

\begin{thebibliography}{112}
\providecommand{\natexlab}[1]{#1}
\providecommand{\url}[1]{\texttt{#1}}
\expandafter\ifx\csname urlstyle\endcsname\relax
  \providecommand{\doi}[1]{doi: #1}\else
  \providecommand{\doi}{doi: \begingroup \urlstyle{rm}\Url}\fi

\bibitem[Weller et~al.(1998)Weller, Tabor, Jasak, and Fureby]{Weller1998}
H.~G. Weller, G.~Tabor, H.~Jasak, and C.~Fureby.
\newblock A tensorial approach to computational continuum mechanics using
  object orientated techniques.
\newblock \emph{Computers in Physics}, 12\penalty0 (6):\penalty0 620--631,
  1998.

\bibitem[Project(2016)]{GNUGPL}
GNU Project.
\newblock \url{http://www.gnu.org/copyleft/gpl.html}, 2016.

\bibitem[Jakirlic et~al.(2017)Jakirlic, Kutej, Unterlechner, and
  Tropea]{Jakirlic2017}
S.~Jakirlic, L.~Kutej, P.~Unterlechner, and C.~Tropea.
\newblock Critical assessment of some popular scale-resolving turbulence models
  for vehicle aerodynamics.
\newblock \emph{{SAE} Int. J. Passeng. Cars - Mech. Syst.}, 10\penalty0
  (1):\penalty0 235--250, 2017.
\newblock \doi{10.4271/2017-01-1532}.

\bibitem[Islam et~al.(2009)Islam, Decker, de~Villiers, Jackson, Gines, Grahs,
  Gitt?Gehrke, and i~Font]{Islam2009}
M.~Islam, F.~Decker, E.~de~Villiers, A.~Jackson, J.~Gines, T.~Grahs,
  A.~Gitt?Gehrke, and J.~Comas i~Font.
\newblock Application of detached?eddy simulation for automotive aerodynamics
  development.
\newblock \emph{{SAE} International}, 2009.
\newblock \doi{10.4271/2009-01-0333}.

\bibitem[Jasak and Weller(2000{\natexlab{a}})]{Jasak2000:linearElasticity}
H.~Jasak and H.~G. Weller.
\newblock Application of the finite volume method and unstructured meshes to
  linear elasticity.
\newblock \emph{International Journal for Numerical Methods in Engineering},
  pages 267--287, 2000{\natexlab{a}}.

\bibitem[Demird\v{z}i\'{c} et~al.(1988)Demird\v{z}i\'{c}, Martinovi\'{c}, and
  Ivankovi\'{c}]{Demirdzic1988}
I.~Demird\v{z}i\'{c}, D.~Martinovi\'{c}, and A.~Ivankovi\'{c}.
\newblock Numerical simulation of thermal deformation in welded workpiece (in
  {C}roatian).
\newblock \emph{Zavarivanje}, 31:\penalty0 209--219, 1988.

\bibitem[Demird\v{z}i\'{c} and Muzaferija(1995)]{Demirdzic1995}
I.~Demird\v{z}i\'{c} and S.~Muzaferija.
\newblock Numerical method for coupled fluid flow, heat transfer and stress
  analysis using unstructured moving meshes with cells of arbitrary topology.
\newblock \emph{Computer Methods in Applied Mechanics and Engineering},
  125\penalty0 (1-4):\penalty0 235--255, 1995.

\bibitem[Demird\v{z}i\'{c} et~al.(1997)Demird\v{z}i\'{c}, Muzaferija, and
  Peri\'{c}]{Demirdzic1997}
I.~Demird\v{z}i\'{c}, S.~Muzaferija, and M.~Peri\'{c}.
\newblock Benchmark solutions of some structural analysis problems using the
  finite-volume method and multgrid acceleration.
\newblock \emph{International Journal for Numerical Methods in Engineering},
  40:\penalty0 1893--1908, 1997.

\bibitem[Tukovi\'{c} and
  Jasak(2007{\natexlab{a}})]{Tukovic2007:updatedLagrangainFAMENA}
{\v{Z}}.~Tukovi\'{c} and H.~Jasak.
\newblock Updated {Lagrangian} finite volume solver for large deformation
  dynamic response of elastic body.
\newblock \emph{Transactions of FAMENA}, 1\penalty0 (31):\penalty0 1--16,
  2007{\natexlab{a}}.

\bibitem[Cardiff et~al.(2014{\natexlab{a}})Cardiff, Kara\v{c}, and
  Ivankovi\'{c}]{Cardiff2014:orthotropicPaper}
P.~Cardiff, A.~Kara\v{c}, and A.~Ivankovi\'{c}.
\newblock A large strain finite volume method for orthotropic bodies with
  general material orientations.
\newblock \emph{Computer Methods in Applied Mechanics and Engineering},
  268:\penalty0 318--335, 2014{\natexlab{a}}.
\newblock \doi{10.1016/j.cma.2013.09.008}.

\bibitem[Cardiff et~al.(2016{\natexlab{a}})Cardiff, Tukovi\'{c}, De~Jaeger,
  Clancy, and Ivankovi\'{c}]{Cardiff2016:metalForming}
P.~Cardiff, Tukovi\'{c}, P.~De~Jaeger, M.~Clancy, and A.~Ivankovi\'{c}.
\newblock A lagrangian cell-centred finite volume method for metal forming
  simulation.
\newblock \emph{International journal for numerical methods in engineering},
  109\penalty0 (13):\penalty0 1777--1803, 2016{\natexlab{a}}.
\newblock \doi{10.1002/nme.5345}.

\bibitem[Leonard et~al.(2012)Leonard, Murphy, Kara\v{c}, and
  Ivankovi\'{c}‡]{Leonard2012}
M.~Leonard, N.~Murphy, A.~Kara\v{c}, and A.~Ivankovi\'{c}‡.
\newblock A numerical investigation of spherical void growth in an
  elastic-"plastic continuum.
\newblock \emph{Computational Materials Science}, 64:\penalty0 38 -- 40, 2012.
\newblock \doi{http://dx.doi.org/10.1016/j.commatsci.2012.04.015}.

\bibitem[Tang et~al.(2015)Tang, Hededal, and Cardiff]{Tang2015}
T.~Tang, O.~Hededal, and P.~Cardiff.
\newblock On finite volume method implementation of poro-elasto-plasticity soil
  model.
\newblock \emph{International Journal for Numerical and Analytical Methods in
  Geomechanics}, 39\penalty0 (13):\penalty0 1410--1430, 2015.

\bibitem[Safari et~al.(2016)Safari, Tukovic, Cardiff, Walter, Casey, and
  Ivankovic]{Safari2016}
A.~Safari, Z.~Tukovic, P.~Cardiff, M.~Walter, E.~Casey, and A.~Ivankovic.
\newblock Interfacial separation of a mature biofilm from a glass surface – a
  combined experimental and cohesive zone modelling approach.
\newblock \emph{Journal of the Mechanical Behavior of Biomedical Materials},
  54:\penalty0 205 -- 218, 2016.
\newblock \doi{http://dx.doi.org/10.1016/j.jmbbm.2015.09.013}.
\newblock URL
  \url{http://www.sciencedirect.com/science/article/pii/S1751616115003495}.

\bibitem[Cardiff et~al.(2012)Cardiff, Kara\v{c}, and
  Ivankovi\'{c}]{Cardiff2012:contactPaper}
P.~Cardiff, A.~Kara\v{c}, and A.~Ivankovi\'{c}.
\newblock Development of a finite volume contact solver based on the penalty
  method.
\newblock \emph{Computational Material Science}, 64:\penalty0 283 -- 284, 2012.

\bibitem[Tukovi\'{c}(2010)]{Tukovic2010acp}
{\v{Z}}.~Tukovi\'{c}.
\newblock Arbitrary crack propagation model in {OpenFOAM}.
\newblock Faculty of Mechanical Engineering and Naval Architecture, University
  of Zagreb, In association with Mechanical \& Materials Engineering,
  University College Dublin, 2010.
\newblock Internal Report.

\bibitem[Carolan et~al.(2013)Carolan, Tukovi\'{c}, Murphy, and
  Ivankovi\'{c}]{Carolan2013}
D.~Carolan, \v{Z}. Tukovi\'{c}, N.~Murphy, and A.~Ivankovi\'{c}.
\newblock Arbitrary crack propagation in multi-phase materials using the finite
  volume method.
\newblock \emph{Computational Materials Science}, 69:\penalty0 153--159, 2013.

\bibitem[Cardiff et~al.(2014{\natexlab{b}})Cardiff, Kara\v{c}, FitzPatrick,
  Flavin, and Ivankovi\'{c}]{Cardiff2014:mappedMusclesPaper}
P.~Cardiff, A.~Kara\v{c}, D.~FitzPatrick, R.~Flavin, and A.~Ivankovi\'{c}.
\newblock Development of mapped stress-field boundary conditions based on a
  {Hill}-type muscle model.
\newblock \emph{International Journal for Numerical Methods in Biomedical
  Engineering}, 2014{\natexlab{b}}.
\newblock \doi{10.1002/cnm}.

\bibitem[Cardiff et~al.(2014{\natexlab{c}})Cardiff, Kara\v{c}, FitzPatrick,
  Flavin, and Ivankovi\'{c}]{Cardiff2014:hipPaper}
P.~Cardiff, A.~Kara\v{c}, D.~FitzPatrick, R.~Flavin, and A.~Ivankovi\'{c}.
\newblock Development of a hip joint model for finite volume simulations.
\newblock \emph{Journal of Biomechanical Engineering}, 136:\penalty0 1--8,
  2014{\natexlab{c}}.
\newblock \doi{10.1115/1.4025776}.

\bibitem[Tukovi\'{c} et~al.(2012)Tukovi\'{c}, Ivankovi\'{c}, and
  Kara\v{c}]{Tukovic2012}
{\v{Z}}.~Tukovi\'{c}, A.~Ivankovi\'{c}, and A.~Kara\v{c}.
\newblock Finite volume stress analysis in multi-material linear elastic body.
\newblock \emph{International Journal for Numerical Methods in Engineering},
  2012.
\newblock \doi{10.1002/nme}.

\bibitem[Cardiff et~al.(2016{\natexlab{b}})Cardiff, Tukovi\'{c}, Jasak, and
  Ivankovi\'{c}]{Cardiff2016:blockCoupled}
P.~Cardiff, Tukovi\'{c}, H.~Jasak, and A.~Ivankovi\'{c}.
\newblock A block-coupled finite volume methodology for linear elasticity and
  unstructured meshes.
\newblock \emph{Computers and Structures}, 175:\penalty0 100--122,
  2016{\natexlab{b}}.
\newblock \doi{10.1016/j.compstruc.2016.07.004}.

\bibitem[Elsafti and Oumeraci(2016)]{Elsafti2016}
H.~Elsafti and H.~Oumeraci.
\newblock A numerical hydro-geotechnical model for marine gravity structures.
\newblock \emph{Computers and Geotechnics}, 79:\penalty0 105 -- 129, 2016.
\newblock \doi{http://dx.doi.org/10.1016/j.compgeo.2016.05.025}.

\bibitem[Haider et~al.(2017)Haider, Lee, Gil, and Bonet]{Haider2017}
J.~Haider, C.~H. Lee, A.~J. Gil, and J.~Bonet.
\newblock A first-order hyperbolic framework for large strain computational
  solid dynamics: An upwind cell centred total lagrangian scheme.
\newblock \emph{International Journal for Numerical Methods in Engineering},
  109\penalty0 (3):\penalty0 407--456, 2017.
\newblock ISSN 1097-0207.
\newblock \doi{10.1002/nme.5293}.
\newblock URL \url{http://dx.doi.org/10.1002/nme.5293}.

\bibitem[Greenshields et~al.(1999)Greenshields, Venizelos, and
  Ivankovi\'{c}]{Greenshields1999}
C.J. Greenshields, H.~Venizelos, and A.~Ivankovi\'{c}.
\newblock The finite volume method for coupled fluid flow and stress analysis.
\newblock \emph{Computer Modeling and Simulation in Engineering}, 4\penalty0
  (3):\penalty0 213 -- 218, 1999.

\bibitem[Ivankovi\'{c} et~al.(2002)Ivankovi\'{c}, Karac, and
  Dendrinos]{Ivankovic2002}
A.~Ivankovi\'{c}, A.~Karac, and E.~Dendrinos.
\newblock Towards early diagnosis of atherosclerosis: the finite volume method
  for fluid-structure interaction.
\newblock \emph{Biorheology}, 39:\penalty0 401 -- 407, 2002.

\bibitem[Kara\v{c} and Ivankovi\'{c}(2002)]{Karac2002}
A.~Kara\v{c} and A.~Ivankovi\'{c}.
\newblock Drop impact of fluid-filled plastic containers: Finite volume method
  for coupled fluid-structure-fracture problems.
\newblock In \emph{Fifth World Congress on Computational Mechanics}, Vienna,
  Austria, 2002.

\bibitem[Tukovi\'{c} and Jasak(2007{\natexlab{b}})]{Tukovic2007:fsi}
{\v{Z}}.~Tukovi\'{c} and H.~Jasak.
\newblock {FVM} for fluid-structure interaction with large structural
  displacements.
\newblock In \emph{2$^{th}$ {OpenFOAM} Workshop}, Zagreb, Croatia,
  2007{\natexlab{b}}.

\bibitem[Kara\v{c} and Ivankovi\'{c}(2009)]{Karac2009}
A.~Kara\v{c} and A.~Ivankovi\'{c}.
\newblock Investigating the behaviour of fluid-filled polyethylene containers
  under base drop impact: A combined experimental/numerical approach.
\newblock \emph{International Journal of Impact Engineering}, 36\penalty0
  (4):\penalty0 621--631, 2009.

\bibitem[Degroote et~al.(2009)Degroote, Bathe, and Vierendeels]{Degroote2009}
J.~Degroote, K.-J. Bathe, and J.~Vierendeels.
\newblock Performance of a new partitioned procedure versus a monolithic
  procedure in fluid--structure interaction.
\newblock \emph{Computers and structures}, 87:\penalty0 793--801, 2009.

\bibitem[Tukovi\'{c} et~al.(2014)Tukovi\'{c}, Cardiff, Ivankovi\'{c}, and
  Kara\v{c}]{Tukovic2014}
{\v{Z}}.~Tukovi\'{c}, P.~Cardiff, A.~Ivankovi\'{c}, and A.~Kara\v{c}.
\newblock Openfoam library for fluid structure interaction.
\newblock In \emph{9$^{th}$ {OpenFOAM} Workshop, Zagreb, Croatia}, Zagreb,
  Croatia, 2014.

\bibitem[Cardiff et~al.(2015)Cardiff, Manchanda, Bryant, Lee, Ivankovi\'{c},
  and Sharma]{Cardiff2015:hydraulicFractures}
P.~Cardiff, R.~Manchanda, E.~C. Bryant, D.~Lee, A.~Ivankovi\'{c}, and M.~M.
  Sharma.
\newblock Simulation of fractures in {OpenFOAM}: From adhesive joints to
  hydraulic fractures.
\newblock In \emph{10$^{th}$ {OpenFOAM} Workshop}, University of Michigan, Ann
  Arbor, MI, USA, 2015.

\bibitem[Gillebaart et~al.(2016)Gillebaart, S., van Zuijlen, and
  Bijl]{Gillebaart2016}
T.~Gillebaart, Blom~D. S., A.~H. van Zuijlen, and H.~Bijl.
\newblock Time consistent fluid structure interaction on collocated grids for
  incompressible flow.
\newblock \emph{Computer Methods in Applied Mechanics and Engineering},
  298\penalty0 (0):\penalty0 159--182, 2016.

\bibitem[\v{S}ekutkovski et~al.(2016)\v{S}ekutkovski, Kosti\'{c},
  Simonovi\'{c}, Cardiff, and Jazarevi\'{c}]{Sekutkovski2016}
B.~\v{S}ekutkovski, I.‡ Kosti\'{c}, A.‡ Simonovi\'{c}, P.~Cardiff, and V.‡
  Jazarevi\'{c}.
\newblock Three-dimensional fluid-structure interaction simulation with a
  hybrid {RANS-€"LES} turbulence model for applications in transonic flow
  domain.
\newblock \emph{Aerospace Science and Technology}, 49:\penalty0 1 -- 16, 2016.

\bibitem[Tukovi\'{c} et~al.(2017)Tukovi\'{c}, Kara\v{c}, Cardiff, Jasak, and
  Ivankovi\'{c}]{Tukovic2017:fsi}
{\v{Z}}.~Tukovi\'{c}, A.~Kara\v{c}, P.~Cardiff, H.~Jasak, and A.~Ivankovi\'{c}.
\newblock Parallel unstructured finite-volume method for fluid-structure
  interaction in {OpenFOAM}.
\newblock \emph{Transactions of FAMENA}, 0:\penalty0 0, 2017.
\newblock Under review.

\bibitem[Project(2017)]{OpenFOAMExtend}
The {FOAM}~Extend Project.
\newblock \url{http://foam-extend.org/}, 2017.

\bibitem[Repository(2017)]{OpenFOAMCommunityRepo}
The {OpenFOAM}~Community Repository.
\newblock \url{https://develop.openfoam.com}, 2017.

\bibitem[Archer(1996)]{Archer1996}
G.~C. Archer.
\newblock \emph{Object-Oriented Finite Element Analysis}.
\newblock PhD thesis, UNIVERSITY of CALIFORNIA at BERKELEY, 1996.

\bibitem[van Riesen et~al.(2004)van Riesen, Monzel, Kaehler, Schlensok, and
  Henneberger]{iMOOSE2004}
D.~van Riesen, C.~Monzel, C.~Kaehler, C.~Schlensok, and Gerhard Henneberger.
\newblock i{MOOSE} - {A}n open-source environment for finite-element
  calculations.
\newblock \emph{IEEE transactions on magnetics}, 40\penalty0 (2, Part
  2):\penalty0 1390--1393, 2004.

\bibitem[Bangerth et~al.(2007)Bangerth, Hartmann, and Kanschat]{Bangerth2007}
W.~Bangerth, R.~Hartmann, and G.~Kanschat.
\newblock Deal.ii\—a general-purpose object-oriented finite element
  library.
\newblock \emph{ACM Trans. Math. Softw.}, 33\penalty0 (4), 2007.
\newblock \doi{10.1145/1268776.1268779}.

\bibitem[Palacios et~al.(2013)Palacios, Alonso, Duraisamy, Colonno, Hicken,
  Aranake, Campos, Copeland, Economon, Lonkar, Lukaczyk, and Taylor]{SU2:2013}
F.~Palacios, J.~Alonso, K.~Duraisamy, M.~Colonno, J.~Hicken, A.~Aranake,
  A.~Campos, S.~Copeland, T.~Economon, A.~Lonkar, T.~Lukaczyk, and T.~Taylor.
\newblock Stanford university unstructured (su2): An open-source integrated
  computational environment for multi-physics simulation and design.
\newblock In \emph{51$6{st}$ {AIAA} Aerospace Sciences Meeting including the
  New Horizons Forum and Aerospace Exposition}, Grapevine (Dallas/Ft. Worth
  Region), Texas, 2013.

\bibitem[Jacobsen et~al.(2012)Jacobsen, Fuhrman, and Freds{\o}e]{Jacobsen2012}
N.~G. Jacobsen, D.~R. Fuhrman, and J.~Freds{\o}e.
\newblock A wave generation toolbox for the open-source {CFD} library:
  {OpenFOAM}.
\newblock \emph{International Journal for Numerical Methods in Fluids},
  70\penalty0 (9):\penalty0 1073--1088, 2012.
\newblock \doi{10.1002/fld.2726}.

\bibitem[Horgue et~al.(2015)Horgue, Soulaine, Franc, Guibert, and
  Debenest]{Horgue2015}
P.~Horgue, C.~Soulaine, J.~Franc, R.~Guibert, and G.~Debenest.
\newblock An open-source toolbox for multiphase flow in porous media.
\newblock \emph{Computer Physics Communications}, 187:\penalty0 217 -- 226,
  2015.
\newblock \doi{http://dx.doi.org/10.1016/j.cpc.2014.10.005}.

\bibitem[Larman(1998)]{Larman1998}
C.~Larman.
\newblock \emph{Applying UML and Patterns to display class structures and
  relationships}.
\newblock Prentice-Hall, Englewood Cliffs, New Jersey, 1998.

\bibitem[Gamma et~al.(1995)Gamma, Helm, Johnson, and Vlissides]{Gamma1995}
E.~Gamma, R.~Helm, R.~Johnson, and J.~Vlissides.
\newblock \emph{Design Patterns: Elements of Reusable Object-oriented
  Software}.
\newblock Addison-Wesley Longman Publishing Co., Inc., Boston, MA, USA, 1995.
\newblock ISBN 0-201-63361-2.

\bibitem[Foundation(2017)]{OpenFOAM:codingStyle}
{OpenFOAM} Foundation.
\newblock Coding style guide.
\newblock \url{http://openfoam.org/dev/coding-style-guide}, 2017.
\newblock Accessed 10$^{th}$ February 2017.

\bibitem[Jasak(1996)]{Jasak1996}
H.~Jasak.
\newblock \emph{Error Analysis and Estimation for the Finite Volume Method with
  Applications to Fluid Flows}.
\newblock PhD thesis, Imperial College London, 1996.

\bibitem[Ferziger and Peric(2002)]{Ferziger2002}
J.~H. Ferziger and M.~Peric.
\newblock \emph{Computational methods for fluid dynamics}.
\newblock Springer, 3$^{rd}$ edition, 2002.

\bibitem[Bathe(1996)]{Bathe1996}
K.~J. Bathe.
\newblock \emph{Finite element procedures}.
\newblock Prentice-Hall, New Jersey, 1996.

\bibitem[Maneeratana(2000)]{Maneeratana2000}
K.~Maneeratana.
\newblock \emph{Development of the finite volume method for non-linear
  structural applications}.
\newblock PhD thesis, Imperial College London, 2000.

\bibitem[Maneeratana and
  Ivankovi\'{c}(1999{\natexlab{a}})]{Maneeratana1999:ACME}
K.~Maneeratana and A.~Ivankovi\'{c}.
\newblock Finite volume method for geometrically nonlinear stress analysis
  applications.
\newblock In \emph{7$^{th}$ Annual {ACME} Conference}, 1999{\natexlab{a}}.

\bibitem[Maneeratana and
  Ivankovi\'{c}(1999{\natexlab{b}})]{Maneeratana1999:ECCM}
K.~Maneeratana and A.~Ivankovi\'{c}.
\newblock Finite volume method for structural applications involving material
  and geometrical non-linearities.
\newblock In \emph{Proceedings of European Conference on Computational
  Mechanics, ECCM'99}, 1999{\natexlab{b}}.

\bibitem[Maneeratana and Ivankovi\'{c}(1999{\natexlab{c}})]{Maneeratana1999}
K.~Maneeratana and A.~Ivankovi\'{c}.
\newblock Finite volume method for large deformation with linear hypoelastic
  materials.
\newblock In \emph{In Finite Volumes for Complex Applications II, Vilsmeier R,
  Benkhaldoun F, Hanel D (eds). HERMES Science Publications}, pages 459--466,
  1999{\natexlab{c}}.

\bibitem[Clifford(2011)]{Clifford2011}
I.~Clifford.
\newblock Block-coupled simulations using {OpenFOAM}.
\newblock In \emph{6$^{th}$ {OpenFOAM} Workshop}, Penn State, USA, 2011.

\bibitem[Das et~al.(2011)Das, Mathur, and Murthy]{Das2011}
S.~Das, S.~R. Mathur, and J.~Y. Murthy.
\newblock An unstructured finite-volume method for structure-electrostatics
  interactions in {MEMS}.
\newblock \emph{Numerical Heat Transfer, Part B: Fundamentals: An International
  Journal of Computation and Methodology}, 60:\penalty0 425--451, 2011.

\bibitem[Jasak and Weller(2000{\natexlab{b}})]{Jasak2000:contact}
H.~Jasak and H.~Weller.
\newblock Finite volume methodology for contact problems of linear elastic
  solids.
\newblock In \emph{Proceedings of 3$^{rd}$ International Conference of Croatian
  Society of Mechanics}, pages 253--260, Cavtat/Dubrovnik, Crotatia,
  2000{\natexlab{b}}.

\bibitem[Ivankovi\'{c} et~al.(1994)Ivankovi\'{c}, Demird\v{z}i\'{c}, Williams,
  and Leevers]{Ivankovic1994}
A.~Ivankovi\'{c}, I.~Demird\v{z}i\'{c}, J.~G. Williams, and P.~S. Leevers.
\newblock Application of the finite volume method to the analysis of dynamic
  fracture problems.
\newblock \emph{International journal of fracture}, 66\penalty0 (4):\penalty0
  357--371, 1994.

\bibitem[Ivankovi\'{c} et~al.(1997)Ivankovi\'{c}, Muzaferija, and
  Demird\v{z}i\'{c}]{Ivankovic1997}
A.~Ivankovi\'{c}, A.~Muzaferija, and I.~Demird\v{z}i\'{c}.
\newblock Finite volume method and multigrid acceleration in modelling of rapid
  crack propagation in full-scale pipe test.
\newblock \emph{Computational mechanics}, 20\penalty0 (1-2):\penalty0 46--52,
  1997.

\bibitem[Ivankovi\'{c} and Venizelos(1998)]{Ivankovic1998}
A.~Ivankovi\'{c} and G.P. Venizelos.
\newblock Rapid crack propagation in plastic pipe: predicting full-scale
  critical pressure from s4 test results.
\newblock \emph{Engineering Fracture Mechanics}, 59\penalty0 (5):\penalty0 607
  -- 622, 1998.

\bibitem[Demird\v{z}i\'{c} et~al.(2000)Demird\v{z}i\'{c}, Horman, and
  Martinovi\'{c}]{Demirdzic2000}
I.~Demird\v{z}i\'{c}, I.~Horman, and D.~Martinovi\'{c}.
\newblock Finite volume analysis of stress and deformation in
  hygro-thermo-elastic orthotropic body.
\newblock \emph{Computer Methods in Applied Mechanics and Engineering},
  190:\penalty0 1221--1232, 2000.

\bibitem[Demird\v{z}i\'{c} et~al.(2005)Demird\v{z}i\'{c}, D\v{z}afarovi\'{c},
  and I.~Ivankovi\'{c}]{Demirdzic2005}
I.~Demird\v{z}i\'{c}, E.~D\v{z}afarovi\'{c}, and A.~I.~Ivankovi\'{c}.
\newblock Finite-volume approach to thermoviscoelasticity.
\newblock \emph{Numerical Heat Transfer, Part B: Fundamentals}, 47\penalty0
  (3):\penalty0 213--237, 2005.

\bibitem[Ba\v{s}i\'{c} et~al.(2005)Ba\v{s}i\'{c}, Demird\v{z}i\'{c}, and
  Muzaferija]{Basic2005}
H.~Ba\v{s}i\'{c}, I.~Demird\v{z}i\'{c}, and S.~Muzaferija.
\newblock Finite volume method for simulation of extrusion processes.
\newblock \emph{International Journal for Numerical Methods in Engineering},
  62\penalty0 (4):\penalty0 475--494, 2005.
\newblock \doi{10.1002/nme.1168}.
\newblock URL \url{http://dx.doi.org/10.1002/nme.1168}.

\bibitem[Bijelonja et~al.(2005)Bijelonja, Demird\v{z}i\'{c}, and
  Muzaferija]{Bijelonja2005}
I.~Bijelonja, I.~Demird\v{z}i\'{c}, and S.~Muzaferija.
\newblock A finite volume method for large strain analysis of incompressible
  hyperelastic materials.
\newblock \emph{International Journal for Numerical Methods in Engineering},
  64\penalty0 (12):\penalty0 1594--1609, 2005.

\bibitem[Bijelonja et~al.(2006)Bijelonja, Demird\v{z}i\'{c}, and
  Muzaferija]{Bijelonja2006}
I.~Bijelonja, I.~Demird\v{z}i\'{c}, and S.~Muzaferija.
\newblock A finite volume method for incompressible linear elasticity.
\newblock \emph{Computer Methods in Applied Mechanics and Engineering},
  195\penalty0 (44-47):\penalty0 6378--6390, 2006.

\bibitem[Demird\v{z}i\'{c}(2015)]{Demirdzic2015:fourthOrder}
I.~Demird\v{z}i\'{c}.
\newblock A fourth-order finite volume method for structural analysis.
\newblock \emph{Applied Mathematical Modelling}, 000:\penalty0 1--11, 2015.

\bibitem[Beale and Elias(1990)]{Beale1990}
S.B. Beale and S.R. Elias.
\newblock Stress distribution in a plate subject to uniaxial loading,
  {PHOENICS}.
\newblock \emph{J. Comput. Fluid Dyn. Appl.}, 3\penalty0 (3):\penalty0
  255--287, 1990.

\bibitem[Trangenstein(1991)]{Trangenstein1991}
J.~A. Trangenstein.
\newblock The riemann problem for longitudinal motion in an elastic-plastic
  bar.
\newblock \emph{Comput. Mech.}, 12:\penalty0 180--207, 1991.

\bibitem[Spalding(1993)]{Spalding1993}
D.B. Spalding.
\newblock Simulation of fluid flow, heat transfer and solid deformation
  simultaneously.
\newblock In \emph{NAFEMS Conference No. 4}, Brighton, UK, 1993.

\bibitem[Hattel and Hansen(1995)]{Hattel1995}
J.H. Hattel and P.N. Hansen.
\newblock A control volume-based finite difference method for solving the
  equilibrium equations in terms of displacements.
\newblock \emph{Appl. Math. Model.}, 19:\penalty0 210--243., 1995.

\bibitem[Trangenstein(1994)]{Trangenstein1994}
J.~A. Trangenstein.
\newblock Second-order godunov algorithm for two-dimensional solid mechanics.
\newblock \emph{Comput. Mech.}, 13:\penalty0 343--359, 1994.

\bibitem[Fryer et~al.(1991)Fryer, Bailey, Cross, and Lai]{Fryer1991}
Y.~D. Fryer, C.~Bailey, M.~Cross, and C.~H. Lai.
\newblock A control volume procedure for solving elastic stress-strain
  equations on an unstructured mesh.
\newblock \emph{Applied mathematical modelling}, 15\penalty0 (11-12):\penalty0
  639--645, 1991.

\bibitem[Wheel(1996)]{Wheel1996}
M.~A. Wheel.
\newblock A geometrically versatile finite volume formulation for plane
  elastostatic stress analysis.
\newblock \emph{The Journal of Strain Analysis for Engineering Design},
  31\penalty0 (2):\penalty0 111--116, 1996.

\bibitem[Slone et~al.(2002)Slone, Pericleous, Bailey, and Cross]{Slone2002}
A.K. Slone, K.~Pericleous, C.~Bailey, and M.~Cross.
\newblock Dynamic fluid structure interaction using finite volume unstructured
  mesh procedures.
\newblock \emph{Computers and Structures}, 80\penalty0 (5–6):\penalty0 371 --
  390, 2002.
\newblock ISSN 0045-7949.
\newblock \doi{http://dx.doi.org/10.1016/S0045-7949(01)00177-8}.

\bibitem[Pan et~al.(2010)Pan, Wheel, and Qin]{Pan2010}
W.~Pan, M.A. Wheel, and Y.~Qin.
\newblock Six-node triangle finite volume method for solids with a rotational
  degree of freedom for incompressible material.
\newblock \emph{Computers and Structures}, 88\penalty0 (23-"24):\penalty0
  1506--1511, 2010.
\newblock ISSN 0045-7949.
\newblock \doi{http://dx.doi.org/10.1016/j.compstruc.2010.08.001}.
\newblock Special Issue: Association of Computational Mechanics, United
  Kingdom.

\bibitem[Cavalcante and Pindera(2012)]{Cavalcante2012:generalized}
M.~A.~A. Cavalcante and M.-J. Pindera.
\newblock Generalized finite-volume theory for elastic stress analysis in solid
  mechanics--part {I}: Framework.
\newblock \emph{Journal of Applied Mechanics}, 79\penalty0 (5), 2012.
\newblock \doi{10.1115/1.4006805}.

\bibitem[Nordbotten(2014)]{Nordbotten2014}
J.~M. Nordbotten.
\newblock Cell-centered finite volume discretizations for deformable porous
  media.
\newblock \emph{International journal for numerical methods in engineering},
  2014.
\newblock \doi{10.1002/nme.4734}.

\bibitem[Simo and Hughes(1998)]{Simo1998}
J.~C. Simo and T.~J.~R. Hughes.
\newblock \emph{Computational Inelasticity}, volume~7.
\newblock Springer-Verlag, New York, 1998.

\bibitem[Belytschko et~al.(2000)Belytschko, Liu, and Moran]{Belytschko2000}
T.~Belytschko, W.~Kam Liu, and B.~Moran.
\newblock \emph{Nonlinear Finite Elements for Continua and Structures}.
\newblock Wiley, 2000.

\bibitem[Jacobs(1980)]{Jacobs1980}
D.~A.~H. Jacobs.
\newblock Preconditioned conjugate gradient methods for solving systems of
  algebraic equations.
\newblock \emph{Central Electricity Research Laboratories Report},
  RD/L/N193/80, 1980.

\bibitem[Guennebaud et~al.(2010)Guennebaud, Jacob, et~al.]{eigenweb}
Ga\"{e}l Guennebaud, Beno\^{i}t Jacob, et~al.
\newblock Eigen v3.
\newblock http://eigen.tuxfamily.org, 2010.

\bibitem[Amestoy et~al.(2001)Amestoy, Duff, Koster, and L'Excellent]{MUMPS}
P.~R. Amestoy, I.~S. Duff, J.~Koster, and J.-Y. L'Excellent.
\newblock A fully asynchronous multifrontal solver using distributed dynamic
  scheduling.
\newblock \emph{SIAM Journal on Matrix Analysis and Applications}, 23\penalty0
  (1):\penalty0 15--41, 2001.

\bibitem[Balay et~al.(2016)Balay, Abhyankar, Adams, Brown, Brune, Buschelman,
  Dalcin, Eijkhout, Gropp, Kaushik, Knepley, McInnes, Rupp, Smith, Zampini,
  Zhang, and Zhang]{petsc-web-page}
S.~Balay, S.~Abhyankar, M.~F. Adams, J.~Brown, P.~Brune, K.~Buschelman,
  L.~Dalcin, V.~Eijkhout, W.~D. Gropp, D.~Kaushik, M.~G. Knepley, L.~C.
  McInnes, K.~Rupp, B.~F. Smith, S.~Zampini, H.~Zhang, and H.~Zhang.
\newblock {PETS}c {W}eb page.
\newblock \url{http://www.mcs.anl.gov/petsc}, 2016.
\newblock URL \url{http://www.mcs.anl.gov/petsc}.

\bibitem[Heroux et~al.(2005)Heroux, Bartlett, Howle, Hoekstra, Hu, Kolda,
  Lehoucq, Long, Pawlowski, Phipps, Salinger, Thornquist, Tuminaro,
  Willenbring, Williams, and Stanley]{Trilinos}
Michael~A. Heroux, Roscoe~A. Bartlett, Vicki~E. Howle, Robert~J. Hoekstra,
  Jonathan~J. Hu, Tamara~G. Kolda, Richard~B. Lehoucq, Kevin~R. Long, Roger~P.
  Pawlowski, Eric~T. Phipps, Andrew~G. Salinger, Heidi~K. Thornquist, Ray~S.
  Tuminaro, James~M. Willenbring, Alan Williams, and Kendall~S. Stanley.
\newblock An overview of the trilinos project.
\newblock \emph{ACM Trans. Math. Softw.}, 31\penalty0 (3):\penalty0 397--423,
  2005.
\newblock ISSN 0098-3500.
\newblock \doi{http://doi.acm.org/10.1145/1089014.1089021}.

\bibitem[Muzaferija(1994)]{Muzaferija1994}
S.~Muzaferija.
\newblock \emph{Adaptive Finite Volume Method Flow Prediction Using
  Unstructured Meshes and Multigrid Approach}.
\newblock British Thesis Service. University of London, 1994.

\bibitem[Fainberg and Leister(1996)]{Fainberg1996}
J.~Fainberg and H.~J. Leister.
\newblock Finite volume multigrid solver for thermo-elastic analysis in
  anisotropic materials.
\newblock \emph{Comput. Methods. Appl. Mech. Engrg.}, 137:\penalty0 167--174,
  1996.

\bibitem[{Rhie} and {Chow}(1983)]{Rhie1983}
C.~M. {Rhie} and W.~L. {Chow}.
\newblock {Numerical study of the turbulent flow past an airfoil with trailing
  edge separation}.
\newblock \emph{AIAA Journal}, 21:\penalty0 1525--1532, November 1983.
\newblock \doi{10.2514/3.8284}.

\bibitem[Lee et~al.(2015)Lee, Cardiff, Bryant, Manchanda, Wang, and
  Sharma]{Lee2015}
D.~Lee, P.~Cardiff, E.~C. Bryant, R.~Manchanda, H.~Wang, and M.~M. Sharma.
\newblock A new model for hydraulic fracture growth in unconsolidated sands
  with plasticity and leak-off.
\newblock In \emph{{SPE} Annual Technical Conference and Exhibition, 28-30
  September}, 2015.
\newblock \doi{10.2118/174818-MS}.

\bibitem[{SimScale}(2017)]{SimScaleDoc:sphericalPresureVessel}
{SimScale}.
\newblock Validation cases: Design analysis of spherical pressure vessel.
\newblock MIT OpenCourseWare, 2017.

\bibitem[Afkar et~al.(2014)Afkar, Camari, and Paykani]{Afkar2014}
A.~Afkar, M.~N. Camari, and A.~Paykani.
\newblock Design and analysis of a spherical pressure vessel using finite
  element method.
\newblock \emph{World Journal of Modelling and Simulation}, 10\penalty0
  (2):\penalty0 126--135, 2014.

\bibitem[Foundation(2014)]{OpenFOAM.org}
The~{OpenFOAM} Foundation.
\newblock \url{http://www.openfoam.org}, 2014.

\bibitem[Roache(1997)]{Roache1997}
P.~J. Roache.
\newblock Quantification of uncertainty in computational fluid dynamics.
\newblock \emph{Annual review of fluid mechanics}, 29:\penalty0 123--60, 1997.

\bibitem[Karypis and Kumar(1999)]{Karypis1999}
G.~Karypis and V.~Kumar.
\newblock A fast and highly quality multilevel scheme for partitioning
  irregular graphs.
\newblock \emph{{SIAM} Journal on Scientific Computing}, 20\penalty0
  (1):\penalty0 359--392, 1999.

\bibitem[for Finite Element~Methods and (U.K.)(1990)]{NAFEMS}
National~Agency for Finite Element~Methods and Standards (U.K.).
\newblock \emph{The Standard NAFEMS Benchmarks}.
\newblock NAFEMS, 1990.

\bibitem[for Finite Element~Methods and (U.K.)(2006)]{NAFEMS:contact}
National~Agency for Finite Element~Methods and Standards (U.K.).
\newblock \emph{NAFEMS Advanced Finite Element Contact Benchmarks, Benchmark 2,
  NAFEMS Publication Ref: R0094}.
\newblock NAFEMS, 2006.

\bibitem[Code\_Aster(2017)]{CodeAster:documentation}
Code\_Aster.
\newblock \url{http://www.code-aster.org}, 2017.
\newblock Accessed April 2017.

\bibitem[FreeCAD(2017)]{FreeCAD}
FreeCAD.
\newblock Freecad documentation.
\newblock \url{http://www.freecadweb.org}, 2017.

\bibitem[Geuzaine and Remacle(2009)]{Geuzaine2009}
C.~Geuzaine and J.-F. Remacle.
\newblock Gmsh: A 3-d finite element mesh generator with built-in pre- and
  post-processing facilities.
\newblock \emph{International Journal for Numerical Methods in Engineering},
  79\penalty0 (11):\penalty0 1309--1331, 2009.
\newblock \doi{10.1002/nme.2579}.

\bibitem[Jureti\'{c} et~al.(2015)Jureti\'{c}, Baburi\'{c}, and
  Lugari\'{c}]{cfmeshweb}
F.~Jureti\'{c}, M.~Baburi\'{c}, and T.~Lugari\'{c}.
\newblock cfmesh v1.1.
\newblock http://www.c-fields.com, 2015.

\bibitem[WebPlotDigitizer(2017)]{WebPlotDigitizer}
WebPlotDigitizer.
\newblock Web based tool to extract data from plots, images, and maps.
\newblock \url{http://arohatgi.info/WebPlotDigitizer}, 2017.

\bibitem[Comsol(2017)]{Comsol:documentation}
Comsol.
\newblock Comsol multi-physics 5.2a documentation.
\newblock \url{http://www.comsol.com}, 2017.

\bibitem[Zienkiewicz and Taylor(2000)]{Zienkiewicz2000}
O.~C. Zienkiewicz and R.~L. Taylor.
\newblock \emph{The finite element method: volume 2 solid mechanics}.
\newblock Butterworth-Heineman, 5$^{th}$ edition, 2000.

\bibitem[Turek and Hron(2006)]{Turek2006}
S.~Turek and J.~Hron.
\newblock \emph{Proposal for Numerical Benchmarking of Fluid-Structure
  Interaction between an Elastic Object and Laminar Incompressible Flow}.
\newblock Springer Berlin Heidelberg, Berlin, Heidelberg, 2006.
\newblock \doi{http://dx.doi.org/10.1007/3-540-34596-5\_15}.

\bibitem[Jureti\'{c}(2016)]{cfMesh}
F.~Jureti\'{c}.
\newblock {cfMesh} meshing library, {Creative Fields}.
\newblock \url{https://cfmesh.com}, 2016.
\newblock Accessed 7$^{th}$ July 2017.

\bibitem[Idelsohn and O\~{n}ate(1994)]{Idelsohn1994}
S.~R. Idelsohn and E.~O\~{n}ate.
\newblock Finite volumes and finite elements: Two 'good friends'.
\newblock \emph{International Journal for Numerical Methods in Engineering},
  37\penalty0 (19), 1994.
\newblock \doi{10.1002/nme.1620371908}.

\bibitem[Kluth and Despr\'{e}s(2010)]{Kluth2010}
G.~Kluth and B.~Despr\'{e}s.
\newblock Discretization of hyperelasticity on unstructured mesh with a
  cell-centered {Lagrangian} scheme.
\newblock \emph{Journal of Computational Physics}, 229\penalty0 (24):\penalty0
  9092 -- 9118, 2010.
\newblock \doi{http://dx.doi.org/10.1016/j.jcp.2010.08.024}.

\bibitem[Lee et~al.(2013)Lee, J., and J.]{Lee2013}
C.~H. Lee, Gil~A. J., and Bonet J.
\newblock Development of a cell centred upwind finite volume algorithm for a
  new conservation law formulation in structural dynamics.
\newblock \emph{Comput. Struct.}, 118:\penalty0 13--38, 2013.

\bibitem[Aguirre et~al.(2015)Aguirre, Gil, Bonet, and Lee]{Aguirre2015}
M.~Aguirre, A.~J. Gil, J.~Bonet, and C.~H. Lee.
\newblock An upwind vertex centred finite volume solver for {Lagrangian} solid
  dynamics.
\newblock \emph{J. Comput. Phys.}, 300:\penalty0 387--422, 2015.

\bibitem[Demird\v{z}i\'{c} and Ivankovi\'{c}(1997)]{Demirdzic1997:plates}
I.~Demird\v{z}i\'{c} and A.~Ivankovi\'{c}.
\newblock Finite volume approach to modelling of plates.
\newblock In \emph{Proceedings of the 2$^{nd}$ Congress of Croatian Society of
  Mechanics}, pages 101--108, Brac, Croatia, 1997.

\bibitem[Wheel(1997)]{Wheel1997}
M.~A. Wheel.
\newblock A finite volume method for analysing the bending deformation of thick
  and thin plates.
\newblock \emph{Computer Methods in Applied Mechanics and Engineering},
  147\penalty0 (1-2):\penalty0 199--208, 1997.

\bibitem[Fallah and Ebrahimnejad(2014)]{Fallah2014}
N.~Fallah and M.~Ebrahimnejad.
\newblock Finite volume analysis of adaptive beams with piezoelectric sensors
  and actuators.
\newblock \emph{Applied Mathematical Modelling}, 38\penalty0 (2):\penalty0 722
  -- 737, 2014.
\newblock ISSN 0307-904X.
\newblock \doi{http://dx.doi.org/10.1016/j.apm.2013.07.004}.
\newblock URL
  \url{http://www.sciencedirect.com/science/article/pii/S0307904X13004484}.

\bibitem[Golubovi\'{c} et~al.(2017)Golubovi\'{c}, Demird\v{z}i\'{c}, and
  Muzaferija]{Golubovic2017}
A.~Golubovi\'{c}, I.~Demird\v{z}i\'{c}, and S.~Muzaferija.
\newblock Finite volume analysis of laminated composite plates.
\newblock \emph{International Journal for Numerical Methods in Engineering},
  109\penalty0 (11):\penalty0 1607--1620, 2017.
\newblock ISSN 1097-0207.
\newblock \doi{10.1002/nme.5347}.
\newblock URL \url{http://dx.doi.org/10.1002/nme.5347}.

\bibitem[Tukovi\'{c} and Jasak(2008)]{Tukovic2008}
{\v{Z}}.~Tukovi\'{c} and H.~Jasak.
\newblock Simulation of free-rising bubble with soluble surfactant using moving
  mesh finite volume / area method.
\newblock In \emph{6th International Conference on CFD in Oil and Gas,
  Metallurgical and Process Industries}, pages 1--11, Trondheim, Norway, 2008.

\bibitem[Tukovi\'{c} and Jasak({2012})]{Tukovic2012:finiteArea}
{\v{Z}}.~Tukovi\'{c} and H.~Jasak.
\newblock {A moving mesh finite volume interface tracking method for surface
  tension dominated interfacial fluid flow}.
\newblock \emph{{Computers \& Fluids}}, {55}:\penalty0 {70--84}, {2012}.
\newblock \doi{{10.1016/j.compfluid.2011.11.003}}.

\end{thebibliography}

\end{document}